\documentclass[a4paper,11pt,leqno]{amsart}
\usepackage[latin1]{inputenc}
\usepackage{amsmath} 
\usepackage{mathtools} 
\usepackage{amssymb}
\usepackage{enumerate} 
\usepackage{xspace} 
\usepackage{amsfonts}
\usepackage{latexsym} 
\usepackage{amsthm} 
\usepackage[backrefs]{amsrefs}
\usepackage[all]{xypic} 
\usepackage{hyperref}

\numberwithin{equation}{section}

\newtheorem{lemma}{Lemma}[section]
\newtheorem{proposition}[lemma]{Proposition}
\newtheorem{theorem}[lemma]{Theorem}
\newtheorem{corollary}[lemma]{Corollary}
\newtheorem{claim}{Claim}
\newtheorem*{clai}{Claim}

\theoremstyle{definition}
\newtheorem{definition}[lemma]{Definition}
\newtheorem{example}[lemma]{Example}

\newtheorem{remark}[lemma]{Remark}

\def\C{\mathbb{C}}
\def\N{\mathbb{N}}
\def\No{{\N_0}}
\def\Z{\mathbb{Z}}

\def\O{\mathcal{O}}

\def\X{\hat{X}}
\def\H{\mathsf{H}}
\def\A{\ba{X}}
\def\Aa{\mathcal{A}}
\def\DX{\mathcal{D}_{\OSS}}

\DeclareMathOperator{\proj}{proj}
\DeclareMathOperator{\spa}{span}
\DeclareMathOperator{\Id}{Id}

\newcommand{\pxa}{\bigl(\XA,\theta_\XA,\fg{\ind}\bigr)}
\newcommand{\cspxa}{\cs\bigl(\XA,\theta_\XA,\fg\ind\bigr)}
\newcommand{\xanp}{\bigl(\XA,\theta_\XA,\fg{n}\bigr)}
\newcommand{\K}{\mathbf{K}}
\newcommand{\indd}{\N^2_0}
\newcommand{\lee}{\preceq}
\newcommand{\gee}{\succeq}
\newcommand{\BO}{\Bx(\O_A)}
\newcommand{\Bx}{\mathcal{B}}
\newcommand{\B}{\Bx(X,\theta,G)}

\newcommand{\cy}[2]{C(#1,#2)}
\newcommand{\spc}{\overline{\spa}}
\newcommand{\coni}{(*)}

\newcommand{\cs}{C^*}
\newcommand{\inv}{^{-1}}
\newcommand{\ba}[1]{\bar{#1}}
\newcommand{\h}[1]{\hat{#1}}
\newcommand{\fg}[1]{\mathbb{F}_{#1}}
\newcommand{\al}{\mathfrak a}
\newcommand{\Past}{{\mathcal P}}
\newcommand{\equ}[2]{\hspace{2pt}_{#1}\hspace{-3pt}\sim_{#2}}
\newcommand{\eq}{\equ{k}{l}}
\newcommand{\OSS}[1][]{{\mathsf{X}^{+}_{#1}}}
\newcommand{\TSS}{\mathsf{X}}

\newcommand{\laOSS}{\mathcal{L}(\OSS)}
\newcommand{\ind}{\mathcal{I}}
\newcommand{\ec}[3]{{}_{#1}\hspace{-1pt}[#2]_{#3}}
\newcommand{\osh}{\sigma}
\newcommand{\tsh}{\tau}
\newcommand{\msh}{\sigma}
\newcommand{\alwords}{\al^*}
\newcommand{\XA}{{\OSS[A]}}

\newcommand{\Dx}{\mathcal{D}_{\OSS}}
\newcommand{\BA}{\Bx(\XA,\theta_\XA,\fg\ind)}

\newcommand{\cset}{\mathcal{F}}
\newcommand{\boomap}[1][(X,\theta,G)]{\phi_{#1}}

\newcommand{\OSD}{\mathrm{X^+}}
\newcommand{\osd}{\sigma}
\newcommand{\TSD}{\mathrm{X}}
\newcommand{\tsd}{\tau}
\newcommand{\csp}[1]{\cs(#1)}

\newcommand{\BXY}{\Bx(X\setminus Y,\theta_{|X\setminus Y},G)}
\newcommand{\crosmap}{\hspace{.3ex}\eta_\rtimes}
\newcommand{\tcrosmap}{\tilde{\eta}_\rtimes}
\newcommand{\fw}{\mathcal{W}}
\newcommand{\find}{\mathcal{F}}
\newcommand{\ideal}{\mathcal{I}}
\newcommand{\uset}[1][\ideal]{U(#1)}
\newcommand{\ioss}{\eta_\O}
\newcommand{\pnes}{\mathcal{E}}
\newcommand{\spro}{\pi}
\newcommand{\boof}{\eta}
\newcommand{\nTSS}{n_\TSS}
\newcommand{\emptyword}{\varepsilon}
\newcommand{\rte}{\mathcal{J}_\TSS}
\newcommand{\jj}{\mathbf{j}}
\newcommand{\xx}{\mathbf{x}}

\renewcommand{\mu}{u}
\renewcommand{\nu}{v}

\providecommand{\norm}[1]{[#1]} 
\providecommand{\abs}[1]{\lvert#1\rvert}

\begin{document}
\title[Dynamics, Boolean algebras and $\cs$-algebras]{Symbolic
  dynamics, partial dynamical systems, Boolean algebras and
  $\cs$-algebras generated by partial isometries}
\author{Toke Meier Carlsen}
\date{\today}
\address{Mathematisches Institut\\
Einsteinstra\ss e 62\\
48149 M\" unster\\
Germany}
\thanks{This research has been supported by the EU-Network Quantum
  Spaces - Noncommutative Geometry (HPRN-CT-2002-00280).} 
\curraddr{Department of Mathematics\\
  University of Newcastle\\
  NSW 2308\\
  Australia}
\email{toke@math.uni-muenster.de} 

\begin{abstract}
  We associate to each discrete partial dynamical system a universal
  $\cs$-algebra generated by partial isometries satisfying relations
  given by a Boolean algebra connected to the discrete partial
  dynamical system in question. We show that for symbolic dynamical
  systems like one-sided and two-sided shift spaces and topological
  Markov chains with an arbitrary state space the $\cs$-algebras
  usually associated to them can be obtained in this way.

  As a consequence of this, we will be able to show that for two-sided
  shift spaces having a certain property, the crossed product of the
  two-sided shift space is a quotient of the $\cs$-algebra associated
  to the corresponding one-sided shift space.
 \end{abstract}

\maketitle

\section{Introduction}

The history of associating $\cs$-algebras to symbolic dynamical
systems is long and successful.  

A good example of this is the crossed product of infinite minimal
two-sided shift spaces which in \cite{MR1363826} was used to classify infinite
minimal shift spaces up till strong orbit equivalence and flip
conjugacy (it is actually done for a bigger class of dynamical
systems, namely Cantor systems, but we will in this paper only concern
ourselves with symbolic dynamical systems). 

Another very important example of a class of $\cs$-algebras associated
to symbolic dynamical systems is the Cuntz-Krieger algebras
\cite{MR561974}, which in a natural way can be viewed as
$\cs$-algebras associated to topological Markov chains with finite
state space. The Cuntz-Krieger algebras have proved to be very
important examples in the theory of $\cs$-algebras and have also let
to invariants of shift of finite type such as the dimension group
(cf. \cite{MR561973} and \cite{MR561974}). 

The Cuntz-Krieger algebras have been generalized in many different
ways. We will in this paper focus on two of those. The first one is
due to Exel and Lace, who in \cite{MR2000i:46064} have generalized the
Cuntz-Krieger algebras to topological Markov chains with arbitrary
state space. The other generalization is due to Matsumoto,
who in \cite{MR1454478} associated a $\cs$-algebra to every shift
space (called a subshift in that paper). Matsumoto's original
construction associated a $\cs$-algebra to every \emph{two-sided}
shift space, but it is more natural to view it as a way to associated
a $\cs$-algebra to every \emph{one-sided} shift space (cf. \cite{tmc}
and \cite{CS}).  

Topological Markov chains with finite state space are examples of
one-sided shift spaces, and it turns out that the $\cs$-algebras
associated to these kind of shift spaces are Cuntz-Krieger algebras,
so in this way the class of $\cs$-algebras associated with shift
spaces is a generalization of the class of Cuntz-Krieger algebras
(cf. \cite{CS}*{Section 8}).  

Thus we have three different classes of $\cs$-algebras associated to
symbolic dynamical system, namely crossed products of two-sided shift
spaces, Exel and Laca's generalization of Cuntz-Krieger algebras and
$\cs$-algebras associated to one-sided shift spaces. The main purpose
of this paper is to unify these three constructions to one.
 
My original motivation for written this paper was to prove Theorem
\ref{theorem:toke}, which for shift spaces having a certain property
relates the crossed product of the two-sided shift space and the
$\cs$-algebra associated to the corresponding one-sided shift
space. Doing this I found that the crossed product of two-sided shift
spaces and the $\cs$-algebra associated to one-sided shift spaces have
a common structure, which I also found in Exel and Laca's
generalization of Cuntz-Krieger algebras (which we for simplicity from
now on just will call Cuntz-Krieger algebras). This structure can be
described by partial representations of groups and Boolean
algebras. More formal, what we will do is to associate to every so
called discrete partial dynamical system a $\cs$-algebra and then to
every symbolic dynamical system (both one-sided and two-sided)
associate a discrete partial dynamical system in such a way that the
$\cs$-algebras we get in this way for one-sided shift spaces,
two-sided shift spaces and topological Markov chains are canonical
isomorphic to the $\cs$-algebra associated to the one-sided shift
space, the crossed product of the two-sided shift space and the
unitization of the Cuntz-Krieger of the transition matrix of the
topological Markov chain, respectively.

This construction is very natural, and I hope that beside the
benefits from having a unified construction of these different classes
of $\cs$-algebras associated to different kinds of symbolic dynamical
systems, the construction will also clarify in which way the
structure of the symbolic dynamical system is reflected in the
structure of the associated $\cs$-algebra.

The paper is organized as follows: In Section \ref{sec:notat}, we will
shortly introduce some notation which will be used throughout the
paper, in Section \ref{ekse} \emph{discrete partial dynamical
  systems} will be defined, and we will see how we from one- and
two-sided symbolic dynamical systems can construct discrete partial
dynamical systems. We will then in Section \ref{sec:cs-algebra} define
the $\cs$-algebra of a discrete partial dynamical system and show some
basic properties of it. In Section \ref{sec:construction} we will show
the $\cs$-algebra of a discrete partial dynamical system can be
constructed as a crossed product of a $\cs$-partial dynamical systems,
and we will in Section \ref{sec:repr} construct a representation of
the $\cs$-algebra of a discrete partial dynamical system as operators
on a Hilbert space. Section \ref{retu} is the main section of this
paper; here we will show that the crossed product of a two-sided shift
space, the $\cs$-algebra of a one-sided shift space, and Cuntz-Krieger
algebras can be obtained as $\cs$-algebras of discrete partial
dynamical systems. We will then in Section \ref{sec:ideal} describe
the ideal structure of the $\cs$-algebra of a discrete partial
dynamical system and then use this description to prove the above
mentioned Theorem \ref{theorem:toke}. The paper finish with three
appendices in which we will give a short introduction to partial
representations of groups, Boolean algebras and crossed products of
$\cs$-partial dynamical systems.

A previous version of this paper appeared in my Ph.D thesis
\cite{phd}. Unfortunately that version contained a lot of mistakes, which
hopefully have been fixed in this version. The most notable of these
mistakes was that I claimed that the $\cs$-algebra of a higher rank graph
could be constructed as the $\cs$-algebra of a discrete partial
dynamical system. The proof of this is however false, but it is
possible to construct the $\cs$-algebras of a higher rank graph
in a very similar way by using an action of an (discrete) inverse
semigroup instead of a partial action of a (discrete) group. This will
be proved in a forthcoming paper by Gwion Evans and the author.

\subsection*{Acknowledgment}
The process of writing this paper has been very long. I started on it
will I was a Ph.D-student at the University of Copenhagen, I continued
working on it while I was a post Doc at the Institut Mittag-Leffler and
at the Norwegian University of Science and Technology and I finally
finished it as a post Doc at Universit\" at M\" unster. I wish to
thank all members of the operator algebra groups at these places
for their kind hospitality and especially S\o ren Eilers, Christian Skau
and Joachim Cuntz. I also wish to thank Aidan Sims for pointing out
the above mentioned mistake. 

\section{Notation and preliminaries} \label{sec:notat}
Throughout this paper, $e$ will denote the \emph{neutral element} of a
given group. We will by $\Z$ denote \emph{the set of integers}, $\N_0$
will denote \emph{the set of non-negative integers} and $\N$ will
denote \emph{the set of positive integers}. 

If $X$ is a set, then we will by $\Id_X$ denote \emph{the identity
map on $X$}. For a map $\sigma:X\to X$, we will for every $k\in\N$ by
$\sigma^k$ denote \emph{the $k$-times composition of $\sigma$ with itself},
and we will set $\sigma^0=\Id_X$. If $\sigma$ is invertible, then we
will for every $k\in\N$ by $\sigma^{-k}$ denote the map
$(\sigma\inv)^k$. If $A$ is a subset of $X$, then $1_A$ denotes \emph{the
characteristic function}:
\begin{equation*}
  1_A(x)=
  \begin{cases}
    1&\text{if }x\in A,\\
    0&\text{if }x\notin A.
  \end{cases}
\end{equation*}
If $\theta$ is a map defined on $A$, then we will for another subset
$B$ of $X$ by $\theta(B)$ mean $\theta(A\cap B)$, and we will by
$\theta_{|B}$ denote the restriction of $\theta$ to $A\cap B$. 

If $C$ is a subset of a vector space, then we will use $\spa(C)$ to
denote \emph{the linear span of $C$}, and if the vector space also comes with a
topology, then $\spc(C)$ will denote \emph{the closure of the linear
  span of $C$}. For a subset $D$ of a $\cs$-algebra we denote
\emph{the $\cs$-subalgebra generated by $D$}, by $\cs(D)$.

When $\al$ is a set, then we will by $\fg\al$ denote \emph{the free group
generated by $\al$}. We will regard $\al$ as a subset of $\fg\al$ and
denote the subset $\bigl\{a\inv\mid a\in\al\bigr\}$ of $\fg\al$ by
$\al\inv$. We say that an element $g\in\fg\al$ \emph{is written
in reduced form $b_1b_2\dotsm b_k$} if $b_1,b_2,\dotsc,b_k\in\al\cup\al\inv$ and
$b_j=b\Rightarrow b_{j+1}\ne b\inv$ for every $j\in\{1,2,\dotsc,k-1\}$
and every $b\in\al\cup\al\inv$. We let $[\cdot]$ be the unique
homomorphism from $\fg\al$ to $\Z$ such that $[a]=1$ for every $a\in\al$. 

We will by $\al^*$ denote \emph{the set of finite words of symbols from
$\al$}. If $\mu\in \al^*$, then we will by $\abs{\mu}$ denote \emph{the
length of $\mu$} (i.e., the number of symbols in $\mu$), by $\mu_1$ the
first letter (the leftmost) letter of $\mu$, by $\mu_2$ the second
letter of $\mu$, and so on till $\mu_{|\mu|}$ which denotes the last
(the rightmost) letter of $\mu$. Thus $\mu=\mu_1\mu_2\dotsm
\mu_{|\mu|}$. We will identify $\al^*$ with the positive cone of $\fg\al$ 
which is defined to be the unital sub-semigroup of $\fg\al$ generated by
$\al$. Under this identification, the unit element of $\fg\al$ is
equal to the \emph{empty word}, which we will denote by $\emptyword$,
and $\abs{\mu}$ is equal to $\norm{\mu}$ for all $\mu\in\al^*$.

We will denote the set of one-sided, respectively two-sided, infinite
sequences in $\al$ by $\al^\No$, respectively $\al^\Z$. We will
often denote an element $x=(x_n)_{n\in\N}$ of $\al^\No$ by
\begin{equation*}
  x_0x_1\dotsm ,
\end{equation*}
and if $\mu\in\al^*$, then we will by $\mu x$ denote the sequence
\begin{equation*}
  \mu_1\mu_2\dotsm \mu_{|\mu|}x_0x_1\dotsm .
\end{equation*}
We will also often for a sequence $x$ belonging to either $\al^\No$
or $\al^\Z$ and for integers $k< l$ belonging to the appropriate
index set denote $x_kx_{k+1}\dotsm x_{l-1}$ by $x_{[k,l[}$ and regard
it as an element of $\al^*$. Likewise will $x_{[k,\infty[}$ denote
the element
\begin{equation*}
  x_kx_{k+1}\dotsm
\end{equation*}
of $\al^\N$.

We say that an element $(x_n)_{n\in\Z}\in\al^\Z$ is \emph{periodic} if
there exists an $m\in\N$ such that $x_{n+m}=x_n$ for all $n\in\Z$, and
we say that an element $(x_n)_{n\in\No}\in\al^{\No}$ is \emph{eventually
  periodic} if there exist $m,N\in\N$ such that $x_{n+m}=x_n$ for $n>N$.

\section{Discrete partial dynamical systems} \label{ekse}
We are now going to define \emph{partial actions} and \emph{discrete
  partial dynamical systems} and look at some ways to construct
discrete partial dynamical systems. 

Partial actions have been defined and studied by Ruy Exel in
\cite{MR1469405}, where he for any given group constructed an inverse
semigroup such that there is a one-to-one correspondence between the actions of
the inverse semigroup and the partial actions of the group.

\begin{definition}
  Given a group $G$ and a set $X$, a \emph{partial action} $\theta$ of
  $G$ on $X$ is a pair
  \begin{equation*}
    \bigl((D_g)_{g\in G},(\theta_g)_{g\in G}\bigr),
  \end{equation*}
  where for each $g\in G$, $D_g$ is a subset of $X$ and $\theta_g$ is
  a bijective map from $D_{g\inv}$ to $D_g$, satisfying for all $h$
  and $i$ in $G$:
  \begin{subequations}    
    \begin{align}
      &D_e=X \text{ and } \theta_e \text{ is the identity map on }X,\\
      &\theta_h(D_i)=D_h\cap D_{hi}, \label{eq:paractb}\\
      &\theta_h(\theta_i(x))=\theta_{hi}(x) \text{ for }x\in
      D_{i\inv}\cap D_{i\inv h\inv}.
    \end{align}
  \end{subequations}
  We will call the family $(D_g)_{g\in G}$ the \emph{domains} of
  $\theta$ and the family $(\theta_g)_{g\in G}$ the \emph{partial
    one-to-one maps} of $\theta$.

  The triple $(X,\theta,G)$ is called a \emph{discrete partial
    dynamical system}.
\end{definition}

The reason that we have chosen to call such a triple a \emph{discrete}
partial dynamical system is that we do not require any structure of
the system $(X,\theta,G)$ other than the above mentioned. It is in many
cases (see for example \cite{MR2003f:46108}) natural to ask for the set
$X$ to be a topological space, the domains $(D_g)_{g\in G}$ open
subsets, and the maps $(\theta_g)_{g\in G}$ homeomorphisms, but we
will in this paper only consider the discrete case, where we do not
require such a structure. 

A very simple example of a discrete partial dynamical systems is if we for a
given group $G$ and a given set $X$ for every $g\in G$ let $D_g=X$ and
let $\theta_g$ be the identity map on $X$. We get slightly more
interesting examples if we consider group actions:
\begin{example} \label{group}
  Let $G$ be a group, $X$ a set and $\theta$ an action of $G$ on $X$,
  i.e., $\theta_g$ is for every $g\in G$ a map from $X$ to $X$ such that
  \begin{subequations}    
    \begin{align}
      &\theta_e \text{ is the identity map on } X,\\
      &\theta_h\circ\theta_i=\theta_{hi} \text{ for every } h,i\in G.
    \end{align}
  \end{subequations}
  If we for every $g\in G$ let $D_g=X$, then $\theta=\bigl((D_g)_{g\in
    G},(\theta_g)_{g\in G}\bigr)$ is a partial action of $G$ on $X$, and
  $(X,\theta,G)$ a discrete partial dynamical system.
\end{example}

As mentioned in the Introduction, we will in this paper mainly concern
ourselves with partial dynamical systems which come from symbolic
dynamical systems. 
We will now show how to get partial dynamical systems from symbolic
dynamical systems; we will first see how to define a partial dynamical
system from a \emph{one-sided} symbolic dynamical system, and
then how to define a partial dynamical system from a \emph{two-sided}
symbolic dynamical system.  

\subsection{One-sided symbolic dynamical systems}\label{onesided}
Let $(\OSD,\osd)$ be a one-sided symbolic dynamical
system over the alphabet $\al$. That is: $\al$ is a set (finite or
infinite), $\osd:\al^{\No}\to \al^{\No}$ is the map 
\begin{equation*}
  x_0x_1x_2\dotsm\mapsto x_1x_2\dotsm
\end{equation*}
and $\OSD$ is a subset of $\al^{\No}$ such
that $\osd(\OSD)\subseteq \OSD$. To turn $(\OSD,\osd)$ into a
partial dynamical system we restrict $\osd$ to subsets of $\OSD$
such that $\osd$ is injective on these subsets. This is done in this way:

Let for every $a\in \al$, $D_a$ be the subset $\bigl\{(x_n)_{n\in\No}\in
\OSD\mid x_0=a\bigr\}$, $D_{a\inv}$ be the subset $\osd(D_a)$,
$\theta_a:D_{a\inv}\to D_a$ be the map 
\begin{equation*}
  x\mapsto ax,
\end{equation*}
and $\theta_{a\inv}:D_a\to D_{a\inv}$ be the map
\begin{equation*}
  x\mapsto \osd(x).
\end{equation*}
Let $\fg{\al}$ be the free group generated by $\al$ and let for every
$g\in \fg{\al}$ written in the reduced form $b_1b_2\dotsm b_k$, where
$b_1,b_2,\dotsc ,b_k\in \al\cup\al\inv$, $D_g$ be the subset of $\OSD$
defined by 
\begin{equation*}
  D_g=\theta_{b_1}\circ\theta_{b_2}\circ \dotsb \circ\theta_{b_k}(\OSD),
\end{equation*}
and $\theta_g:D_{g\inv}\to D_g$ be the map defined by
\begin{equation*}
  \theta_g=\theta_{b_1}\circ\theta_{b_2}\circ \dotsb \circ\theta_{b_k}.
\end{equation*}
Then $\theta_\OSD=\bigl((D_g)_{g\in\fg\al},(\theta_g)_{g\in\fg\al}\bigr)$ is a
partial action of $\fg\al$ on $\OSD$, and $(\OSD,\theta_\OSD,\fg{\al})$
is a discrete partial dynamical system.

\begin{definition} \label{onesideddef}
  Let $(\OSD,\osd)$ be a one-sided symbolic dynamical system over the
  alphabet $\al$. Then we
  call the discrete partial dynamical system
  $(\OSD,\theta_\OSD,\fg{\al})$ constructed above \emph{the
  discrete partial dynamical system associated to $(\OSD,\osd)$}.
\end{definition}

\begin{lemma} \label{lemma:onefund}
  Let $(\OSD,\osd)$ be a one-sided symbolic dynamical system over the
  alphabet $\al$ and let $(\OSD,\theta_\OSD,\fg{\al})$ be the discrete
  partial dynamical system associated to $(\OSD,\osd)$ as in
  Definition \ref{onesideddef}. Then the following holds for the
  partial action $\theta_\OSD=\bigl((D_g)_{g\in G},(\theta_g)_{g\in G}\bigr)$:
  \begin{enumerate}
  \item if $g\in G$ and $D_g\ne\emptyset$, then there exist
    $\mu,\nu\in\al^*$ such that $g=\mu\nu\inv$, \label{item:1}
  \item if $\mu,\nu\in\al^*$ and the last (rightmost) letters of $\mu$
    and $\nu$ are not equal (or either $\mu$ or $\nu$ is equal to the
    empty word), then we have that 
    \begin{equation*}
      D_{\nu\mu\inv}=\Bigl\{x\in\OSD\bigm| x_{[0,|\nu|[}=\nu,\ \mu
      x_{[|\nu|,\infty[}\in\OSD\Bigr\},
    \end{equation*}
    and that $\theta_{\nu\mu\inv}$ is the map 
    \begin{equation*}
      \mu x\mapsto \nu x
    \end{equation*}
    from $D_{\mu\nu\inv}$ to $D_{\nu\mu\inv}$, \label{item:2}
  \item if $\mu,\nu\in\al^*$, $|\mu|=|\nu|$ and $\mu\ne\nu$, then
    $D_\mu\cap D_\nu=\emptyset$. \label{item:3}
  \end{enumerate}
\end{lemma}

\begin{proof}
  (\ref{item:2}) can easily be proved by induction over the length of
  $\mu$ and $\nu$. (\ref{item:3}) then follows from (\ref{item:2}),
  and (\ref{item:1}) follows from the definition of $D_g$ and (\ref{item:3}).
\end{proof}

\subsection{Two-sided symbolic dynamical systems} \label{twosided}
Let $(\TSD,\tsd)$ be a two-sided symbolic dynamical system over the
alphabet $\al$. That is: $\al$ is a set (finite or infinite),
$\tsd:\al^{\Z}\to \al^{\Z}$ is the map defined by 
\begin{equation*}
  \bigl(\tsd((z_n)_{n\in\Z})\bigr)_m=z_{m+1}
\end{equation*}
for every $(z_n)_{n\in\Z}\in\al^{\Z}$ and every $m\in\Z$, and $\TSD$ is
a subset of $\al^{\Z}$ such that $\tsd(\TSD)=\TSD$. 

Since $\tsd$ is bijective, $(\tsd^k)_{k\in\Z}$ is an action of $\Z$ on
$\TSD$, so we could turn $(\TSD,\tsd)$ into a discrete partial
dynamical system by the method of Example \ref{group}, but we can also
do it by imitating the method we used to define the discrete partial
dynamical system associated to a one-sided symbolic dynamical system,
and that is what we will do here. 

Let for every $a\in \al$, $D_a$ be the subset
$\bigl\{(z_n)_{n\in\Z}\in \TSD\mid z_0=a\bigr\}$, $D_{a\inv}$ be the
subset $\bigl\{(z_n)_{n\in\Z}\in \TSD\mid z_{-1}=a\bigr\}$,
$\theta_a:D_{a\inv}\to D_a$ be the restriction of $\tsd\inv$ to
$D_{a\inv}$, and $\theta_{a\inv}:D_a\to D_{a\inv}$ be the restriction
of $\tsd$ to $D_a$. 
  
Let $\fg{\al}$ be the free group generated by $\al$, and let for every
$g\in \fg{\al}$ written in the reduced form $b_1b_2\dotsm b_k$, where
$b_1,b_2,\dotsc ,b_k\in \al\cup\al\inv$, $D_g$ be the subset of
$\TSD$ defined by
\begin{equation*}
  D_g=\theta_{b_1}\circ\theta_{b_2}\circ \dotsb \circ\theta_{b_n}(\TSD),
\end{equation*}
and let $\theta_g:D_{g\inv}\to D_g$ be the map defined by
\begin{equation*}
  \theta_g=\theta_{b_1}\circ\theta_{b_2}\circ \dotsb \circ\theta_{b_n}.
\end{equation*} 
Then $\theta_\TSD=\bigl((D_g)_{g\in\fg\al},(\theta_g)_{g\in\fg\al}\bigr)$
is a partial action of $\fg\al$ on $\TSD$, and $(\TSD,\theta_\TSD,\fg{\al})$
is a discrete partial dynamical system. 

\begin{definition} \label{twosideddef}
  Let $(\TSD,\tsd)$ be a two-sided symbolic dynamical system over the
  alphabet $\al$. Then we 
  call the discrete partial dynamical system
  $(\TSD,\theta_\TSD,\fg{\al})$ constructed above \emph{the
  discrete partial dynamical system associated to $(\TSD,\tsd)$}.
\end{definition}

\begin{lemma} \label{lemma:twofund}
  Let $(\TSD,\tsd)$ be a two-sided symbolic dynamical system over the
  alphabet $\al$ and let $(\TSD,\theta_\TSD,\fg{\al})$ be the discrete
  partial dynamical system associated to $(\TSD,\tsd)$ as in
  Definition \ref{twosideddef}. Then the following holds for the
  partial action $\theta_\TSD=\bigl((D_g)_{g\in G},(\theta_g)_{g\in G}\bigr)$:
  \begin{enumerate}
  \item if $g\in G$ and $D_g\ne\emptyset$, then either $g\in\al^*$ or
    $g\inv\in\al^*$, \label{item:t1}
  \item if $\mu\in\al^*$, then we have that
    \begin{equation*}
      D_\mu=\bigl\{z\in\TSD\mid z_{[0,|\mu|[}=u\bigr\},\quad
      D_{\mu\inv}=\bigl\{z\in\TSD\mid z_{[-|\mu|,0[}=u\bigr\}, 
    \end{equation*}
    $\theta_{\mu}$ is the restriction of $\tsd^{-|\mu|}$ to
    $D_{\mu\inv}$ and $\theta_{\mu\inv}$ is the restriction of
    $\tsd^{|\mu|}$ to $D_\mu$, \label{item:t2}
  \item if $\mu,\nu\in\al^*$, $|\mu|=|\nu|$ and $\mu\ne\nu$, then
    $D_\mu\cap D_\nu=\emptyset$ and $D_{\mu\inv}\cap
    D_{\nu\inv}=\emptyset$. \label{item:t3}
  \end{enumerate}
\end{lemma}

\begin{proof}
  (\ref{item:t2}) can easily be proved by induction over the length of
  $\mu$. (\ref{item:t3}) then follows from (\ref{item:t2}),
  and (\ref{item:t1}) follows from the definition of $D_g$ and (\ref{item:t3}).
\end{proof}

\section{The $\cs$-algebra associated to a discrete partial dynamical
  system}  \label{sec:cs-algebra}
The main object of this paper is to associate to every
discrete partial dynamical system $(X,\theta,G)$ a $\cs$-algebra
$\csp{X,\theta,G}$ in such a way that the class of
$\cs$-algebras obtained in this way in a natural way generalizes
Cuntz-Krieger algebras (both for finite and infinite matrices),
$\cs$-algebras associated to one-sided shift spaces and  
crossed products of two-sided shift spaces.

We want the $\cs$-algebra $\csp{X,\theta,G}$ associated to a discrete
partial dynamical system $(X,\theta,G)$ to have the
property that there is a bijective correspondence between the
representations of $\csp{X,\theta,G}$ and certain representations of
$(X,\theta,G)$. 

Let us as a motivating example look at the Cuntz-Krieger algebra $\O_A$ of an
$n\times n$-matrix $A=(A(i,j))_{i,j=1}^n$ with entries in $\{0,1\}$
and without any zero rows. The matrix $A$ gives raise to a one-sided
symbolic dynamical system, namely the topological Markov chain with
transition matrix $A$: 
\begin{equation*}
  \XA=\Big\{(x_n)_{n\in \No}\in\{1,2,,\dotsc,n\}^{\No}\bigm| \forall
  i\in\No:A(x_i,x_{i+1})=1\Big\}.
\end{equation*}
Thus we have by Definition \ref{onesideddef} a partial action $\theta_\XA$ of
the free group $\fg{n}$ generated by $n$ generators on $\XA$. It turns
out the structure of the discrete partial 
dynamical system $\xanp$ is reflected in the structure
of $\O_A$ in the following way: 

Let $(S_i)_{i=1}^n$ denote the generators of $\O_A$. Then the map 
\begin{equation*}
  i\mapsto S_i
\end{equation*}
extends to a partial representation $u$ of $\fg{n}$ on $\O_A$ (see Appendix
\ref{sec:parrep} for a short introduction to partial representations).
Since $u$ is a partial representation, $\bigl(u(g)\bigr)_{g\in\fg{n}}$ is a
family of partial isometries with commuting range projections, so
$\cs\bigl(\{u(g)u(g)^*\mid g\in \fg{n}\}\bigr)$ is a unital
abelian $\cs$-algebra, and so the set of projections in this
$\cs$-algebra is in a natural way a Boolean algebra (see Appendix
\ref{sec:boo} for a definition of this Boolean algebra). 

The discrete partial dynamical system $\xanp$ also gives raise to a Boolean
algebra: namely the Boolean algebra generated by 
all the domains $\{D_g\mid g\in \fg{n}\}$ of the partial action
$\theta_\XA$ (cf. Appendix \ref{sec:boo}). 
It turns out that the map
\begin{equation*}
  D_g\mapsto u(g)u(g)^*
\end{equation*}
extends to a Boolean homomorphism of the Boolean algebra generated by
the domains $\{D_g\mid g\in \fg{n}\}$ of the partial action $\theta_\XA$ to
the Boolean algebra of projections in the $\cs$-subalgebra of $\O_A$
generated by $\{u(g)u(g)^*\mid g\in \fg{n}\}$.
It is natural to view this Boolean homomorphism together with the
partial representation $u$ of $\fg{n}$ as a representation of the
partial dynamical system $\xanp$. We will in Section \ref{retu} see
that this representation of the partial dynamical system $\xanp$
on $\O_A$ is universal and thus completely characterizes $\O_A$.

Similar characterizations hold for Cuntz-Krieger algebras of infinite
matrices, the $\cs$-algebras associated to one-sided shift spaces and
the crossed products of two-sided shift spaces.

Thus we are lead to the following definition of
$\csp{X,\theta,G}$:
\begin{definition} \label{def:ox}
  Let $(X,\theta,G)$ be a discrete partial dynamical system. Then
  $\csp{X,\theta,G}$ is the universal $\cs$-algebra generated by
  a family $(s_g)_{g\in G}$ of elements satisfying:
  \begin{subequations}    
    \begin{align}
       &\parbox[t]{0.8\textwidth}{$(s_g)_{g\in
           G}$ is a family of partial isometries with commuting range
         projections,} \label{eq:ox}\\
      &s_e=1,\label{eq:oxa}\\
      &s_{g\inv}=s_g^* \text{ for every } g\in G,\label{eq:oxb}\\
      &s_hs_i=s_hs_h^*s_{hi} \text{ for every } h,i\in G,\label{eq:oxc}\\
      &\parbox[t]{0.8\textwidth}{the map $D_g\mapsto s_gs_g^*$ extends to a
        Boolean homomorphism from the Boolean algebra $\B$ generated
        by the domains $\{D_g\mid g\in G\}$ of the partial action
        $\theta$ to the Boolean algebra of projections in the unital
        abelian $\cs$-algebra $\cs\bigl(\{s_g
          s_g^*\mid g\in G\}\bigr)$.} \label{eq:oxd}   
    \end{align}
  \end{subequations}
  We will call the family $(s_g)_{g\in G}$ the \emph{generators} of
  $\csp{X,\theta,G}$. 
\end{definition} 

\begin{remark}
It is not immediately clear that a $\cs$-algebra with the
properties mentioned in Definition \ref{def:ox} exists,
but we will in Section \ref{sec:construction} for every  discrete
partial dynamical system constructed such a $\cs$-algebra as a
crossed product of a $\cs$-partial dynamical system.
\end{remark}

\begin{remark}
  It follows from Appendix \ref{sec:parrep} that condition
  \eqref{eq:ox}-\eqref{eq:oxc} are 
  equivalent to the following 3 conditions: 
  \begin{subequations}    
    \begin{align}
      &s_es_e^*=1,\\
      &s_g^*s_g=s_{g\inv}s_{g\inv}^* \text{ for every } g\in G,\\
      &s_hs_is_i^*s_h^*=s_hs_is_{hi}^* \text{ for every } h,i\in G
    \end{align}
  \end{subequations}
  (cf. \cite{MR99a:46121}) and with the following 4 conditions:
  \begin{subequations}    
    \begin{align}
      &\parbox[t]{0.8\textwidth}{$(s_g)_{g\in G}$
        is a family of partial isometries with commuting range
        projections,}\\
      &s_e=1,\\
      &s_{g\inv}=s_g^* \text{ for every }g\in G,\\
      &s_hs_is_{i\inv}=s_{hi}s_{i\inv} \text{ for every }h,i\in G,
    \end{align}
  \end{subequations}
  (cf. \cite{MR2003f:46108}).
\end{remark}

\begin{definition} \label{remarkboomap}
  Let $(X,\theta,G)$ be a discrete partial dynamical system, 
  $(D_g)_{g\in G}$ the domains of $\theta$ and $(s_g)_{g\in G}$ the
  generators of $\csp{X,\theta,G}$. It is clear that the Boolean
  homomorphism from $\B$ to the Boolean algebra of projections in
  $\cs\bigl(\{s_gs_g^*\in\csp{X,\theta,G}\mid g\in G\}\bigr)$ which
  extends the map $D_g\mapsto s_gs_g^*$ is unique. We will denote it
  by $\boomap$.  
\end{definition}
We will in Corollary \ref{boomapinj} see that $\boomap$ is always
injective, and thus that $\csp{X,\theta,G}\ne 0$ (unless $X=\emptyset$).

The following lemma describe how the structure of the discrete partial
dynamical system $(X,\theta,G)$ is reflected in the structure of the
$\cs$-algebra $\csp{X,\theta,G}$.

\begin{lemma} \label{booinvariant}
  If $(X,\theta,G)$ is a discrete partial dynamical system, then we
  have that 
  \begin{equation*}
     \boomap\bigl(\theta_g(A)\bigr)=s_g\boomap(A)s_g^*
  \end{equation*}
  for every $A\in\B$ and every $g\in G$.
\end{lemma}

\begin{proof}
  If $\boomap\bigl(\theta_g(A)\bigr)=s_g\boomap(A)s_g^*$ and
  $\boomap\bigl(\theta_g(B)\bigr)=s_g\boomap(B)s_g^*$, then the
  following series of equalities hold:
  \begin{equation*}
    \begin{split}
      \boomap\bigl(\theta_g(A\cap B)\bigr)&=
      \boomap\bigl(\theta_g(A)\cap\theta_g(B)\bigr)\\
      &= \boomap\bigl(\theta_g(A)\bigr)\boomap\bigl(\theta_g(B)\bigr)\\
      &= s_g\boomap(A)s_g^*s_g\boomap(B)s_g^*\\
      &= s_gs_g^*s_g\boomap(A)\boomap(B)s_g^*\\
      &= s_g\boomap\bigl(A\cap B\bigr)s_g^*,
    \end{split}
  \end{equation*}
  and so does the next:
  \begin{equation*}
    \begin{split}
      \boomap\bigl(\theta_g(X\setminus A)\bigr) &=
      \boomap\bigl(D_g\setminus\theta_g(A)\bigr)\\
      &= s_gs_g^*-s_g\boomap(A)s_g^*\\
      &= s_g\bigl(1-\boomap(A)\bigr)s_g^*\\
      &= s_g\boomap(X\setminus A)s_g^*.
    \end{split}
  \end{equation*}
  So it is enough to show that
  $\boomap\bigl(\theta_g(D_h)\bigr)=s_gs_hs_h^*s_g^*$ for every $g,h\in G$.

  So let $g,h\in G$. Then we have that
  \begin{equation*}
    \begin{split}
      \boomap\bigl(\theta_g(D_h)\bigr)
      &= \boomap\bigl(D_g\cap D_{gh}\bigr)\\
      &= s_gs_g^*s_{gh}s_{gh}^*\\
      &= s_gs_g^*s_{gh}s_{gh}^*s_gs_g^*\\
      &= s_gs_hs_h^*s_g^*,\\
    \end{split}
  \end{equation*}
  which proves the lemma.
\end{proof}

The next lemma, which we state for further reference, is an easy but
useful consequent of the definition of the $\cs$-algebra of a discrete
partial dynamical system.

\begin{lemma} \label{lemma:orto}
  Let $(X,\theta,G)$ be a discrete partial dynamical system and let
  $g,h,i\in G$. If $D_g=\emptyset$, then $s_g=0$, and if $D_h\cap
  D_i=\emptyset$, then $s_h^*s_i=0$.
\end{lemma}

\begin{proof}
  If $D_g=\emptyset$, then $s_g=s_gs_g^*s_g=\boomap\bigl(D_g\bigr)s_g=0$,
  and if $D_h\cap D_i=\emptyset$, then we have that
  \begin{equation*}
    \begin{split}
      s_h^*s_i&=s_h^*s_hs_h^*s_is_i^*s_i\\
      &=s_h^*\boomap\bigl(D_h\bigr)\boomap\bigl(D_i\bigr)s_i\\
      &=s_h^*\boomap\bigl(D_h\cap D_i\bigr)s_i\\
      &=0.
    \end{split}
  \end{equation*}
\end{proof}

\begin{example}
  If $\theta$ is an action of a discrete group $G$ on a set $X\ne
  \emptyset$ and 
  $(X,\theta,G)$ is the discrete partial dynamical system defined in
  Example \ref{group}, then all the domains $(D_g)_{g\in G}$ of the
  partial action $\theta$ are equal to $X$. Thus the Boolean algebra
  generated by these domains only consists of $X$ and $\emptyset$, and so
  $\cs(X,\theta,G)$ is just the group $\cs$-algebra of $G$. 
\end{example}

\section{A construction of $\csp{X,\theta,G}$}
\label{sec:construction} 

We will in this section for every discrete partial dynamical system
$(X,\theta,G)$ construct the $\cs$-algebra $\csp{X,\theta,G}$
(cf. Definition \ref{def:ox}) as a crossed product of a $\cs$-partial
dynamical systems (see Appendix \ref{sec:parcro} for an short
introduction to crossed product of $\cs$-partial dynamical systems,
cf. also \cites{MR2003f:46108,MR99a:46121,MR1331978}). 

We will first from the discrete partial dynamical system
$(X,\theta,G)$ construct the $\cs$-partial dynamical systems
$(\A,G,\ba{\theta})$ which $\csp{X,\theta,G}$ is a crossed product of. 

Let $(D_g)_{g\in G}$ denote the domains and $(\theta_g)_{g\in G}$ the
partial one-to-one maps of $\theta$. Recall that $\B$ is the Boolean
algebra generated by the domains $\{D_g\mid g\in G\}$ of the partial
action $\theta$. Notice that $\theta_h(D_i)=D_h\cap D_{hi}\in \B$ for
every $h,i\in G$, so $\theta_g(A)\in \B$ for every $A\in \B$ and every $g\in G$.
 
Let $\X$ be the dual of $\B$: i.e., $\X$ is the closed subset
\begin{equation*}
  \Bigl\{\phi\in\{0,1\}^{\B}\bigm| \phi \text{ is a Boolean homomorphism}\Bigr\}
\end{equation*}
of the Cantor space $\{0,1\}^{\B}$ endowed with the product topology
of the discrete topology on $\{0,1\}$. Then $\X$ is a totally
disconnected compact Hausdorff space (a \emph{Boolean space},
cf. \cite{MR29:4713}*{\S 18}). 

For each $A\in \B$, let $\h{A}=\{\phi\in \X\mid \phi(A)=1\}$, and
notice that $\h{A}$ is a clopen subset of $\X$. We then have that the
map $A\mapsto \h{A}$ is a Boolean isomorphism between $\B$ and the
Boolean algebra of clopen subsets of $\X$ (cf. \cite{MR29:4713}*{\S18}).

\begin{lemma} \label{separet}
  The system $\{\h{D}_t\mid t\in G\}$ separates points in $\X$.
\end{lemma}

\begin{proof}
  Let $\phi_1,\phi_2\in \X$ and let $\Aa$ be the subset of $\B$
  defined by
  \begin{equation*}
    \Aa=\{A\in \B\mid \phi_1(A)=\phi_2(A)\}.
  \end{equation*}
  Then $\Aa$ is a Boolean subalgebra of $\B$. Assume that the
  following equivalence 
  \begin{equation*}
    \phi_1\in\h{D}_g\iff \phi_2\in \h{D}_g
  \end{equation*}
  holds for every $g\in G$. That means that $D_g\in \Aa$ for every
  $g\in G$ and thus that $\Aa=\B$. But then $\phi_1$ and $\phi_2$ must
  be equal.  
\end{proof}
For each $g\in G$, let $\h{\theta}_g$ be the map given by
\begin{equation*}
  \h{\theta}_g(\phi)(A)=\phi\bigl(\theta_g\inv(A)\bigr)
\end{equation*}
for $A\in \B$ and $\phi\in \h{D}_{g\inv}$. It is easy to check that
$\h{\theta}_g$ is a homeomorphism from $\h{D}_{g\inv}$ to $\h{D}_{g}$
with $\h{\theta}_{g\inv}$ as its inverse (cf. \cite{MR29:4713}*{\S 20}). 

Let $\A=C(\h{X})$ and let for each $g\in G$, $\ba{D}_g$ be the subset
of $\A$ defined by
\begin{equation*}
  \ba{D}_g=\left\{f\in \A\bigm| f_{|\h{X}\setminus \h{D}_g}=0\right\}, 
\end{equation*}
and let $\ba{\theta}_g:\ba{D}_{g\inv}\to\ba{D}_g$ be defined by
\begin{equation*}
  \ba{\theta}_g(f)(\phi)=
  \begin{cases}
    f(\h{\theta}_{g\inv}(\phi))&\text{if }\phi\in\h{D}_g,\\
    0&\text{if }\phi\in \h{X}\setminus\h{D}_g,
  \end{cases}
\end{equation*}
for $f\in \ba{D}_{g\inv}$ and $\phi\in \h{X}$.
Then $\ba{\theta}=\bigl((\ba{D}_g)_{g\in G},(\ba{\theta}_g)_{g\in G}\bigr)$
is a partial action of $G$ on the $\cs$-algebra $\A$
(cf. Appendix \ref{sec:parcro}). Thus $(\A,G,\ba{\theta})$ is a $\cs$-partial
dynamical system. 

\begin{theorem} \label{hoved}
  Let $(X,\theta,G)$ be a discrete partial dynamical
  system and let $\A$ and $\ba{\theta}$ be as defined above. Then the
  partial crossed product $\A\rtimes_{\ba{\theta}}G$ of 
  the $\cs$-partial dynamical system $(\A,G,\ba{\theta})$ is
  generated by a family of elements $(s_g)_{g\in G}$ satisfying condition
  \eqref{eq:ox}--\eqref{eq:oxd} of Definition \ref{def:ox}, and if
  $\Aa$ is another $\cs$-algebra with a family of elements $(S_g)_{g\in G}$
  satisfying condition \eqref{eq:ox}--\eqref{eq:oxd}, then there
  exists a $*$-homomorphism from $\A\rtimes_{\ba{\theta}}G$ to $\Aa$
  which maps $s_g$ to $S_g$ for every $g\in G$.
 \end{theorem}

\begin{proof}
  Let for every $g\in G$, $s_g=\delta_g$, where $\{\delta_g\}_{g\in G}$
  is a described in Appendix \ref{sec:parcro}. Since $(\pi,u)\mapsto
  \pi\times u$ is a bijective correspondence between covariant
  representations of $(\A,G,\ba{\theta})$ and non-degenerated
  representations of $\A\rtimes_{\ba{\theta}}G$, there exists a
  covariant representation $(\pi,u)$ of $(\A,G,\ba{\theta})$ on a
  Hilbert space $\H$ such that $\pi\times u$ is a faithful
  non-degenerate representation of $\A\rtimes_{\ba{\theta}}G$. We then
  have that  $u$ is a partial representation of $G$ and so the family
  $\bigl(u(g)\bigr)_{g\in G}$ satisfies condition
  \eqref{eq:ox}--\eqref{eq:oxc} of Definition \ref{def:ox}.  

  Denote for every subset $B$ of $\A$, the subspace
  \begin{equation*}
    \spc\{\pi(T)\xi\mid T\in B,\ \xi\in \H\} 
  \end{equation*}
  of $\H$ by $[B]$ and the projection of $\H$ onto $[B]$ by $\proj([B])$. 

  Since the map sending an element $A$ of $\B$ to $\h{A}$, the map
  sending a clopen subset $V$ of $\h{X}$ to $\bigl\{f\in \A\bigm|
  f_{|\h{X}\setminus V}=0\bigr\}$, and the map sending an ideal $I$ of $\A$
  to $\proj([I])$ all are Boolean homomorphism, so is the map sending
  an element $A$ of $\B$ to
  $\proj\bigl(\bigl[\bigl\{f\in \A\mid f_{|\h{X}\setminus \h{A}}=
  0\bigr\}\bigr]\bigr)$, 
  and since we have that
  \begin{equation*}
    \begin{split}
      \proj\bigl(\bigl[\bigl\{f\in \A\mid f_{|\h{X}\setminus \h{D}_g}=
      0\bigr\}\bigr]\bigr)
      &=\proj([\ba{D}_g])\\
      &=u(g)u(g)^*
    \end{split}
  \end{equation*}
  for every $g\in G$, $\bigl(u(g)\bigr)_{g\in G}$ satisfies \eqref{eq:oxd} of
  Definition \ref{def:ox}. Since $\pi\times u$ is faithful and
  $\pi\times u(s_g)=u(g)$ for every $g\in G$, $(s_g)_{g\in G}$
  satisfies \eqref{eq:ox}--\eqref{eq:oxd} of Definition \ref{def:ox}. 

  Remember that $\A\rtimes_{\ba{\theta}}G$ is generated by $\A$ and
  $(s_g)_{g\in G}$. It follows from Lemma \ref{separet} and
  the Stone-Weierstrass Theorem that the equality 
  \begin{equation*}
    \spc\{1_{\h{D}_g}\mid g\in G\}=\A
  \end{equation*}
  holds, and since $s_gs_g^*=1_{\h{D}_g}$ for every $g\in G$, this
  shows that $\A\rtimes_{\ba{\theta}}G$ is generated by $(s_g)_{g\in G}$.  

  Now let $\Aa$ be another $\cs$-algebra with a family 
  $(S_g)_{g\in G}$ of elements which satisfies condition
  \eqref{eq:ox}--\eqref{eq:oxd}. Let $\psi$ be a non-degenerate
  faithful representation of $\Aa$ on a Hilbert space $\H$. Since
  $S_e$ is the unit of $\Aa$, $\psi(S_e)\xi=\xi$ for every $\xi\in\H$.

  Let for every $g\in G$, $U(g)=\psi(S_g)$. Then $U$ is a partial
  representation of $G$ on $\H$. Since $\spc\{1_{\h{D}_g}\mid g\in
  G\}=\A$, since the map $A\mapsto \h{A}$ is a Boolean isomorphism between
  $\B$ and the Boolean algebra of clopen subsets of $\X$, and since the map
  $D_g\mapsto S_gS_g^*$ extends to a Boolean homomorphism from $\B$ to
  the set of projections in $\cs\bigl(\{S_g S_g^*\mid g\in G\}\bigr)$,
  there exists by Lemma \ref{lemma:homo} a $*$-homomorphism
  from $\A$ to $\Aa$ which maps $1_{\h{D}_g}$ to $S_gS_g^*$ for every
  $g\in G$. Let us denote the composition of this $*$-homomorphism
  with $\psi$ by $\eta$. Then $\eta$ is a representation of $\A$ on
  $\H$, and since $\eta(1)\xi=\psi(S_e)\xi=\xi$ for every $\xi\in\H$,
  $\eta$ is non-degenerate.

  If $g\in G$, then we have that
  \begin{equation*}
    \begin{split}
      U(g)U(g)^*&=\psi(S_gS_g^*)\\
      &=\eta(1_{\h{D}_g})\\
      &=\proj\Bigl(\spc\bigl\{\eta(1_{\h{D}_g})\xi\mid \xi\in H\bigr\}\Bigr)\\
      &=\proj\bigl([\ba{D}_g]\bigr).    
    \end{split}
  \end{equation*}   
  Let $i,h\in G$. It then follows from \eqref{eq:paractb} that
  $\ba{\theta}_h\bigl(1_{\h{D}_{h\inv}\cap\h{D}_i}\bigr)=
  1_{\h{D}_h\cap\h{D}_{hi}}$. Thus we have that 
  \begin{equation*}
    \begin{split}
      \eta\Bigl(\ba{\theta}_h\bigl(1_{\h{D}_{h\inv}\cap\h{D}_i}\bigr)\Bigr)
      &=\eta\bigl(1_{\h{D}_h\cap\h{D}_{hi}}\bigr)\\ 
      &=U(h)U(h)^*U(hi)U(hi)^*\\
      &=U(h)U(i)U(hi)^*\\
      &=U(h)U(i)U(i\inv)U(i\inv)^*U(i\inv h\inv)\\
      &=U(h)U(h)^*U(h)U(i)U(i\inv)U(h\inv)\\
      &=U(h)\eta\bigl(1_{\h{D}_{h\inv}\cap\h{D}_i}\bigr)U(h\inv),
    \end{split}
  \end{equation*}
  where the third equality follows from \eqref{eq:oxc}, the forth
  from \eqref{eq:ox} and \eqref{eq:oxb}, and the fifth from
  \eqref{eq:ox} and \eqref{eq:oxd}. It follows from Lemma \ref{separet} and
  the Stone-Weierstrass Theorem that 
  $\spc\bigl\{1_{\h{D}_{h\inv}\cap\h{D}_i}\mid i\in
  G\bigr\}=\ba{D}_{h\inv}$, so the above computation shows that the
  equality 
  \begin{equation*}
    \eta\bigl(\ba{\theta}_h(f)\bigr)=U(h)\bigl(\eta(f)\bigr)U(h\inv)
  \end{equation*}
  holds for every $f\in \ba{D}_{h\inv}$. 
  
  Thus $(\eta,U)$ is a covariant representation of
  $(\A,G,\ba{\theta})$, and so there exists a non-degenerate
  representation of $\A\rtimes_{\ba{\theta}}G$ on $\H$ which maps
  $s_g$ to $\psi(S_g)$ for every $g\in G$, and thus a
  $*$-homomorphism from  $\A\rtimes_{\ba{\theta}}G$ to $\Aa$ which maps
  $s_g$ to $S_g$ for every $g\in G$.
\end{proof}

\begin{remark} \label{remark:exellaca}
  Let $\ind$ be an arbitrary index set and let
  $A={A(i,j)}_{i,j\in\ind}$ be a matrix with entries in $\{0,1\}$ and
  having no zero rows. In \cite{MR2000i:46064} Exel and Lace
  associated to $A$ a $\cs$-algebra $\O_A$ (cf. Section
  \ref{subsec:CK}). We will in Theorem \ref{theorem:CK} show that the
  unitization $\widetilde{\O}_A$ of $\O_A$ is isomorphic to
  $\cspxa$ for a certain discrete partial dynamical system $\pxa$. 

  Exel and Laca constructed in \cite{MR2000i:46064} a partial action 
  $\left((\Delta^A_g)_{g\in\fg\ind},(h^A_g)_{g\in\fg\ind}\right)$
  of $\fg\ind$ on a compact topological space
  $\widetilde{\Omega}_A$ such that $\widetilde{\O}_A$ is isomorphic to
  the crossed product of this partial action. 

  If we let $(D_g)_{g\in G}$ and $(\theta_g)_{g\in G}$ denote the domains
  and partial one-to-one maps of $\theta_\XA$, and $\h{\mathsf{X}}_A^+$,
  $\h{D}_g$ and $\h\theta_g$ be as above, then one can show that the map
  \begin{equation*}
    \phi\mapsto \bigl\{g\in\fg\ind\mid \phi(D_g)=1\bigr\}
  \end{equation*}
  is a homeomorphism from $\h{\mathsf{X}}_A^+$ to
  $\widetilde{\Omega}_A$ which for every 
  $g\in\fg\ind$ maps $\h{D}_g$ to $\Delta^A_g$ and intertwines
  $\h\theta_g$ and $h^A_g$. This fact lays the foundation for an
  alternative proof of Theorem \ref{theorem:CK}.
\end{remark}

\section{A representation of $\csp{X,\theta,G}$} \label{sec:repr}
Let $(X,\theta,G)$ be a discrete partial dynamical system and let
$(D_g)_{g\in G}$ denote the domains of $\theta$ and $(s_g)_{g\in G}$ the
generators of $\csp{X,\theta,G}$. If $\pi$ is a representation of
$\csp{X,\theta,G}$, then $\pi$ will induce a Boolean homomorphism from
the Boolean algebra of projections in the unital abelian $\cs$-algebra
$\cs\bigl(\{s_gs_g^*\in\csp{X,\theta,G}\mid g\in G\}\bigr)$ to the
Boolean algebra of projections in the unital abelian $\cs$-algebra 
$\cs\bigl(\{\pi(s_gs_g^*)\mid g\in G\}\bigr)$ which maps $s_gs_g^*$ to
$\pi(s_gs_g^*)$ for every $g\in G$, and by composing this Boolean
homomorphism with $\boomap$, we get a Boolean homomorphism from
$\Bx(X,\theta,G)$ to the Boolean algebra of projections in the unital
abelian $\cs$-algebra $\cs\bigl(\{\pi(s_gs_g^*)\mid g\in G\}\bigr)$
which maps $D_g$ to $\pi(s_gs_g^*)$ for every $g\in G$. 

We will in this section for every discrete partial dynamical system 
$(X,\theta,G)$ construct a representation $\pi$ of $\csp{X,\theta,G}$
such that the Boolean homomorphism mentioned above is injective. As
a corollary of this, we see that the Boolean homomorphism $\boomap$ is
injective and thus that $\csp{X,\theta,G}\ne 0$ (unless $X=\emptyset$).

Let $(X,\theta,G)$ be a discrete partial dynamical system and let
$(D_g)_{g\in G}$ and $(\theta_g)_{g\in G}$ denote the domains and partial
one-to-one maps of $\theta$.
Let $(e_x)_{x\in X}$ be an orthonormal basis for the Hilbert space
$l_2(X)$, and define for each $g\in G$ an operator $S_g$ by letting 
\begin{equation} \label{eq:op}
  S_g\left(\sum_{x\in X}\lambda_xe_x\right)=\sum_{x\in
    D_g}\lambda_{\theta_{g\inv}(x)}e_x 
\end{equation}
for every $\sum_{x\in X}\lambda_xe_x\in l_2(X)$.

It is straight forward to check that these operators are partial
isometries, that $S_e=1$ and that for every $h,i\in G$, $S_h^*=S_{h\inv}$,
$S_hS_i=S_hS_h^*S_{hi}$ and $S_hS_h^*=\proj\bigl(\spc\{e_x\mid x\in
D_h\}\bigr)$, where $\proj\bigl(\spc\{e_x\mid x\in D_h\}\bigr)$ is the
orthogonal 
projection of $l_2(X)$ onto $\spc\{e_x\mid x\in D_h\}$. Since the map
\begin{equation*}
  A \mapsto \proj\bigl(\spc\{e_x\mid x\in A\}\bigr)
\end{equation*}
is a Boolean homomorphism from $\B$ to the Boolean algebra of
projections in the unital abelian $\cs$-algebra 
$\cs\bigl(\{S_gS_g^*\mid g\in G\}\bigr)$, the family $(S_g)_{g\in G}$
of operators satisfies condition \eqref{eq:ox}--\eqref{eq:oxd} of
Definition \ref{def:ox}. 

Thus we have:
\begin{proposition}
  Let $(X,\theta,G)$ be a discrete partial dynamical
  system, let $(s_g)_{g\in G}$ denote the generators of $\csp{X,\theta,G}$
  and $(S_g)_{g\in G}$ the operators defined by \eqref{eq:op}. Then
  there is a $*$-homomorphism from $\csp{X,\theta,G}$ to 
  the $\cs$-algebra of bounded operators on the Hilbert space
  $l_2(X)$, sending $s_g$ to $S_g$ for every $g\in G$.
\end{proposition}

\begin{corollary} \label{boomapinj}
  The Boolean homomorphism $\boomap$ (cf. Definition \ref{remarkboomap})
  is for every discrete partial dynamical system $(X,\theta,G)$ injective.
\end{corollary}

\begin{proof}
  Since the Boolean homomorphism
  \begin{equation*}
    A \mapsto \proj\bigl(\spc\{e_x\mid x\in A\}\bigr)
  \end{equation*}
  from $\B$ to the set of projections in
  $\cs\bigl(\{S_g S_g^*\mid g\in G\}\bigr)$ is injective, so is $\boomap$.
\end{proof}

\section{$\cs$-algebras associated to symbolic dynamical
  systems} \label{retu} 
We will now show that the class of $\cs$-algebras associated to discrete
partial dynamical systems generalizes Cuntz-Krieger algebras (both
for finite and infinite matrices), $\cs$-algebras associated to
one-sided shift spaces and crossed products of two-sided shift spaces.

We will do that by regarding one-sided and two-sided shift spaces and
topological Markov chains as symbolic dynamical systems and thus
associated to them the discrete partial dynamical system of Definition
\ref{onesideddef} and \ref{twosideddef} and then show that the
$\cs$-algebras of these discrete partial dynamical systems are
isomorphic to the $\cs$-algebra associated to the one-sided shift
space,the crossed product of the two-sided shift space, and the
unitization of the Cuntz-Krieger algebra of the transition matrix of
the Markov chain, respectively.

Before we do that let us briefly look at the general structure of
$\csp{X,\theta,G}$ when $(X,\theta,G)$ is the discrete partial
dynamical system associated to a one- or two-sided symbolic dynamical
system as in Definition \ref{onesideddef} and \ref{twosideddef}.

\begin{lemma} \label{lemma:onestruc}
  Let $(\OSD,\osd)$ be a one-sided symbolic dynamical system over the
  alphabet $\al$ and let $(\OSD,\theta_\OSD,\fg{\al})$ be the discrete
  partial dynamical system associated to $(\OSD,\osd)$ as in
  Definition \ref{onesideddef}. Then the following holds for the
  generators $(s_g)_{g\in\fg\al}$ of $\csp{\OSD,\theta_\OSD,\fg{\al}}$:
  \begin{enumerate}
  \item $s_{b_1}s_{b_2}\dotsm s_{b_k}=s_g$ for $g\in
    \fg\al$ written in reduced form, $b_1b_2\dotsm b_k$, \label{enum:oneet}
  \item if $s_g\ne 0$, then there exist $\mu,\nu\in\al^*$ such that
    $g=\mu\nu\inv$, \label{enum:oneto}
  \item $s_\mu^*s_\nu=0$ if $\mu,\nu\in\al^*$, $|\mu|=|\nu|$ and
    $\mu\ne\nu$. \label{enum:onetre}
  \end{enumerate}
\end{lemma}

\begin{proof}
  If $b\in \al\cup\al\inv$ and $g\in \fg\al$ written in reduced form
  does not begin with $b$, then $D_{bg}\subseteq D_b$, so we have that
  \begin{equation*}
    \begin{split}
      s_bs_g&=s_bs_b^*s_{bg}\\
      &= s_bs_b^*s_{bg}s_{bg}^*s_{bg}\\
      &=s_{bg}.
    \end{split}
  \end{equation*}
  This shows that $s_{b_1}s_{b_2}\dotsm s_{b_k}=s_g$ for
  $g\in \fg\al$ written in reduced form $b_1b_2\dotsm b_k$.

   (\ref{enum:oneto}) and (\ref{enum:onetre}) easily follow from Lemma
   \ref{lemma:onefund} and \ref{lemma:orto}.
\end{proof}

\begin{lemma} \label{lemma:twostruc}
  Let $(\TSD,\tsd)$ be a two-sided symbolic dynamical system over the
  alphabet $\al$ and let $(\TSD,\theta_\TSD,\fg{\al})$ be the discrete
  partial dynamical system associated to $(\TSD,\tsd)$ as in
  Definition \ref{twosideddef}. Then the following holds for the
  generators $(s_g)_{g\in\fg\al}$ of $\csp{\TSD,\theta_\TSD,\fg{\al}}$:
  \begin{enumerate}
  \item $s_{b_1}s_{b_2}\dotsm s_{b_k}=s_g$ for $g\in
    \fg\al$ written in the reduced form $b_1b_2\dotsm b_k$, \label{enum:twoet}
  \item if $s_g\ne 0$, then either $g\in\al^*$ or
    $g\inv\in\al^*$, \label{enum:twoto}
  \item $s_h^*s_i=s_hs_i^*=0$ if $h,i\in\fg\al$, $[h]=[i]$ and $h\ne
    i$. \label{enum:twotre}
  \end{enumerate}
\end{lemma}

\begin{proof}
  (\ref{enum:twoet}) can be proved in exactly
  the same way (\ref{enum:oneet}) was proved
  in Lemma \ref{lemma:onestruc}, (\ref{enum:twoto}) easily follows
  from Lemma \ref{lemma:twofund} and \ref{lemma:orto}, and
  (\ref{enum:twotre}) follows from (\ref{enum:twoto}) and Lemma
  \ref{lemma:twofund} and \ref{lemma:orto}. 
\end{proof}

\subsection{Crossed products of two-sided shift
  spaces} \label{sec:cross} 
Let $(\TSS,\tsh)$ be a two-sided shift space over the finite alphabet
$\al$ (cf. \cite{MR1484730} and \cite{MR97a:58050}). That is:
$\tsh:\al^{\Z}\to\al^{\Z}$ is the map defined by
\begin{equation} \label{eq:twoshift}
  \bigl(\tsh((z_n)_{n\in\Z})\bigr)_m=z_{m+1}
\end{equation}
for every $(z_n)_{n\in\Z}\in\al^{\Z}$ and every $m\in\Z$, and $\TSS$
is a closed (in the product topology of the discrete topology of
$\al$) subset of $\al^{\Z}$ such that $\tsh(\TSS)=\TSS$.

Let $\tsh^\star$ be the automorphism on $C(\TSS)$ defined by $f\mapsto
f\circ\tsh$ and let $C(\TSS)\rtimes_{\tsh^\star}\Z$ be the full
crossed product of the $\cs$-dynamical system $C(\TSS,\tsh^\star,\Z)$
(cf. \cite{MR81e:46037}*{7.6.5}). Thus $C(\TSS)\rtimes_{\tsh^\star}\Z$
is the universal $\cs$-algebra generated by a copy of $C(\TSS)$ and an
unitary operator $U$ which satisfies that $UfU^*=f\circ \tsh$ for
every $f\in C(\TSS)$.

Since $(\TSS,\tsh)$ is a two-sided symbolic dynamical system, we can
associate to it the discrete partial dynamical system
$(\TSS,\theta_\TSS,\fg{\al})$ of Definition \ref{twosideddef}. We then
have the following theorem:
 
\begin{theorem} \label{kryds}
  Let $(\TSS,\tsh)$ be a two-sided shift space and let
  $(\TSS,\theta_\TSS,\fg{\al})$ be the discrete partial dynamical system
  associated to $(\TSS,\tsh)$ as done in Definition \ref{twosideddef}. Then
  $\csp{\TSS,\theta_\TSS,\fg{\al}}$ is isomorphic to the crossed product
  $C(\TSS)\rtimes_{\tsh^\star}\Z$. 

  More precisely: if $(D_g)_{g\in\fg\al}$ denotes the domains of
  $\theta_\TSS$, $(s_g)_{g\in\fg\al}$ denotes the generators of
  $\csp{\TSS,\theta_\TSS,\fg{\al}}$ and $U$ is as above, then we have that
  \begin{equation*}
    C(\TSS)=\spc\{1_{D_g}\mid g\in \fg{\al}\},
  \end{equation*} 
  and there exists a $*$-isomorphism from $C(\TSS)\rtimes_{\tsh^\star}\Z$ to
  $\csp{\TSS,\theta_\TSS,\fg{\al}}$ which maps $1_{D_g}$ to $s_gs_g^*$ for
  every $g\in\fg\al$ and $U$ to $\sum_{a\in\al}s_a$.
\end{theorem}

\begin{proof}
  We will first show that the Boolean algebra
  $\Bx(\TSS,\theta_\TSS,\fg\al)$ generated by 
  $\{D_g\mid g\in \fg{\al}\}$ is the Boolean algebra of clopen subsets 
  of $\TSS$. It easily follows from Lemma \ref{lemma:twofund} that 
  every set in $\Bx(\TSS,\theta_\TSS,\fg\al)$ is clopen.  
  In the other direction, we have by Lemma \ref{lemma:twofund} that
  \begin{equation*}
    \bigl\{z\in \TSS\mid z_{[-|u|,|v|[}=uv\bigr\}= D_{u\inv}\cap D_{v} \in
    \Bx(\TSS,\theta_\TSS,\fg\al) 
  \end{equation*}
  for $u,v\in\al^*$, and since  the system
  consisting of sets of this form is a basis for the topology of
  $\TSS$, every clopen set is a finite union of sets of this form
  and thus belongs to $\Bx(\TSS,\theta_\TSS,\fg\al)$.
  
  So it follows from the Stone-Weierstrass theorem that
  \begin{equation*}
    C(\TSS)=\spc\{1_{D_g}\mid g\in \fg{\al}\}. 
  \end{equation*}
  Since $\boomap[(\TSS,\theta_\TSS,\fg\al)]$ (cf. Definition
  \ref{remarkboomap}) 
  is a Boolean homomorphism from $\Bx(\TSS,\theta,\fg\al)$ to the Boolean
  algebra of projections in the unital abelian $\cs$-algebra
  \begin{equation*}
    \cs\bigl(\{s_gs_g^*\in\csp{\TSS,\theta_\TSS,\fg\al}\mid g\in \fg\al\}\bigr) 
  \end{equation*}
  which maps $D_g$ to $s_gs_g^*$, it follows from Lemma \ref{lemma:homo}
  that there exists a $*$-homomorphism $\crosmap$ from $C(\TSS)$ to
  $\csp{\TSS,\theta_\TSS,\fg{\al}}$ which maps $1_{D_g}$ to $s_gs_g^*$
  for every $g\in \fg\al$. Let $u$ be the element $\sum_{a\in\al}s_a$ in
  $\csp{\TSS,\theta_\TSS,\fg{\al}}$. Then we have that
  \begin{equation*}
    \begin{split}
      u\crosmap(1_{D_g})u^*& =\left(\sum_{a\in\al}s_a\right)
      \bigl(s_gs_g^*\bigr) 
      \left(\sum_{a'\in\al}s_{a'}^*\right) \\
      &=\sum_{\substack{a\in\al\\
          a'\in\al}}s_as_gs_g^*s_a^*s_as_{a'}^*s_{a'}s_{a'}^* \\      
      &=\sum_{a\in\al}s_as_gs_g^*s_a^*\\
      &=\sum_{a\in\al}s_as_a^*s_{ag}s_{ag}^*s_{a}s_{a}^*\\
      &=\sum_{a\in\al}\crosmap\bigl(1_{D_a\cap D_{ag}}\bigr)\\
      &=\sum_{a\in\al}\crosmap\bigl(1_{\theta_a(D_g)}\bigr)\\
      &=\crosmap\bigl(1_{\cup_{a\in\al}\theta_a(D_g)}\bigr)\\
      &=\crosmap\bigl(1_{\tsh\inv(D_g)}\bigr)\\
      &=\crosmap(1_{D_g}\circ\tsh)
    \end{split}
  \end{equation*}
  for every $g\in\fg\al$, where the second equality follows from
  \eqref{eq:ox} and \eqref{eq:oxd}, the third from \eqref{eq:ox} and
  Lemma \ref{lemma:twostruc}\eqref{enum:twotre}, the forth from
  \eqref{eq:oxc} and the sixth from \eqref{eq:paractb}. Since
  $C(\TSS)=\spc\{1_{D_g}\mid g\in \fg{\al}\}$, this shows that
  \begin{equation*}
    u\crosmap(f)u^*=\crosmap(f\circ\tsh)
  \end{equation*} 
  for every $f\in C(\TSS)$. 
  
  Thus it follows from the universal property
  of $C(\TSS)\rtimes_{\tsh^\star}\Z$ that there exists a $*$-homomorphism
  $\tcrosmap$ from $C(\TSS)\rtimes_{\tsh^\star}\Z$ to
  $\csp{\TSS,\theta_\TSS,\fg{\al}}$ which is equal to $\crosmap$ on $C(\TSS)$
  and sends $U$ to $u=\sum_{a\in\al}s_a$.
  
  Let us now look at $C(\TSS)\rtimes_{\tsh^\star}\Z$. Let for every
  $a\in\al$, $S_a$ be the element $1_{D_a}U$ of
  $C(\TSS)\rtimes_{\tsh^\star}\Z$, 
  $S_{a\inv}=S_a^*$, and let $S_e=1$ and $S_g=S_{b_1}S_{b_2}\dotsm
  S_{b_k}$, where $g\in \fg\al$ is written in the reduced form $b_1b_2\dotsm 
  b_k$. Then $(S_g)_{g\in\fg{\al}}$ is a family of elements from
  $C(\TSS)\rtimes_{\tsh^\star}\Z$ 
  which clearly satisfies \eqref{eq:oxa} and \eqref{eq:oxb} 
  of Definition \ref{def:ox}. We will now show that $(S_g)_{g\in
    \fg{\al}}$ also satisfies \eqref{eq:ox}, \eqref{eq:oxc} and
  \eqref{eq:oxd} of Definition \ref{def:ox}.
  
  Remember (cf. Section \ref{sec:notat}) that $[\cdot]$ is the unique
  homomorphism from $\fg\al$ to $\Z$ such that $[a]=1$ for every $a\in\al$.
  According to Lemma \ref{lemma:twofund}, the map $\theta_g$ is for every
  $g\in\fg\al$ equal to 
  the restriction of $\tsh^{-[g]}$ to $D_{g\inv}$, so we have that
  \begin{equation*}
    D_h\cap\tsh^{-[h]}(D_i)=\theta_h(D_i)=D_{hi}\cap D_h
  \end{equation*}
  for all $h,i\in\fg\al$.
  
  If $a\in\al$ and $g\in\fg\al$ written in reduced form does
  not begin with $a\inv$, then $D_{ag}\subseteq D_a$ and
  $D_a\cap\tsh^{\inv}(D_g)=D_{ag}\cap D_a=D_{ag}$ and so we have that
  \begin{equation*}
    \begin{split}
      S_a1_{D_g}U^{[g]}&=1_{D_a}U1_{D_g}U^{[g]}\\
      &=1_{D_a}U1_{D_g}U^*U^{[g]+1}\\
      &=1_{D_a\cap\tsh\inv(D_g)}U^{[g]+1}\\
      &=1_{D_{ag}}U^{[g]+1},\\
    \end{split}
  \end{equation*}
  and if $g$ written in reduced form does not begin with an $a$,
  then $D_{a\inv g}\subseteq D_{a\inv}$ and $\tsh(D_a\cap
  D_g)=\theta_{a\inv}(D_g)=D_{a\inv}\cap D_{a\inv g}=D_{a\inv g}$,
  so we have that
  \begin{equation*}
    \begin{split}
      S_{a\inv}1_{D_g}U^{[g]}&=U^*1_{D_a}1_{D_g}U^{[g]}\\
      &=1_{\tsh(D_a\cap D_g)}U^{[g]-1}\\
      &=1_{D_{a\inv g}}U^{[g]-1}.
    \end{split}
  \end{equation*}
  This shows that $S_g=1_{D_g}U^{[g]}$ and thus that $S_gS_g^*=1_{D_g}$
  for every $g\in \fg\al$. Hence $(S_g)_{g\in\fg{\al}}$ satisfies
  \eqref{eq:ox} and \eqref{eq:oxd}  
  of Definition \ref{def:ox}. If $h,i\in\fg\al$, then we have that
  \begin{equation*}
    \begin{split}
      S_hS_i&=1_{D_h}U^{[h]}1_{D_i}U^{[i]}\\
      &=1_{D_h}U^{[h]}1_{D_i}U^{-[h]}U^{[h]}U^{[i]}\\
      &=1_{D_h\cap\tsh^{-[h]}(D_i)}U^{[hi]}\\
      &=1_{D_{hi}\cap D_h}U^{[hi]}\\
      &=1_{D_h}1_{D_{hi}}U^{[hi]}\\
      &=S_hS_h^*S_{hi},
    \end{split}
  \end{equation*}
  which shows that $(S_g)_{g\in\fg{\al}}$ also satisfies
  \eqref{eq:oxc} of Definition \ref{def:ox}.
  
  Thus it follows from the universal property of
  $\csp{\TSS,\theta_\TSS,\fg{\al}}$ that there is a
  $*$-homomorphism $\psi$ from $\csp{\TSS,\theta_\TSS,\fg{\al}}$ to
  $C(\TSS)\rtimes_{\tsh^\star}\Z$ such that
  $\psi(s_g)=S_g=1_{D_g}U^{[g]}$ for every $g\in \fg\al$.
  
  We have that 
  \begin{equation*}
    \psi\bigl(\tcrosmap(U)\bigr)= \psi\bigg(\sum_{a\in\al}s_a\bigg)=
    \sum_{a\in\al}1_{D_a}U=U,
  \end{equation*}
  and that
  \begin{equation*}
    \psi\bigl(\tcrosmap(1_{D_g})\bigr)=
    \psi(s_gs_g^*)= 1_{D_g}U^{[g]}U^{-[g]}1_{D_g}= 1_{D_g}
  \end{equation*}
  for every $g\in \fg\al$, and since $C(\TSS)\rtimes_{\tsh^\star}\Z$
  is generated by $U$ and $\{1_{D_g}\mid g\in\fg{\al}\}$, this shows that
  $\psi\circ\tcrosmap=\Id_{C(\TSS)\rtimes_{\tsh^\star}\Z}$.
  
  We also have that
  \begin{equation*}
    \tcrosmap\bigl(\psi(s_g)\bigr)=\tcrosmap\bigl(1_{D_g}U^{[g]}\bigr)=
    s_gs_g^*\bigl(\sum_{a\in\al}s_a\bigr)^{[g]}
  \end{equation*} 
  for every $g\in \fg\al$. It follows from Lemma
  \ref{lemma:twostruc} that we have that
  \begin{equation*}
    \biggl(\sum_{a\in\al}s_a\biggr)^{[g]}=\sum_{\substack{h\in\fg\al\\ [h]=[g]}}s_h,
  \end{equation*}
  and that 
  \begin{equation*}
    s_gs_g^*s_h=
    \begin{cases}
      s_g&\text{if }g=h,\\
      0&\text{if }g\ne h,
    \end{cases}
  \end{equation*}
  for every $h\in \fg\al$ with $[h]=[g]$. Thus it follows that
  \begin{equation*}
    \tcrosmap(\psi(s_g))=s_gs_g^*\biggl(\sum_{a\in\al}s_a\biggr)^{[g]}=s_g
  \end{equation*}
  for every $g\in\fg\al$, which shows that
  $\tilde{\crosmap}\circ\psi=\Id_{\csp{\TSS,\theta_\TSS,\fg{\al}}}$. 
  
  Hence $\tcrosmap$ is an isomorphism from
  $C(\TSS)\rtimes_{\tsh^\star}\Z$ to $\csp{\TSS,\theta_\TSS,\fg{\al}}$
  which maps $1_{D_g}$ to $s_gs_g^*$ for
  every $g\in\fg\al$ and $U$ to $\sum_{a\in\al}s_a$.
\end{proof}

\subsection{$\cs$-algebras associated to one-sided shift
  spaces} \label{sec:onesided} 
Let $(\OSS,\osh)$ be a one-sided shift space over the finite alphabet
$\al$ (cf. \cite{MR1484730} and \cite{MR97a:58050}*{\S13.8}). That is:
$\osh:\al^{\No}\to\al^{\No}$ is the map
\begin{equation} \label{eq:oneshift}
  x_0x_1x_2\dotsm\mapsto x_1x_2\dotsm,
\end{equation}
and $\OSS$ is a closed (in the product topology of the discrete
topology of $\al$) subset of $\al^{\No}$ such that $\tsh(\OSS)=\OSS$.
  
As far as the author know, Kengo Matsumoto was in
\cite{MR1454478} the first to consider $\cs$-algebras associated
to shift spaces. Matsumoto's construction is however in the opinion of
the author not the optimal one (see \cite{MR2091486} for a discussing
of this matter). In \cite{tmc} the author considered a different
construction of $\cs$-algebras associated to shift spaces which for
some shift spaces gives a slightly different $\cs$-algebra than
Matsumoto's (and the $\cs$-algebra considered in \cite{MR2091486},
cf. \cite{CS}*{Section 7}). We will in this paper work with the
$\cs$-algebra $\O_{\OSS}$  of \cite{tmc} (it is isomorphic to the $\cs$-algebra
$\Dx\negthickspace\rtimes_{\alpha,\mathcal{L}}\N$ considered in \cite{CS}, cf.
\cite{CS}*{Remark 9}). It can be characterized in the following way
(cf. \cite{tmc}*{Remark 7.3}):

Let $\al^*$ denote the set of finite words with letters from $\al$.
For $\mu,\nu\in \al^*$, let $\cy\mu\nu$ be the subset of $\OSS$
defined by 
\begin{equation*}
  \cy\mu\nu=\bigl\{x\in \OSS\bigm| x_{[0,|\nu|[}=\nu,\ \mu
  x_{[|\nu|,\infty[}\in \OSS\bigr\}.
\end{equation*}
We let $\mathfrak{B}(\OSS)$ be the abelian
$C^*$-algebra of all bounded functions on $\OSS$, and $\Dx$ the
$C^*$-subalgebra of $\mathfrak{B}(\OSS)$ generated by
$\bigl\{1_{\cy\mu\nu}\mid \mu,\nu \in\al^*\bigr\}$.

Then the $\cs$-algebra $\O_{\OSS}$ associated to the one-sided shift
space $(\OSS,\osh)$ is the universal $\cs$-algebra generated by a
family of partial isometries $(S_a)_{a\in\al}$ which satisfies that the map
\begin{equation*}
  1_{\cy\mu\nu}\mapsto S_\nu S_\mu^*S_\mu S_\nu^*
\end{equation*}
extends to a $*$-homomorphism from $\Dx$ to $\O_{\OSS}$, where
$S_\mu=S_{\mu_1}S_{\mu_2}\dotsm S_{\mu_k}$ for $\mu=\mu_1\mu_2\dotsm
\mu_k\in\al^*$ with $\mu_1,\mu_2,\dotsc ,\mu_k\in\al$, and
$S_\nu=S_{\nu_1}S_{\nu_2}\dotsm S_{\nu_l}$ for $\nu=\nu_1\nu_2\dotsm
\nu_l\in\al^*$ with $\nu_1,\nu_2,\dotsc ,\nu_l\in\al$. We will
denote this $*$-homomorphism (which in fact is injective) by
$\ioss$.

One should notice (cf. \cite{CS}*{Theorem 12}) that when $\OSS$ is a
topological Markov chain with transition matrix $A$
(cf. Section \ref{subsec:CK}), then $\O_{\OSS}$ is equal to the
Cuntz-Krieger algebra $\O_A$ (or to be more precise: to the universal
Cuntz-Krieger algebra considered by an Huef and Raeburn in
\cite{MR1452183}, cf. also \cite{MR561974}). 

We will for each $a\in \al$, by $\lambda_a$ denote the map on $\DX$
given by
\begin{equation*}
  \lambda_a(f)(x)=
  \begin{cases}
    f(ax)& \text{if }ax\in \OSS,\\
    0& \text{if }ax\notin \OSS,
  \end{cases}
\end{equation*}
and by $\phi_a$ the map on $\DX$ given by
\begin{equation*}
  \phi_a(f)(x)=
  \begin{cases}
    f(\osh(x))& \text{if }x\in D_a,\\
    0& \text{if }x\notin D_a,
  \end{cases}
\end{equation*}
for $f\in\Dx$ and $x\in \OSS$ (cf. \cite{tmc}*{Proposition 4.3 and
  Lemma 8.2}).

Since $(\OSS,\osh)$ is a one-sided symbolic dynamical system, we can
associate to it the discrete partial dynamical system
$(\OSS,\theta_{\OSS},\fg{\al})$ of Definition \ref{onesideddef}. We then
have the following theorem:

\begin{theorem} \label{theorem:one} 
  Let $(\OSS,\osh)$ be a one-sided shift space and let
  $(\OSS,\theta_{\OSS},\fg{\al})$ be the discrete partial dynamical system
  associated to $(\OSS,\osh)$ as done in Definition \ref{onesideddef}.
  Then $\csp{\OSS,\theta_{\OSS},\fg{\al}}$ is isomorphic to the
  $\cs$-algebra $\O_{\OSS}$ associated to $(\OSS,\osh)$.

  More precisely: if $(D_g)_{g\in\fg\al}$ denotes the domains of
  $\theta_{\OSS}$, $(s_g)_{g\in\fg\al}$ denotes the generators of
  $\csp{\OSS,\theta_{\OSS},\fg{\al}}$, and $\DX$, $\ioss$ and
  $(S_a)_{a\in\al}$ are as above, then we have that
  \begin{equation*}
    \spc\{1_{D_g}\mid g\in\fg\al\}=\DX,
  \end{equation*}
  and that there exists a 
  $*$-isomorphism from $\O_{\OSS}$ to $\csp{\OSS,\theta_{\OSS},\fg{\al}}$
  which maps $\ioss(1_{D_g})$ to $s_gs_g^*$ for every $g\in \fg\al$,
  and $S_a$ to $s_a$ for every $a\in\al$.
\end{theorem}

\begin{proof}
  To see that $\spc\{1_{D_g}\mid g\in\fg\al\}=\DX$ notice first that
  it follows from Lemma \ref{lemma:onefund} that
  if $\mu,\nu\in \alwords$, then we have that 
  \begin{equation*}
    \cy\mu\nu=\theta_\nu(D_{\mu\inv})=D_\nu\cap D_{\nu\mu\inv}.
  \end{equation*}
  So $\DX\subseteq \spc\{1_{D_g}\mid g\in\fg\al\}$.
  
  If $A$ is a subset of $\OSS$ such that $1_A\in \DX$, then
  $1_{\theta_a(A)}=\lambda_a(1_A)\in \DX$ and
  $1_{\theta_{a\inv}(A)}=\phi_a(1_A)\in \DX$ for $a\in \al$. So
  since we have that
  \begin{equation*}
    D_g=\theta_{b_1}\circ\theta_{b_2}\circ \dotsb \circ\theta_{b_k}(\OSS)
  \end{equation*}
  for every $g\in \fg{\al}$ written in the reduced form $b_1b_2\dotsm
  b_k$, we have that $\spc\{1_{D_g}\mid g\in\fg\al\}\subseteq \DX$.
  Thus $\DX=\spc\{1_{D_g}\mid g\in\fg\al\}$.
  
  Since $\boomap[(\OSS,\theta_{\OSS},\fg\al)]$ is a Boolean
  homomorphism from $\Bx(\OSS,\theta_{\OSS},\fg\al)$ to the set of
  projections in $\csp{\OSS,\theta_{\OSS},\fg{\al}}$ which maps $D_g$ to
  $s_gs_g^*$, it follows from Lemma \ref{lemma:homo} that there exists
  a $*$-homomorphism $\gamma$ from $\DX$ to
  $\csp{\OSS,\theta_{\OSS},\fg{\al}}$ which maps $1_{D_g}$ to
  $s_gs_g^*$ for every $g\in\fg\al$.
    
  For $g\in\fg\al$ and $a\in\al$, we have that 
  \begin{equation*}
    \begin{split}
      \gamma\left(1_{\theta_a(D_g)}\right)&=\gamma(1_{D_a\cap D_{ag}})\\
      &=s_as_a^*s_{ag}s_{ag}^*\\
      &=s_as_a^*s_{ag}s_{ag}^*s_as_a^*\\
      &=s_as_gs_g^*s_a^*,
    \end{split}
  \end{equation*}
  and that
  \begin{equation*}
    \begin{split}
      \gamma\bigl(1_{\theta_{a\inv}(D_g)}\bigr)
      &=\gamma\bigl(1_{D_{a\inv}\cap D_{a\inv g}}\bigr)\\
      &=s_{a\inv}s_{a\inv}^*s_{a\inv g}s_{a\inv g}^*\\
      &=s_{a\inv}s_{a\inv}^*s_{a\inv g}s_{a\inv g}^*s_{a\inv}s_{a\inv}^*\\
      &=s_{a\inv}s_gs_g^*s_{a\inv}^*.
    \end{split}
  \end{equation*}
  Since $\spc\{1_{D_g}\mid g\in\fg\al\}=\DX$, this shows that
  $\gamma\bigl(\phi_a(f)\bigr)=s_a\gamma(f)s_a^*$ and
  $\gamma\bigl(\lambda_a(f)\bigr)=s_a^*\gamma(f)s_a$ for $f\in\DX$ and
  $a\in\al$. Thus we have that
  \begin{equation*}
    \begin{split}
      \gamma\bigl(1_{\cy\mu\nu}\bigr)
      &=\gamma\bigl(1_{\theta_\nu(D_{\mu\inv})}\bigr)\\
      &= \gamma\bigl(\phi_{\nu_1}\circ\phi_{\nu_2}\dotsb\circ
      \phi_{\mu_k}\circ\lambda_{\mu_l}\circ\dotsb\circ\lambda_{\mu_1}(1)\bigr)\\
      &= s_{\nu_1}s_{\nu_2}\dotsm s_{\nu_l}s_{\mu_l}^*\dotsm
      s_{\mu_1}^*s_{\mu_1}\dotsm s_{\mu_l}s_{\nu_k}^*\dotsm
      s_{\nu_1}^*
    \end{split}    
  \end{equation*}
  for $\mu=\mu_1\mu_2\dotsm \mu_l$ and $\nu=\nu_1\nu_2\dotsm \nu_k$ in
  $\al^*$. Hence it follows from the universal property of $\O_{\OSS}$
  that there is a $*$-homomorphism $\tilde{\gamma}$ from
  $\O_{\OSS}$ to $\csp{\OSS,\theta_{\OSS},\fg{\al}}$ such that
  $\tilde{\gamma}\circ\ioss=\gamma$ (and hence
  $\tilde{\gamma}\bigl(\ioss(1_{D_g})\bigr)=s_gs_g^*$ for every $g\in \fg\al$)
  and $\tilde{\gamma}(S_a)=s_a$ for every $a\in\al$.
  
  Let us now look at $\O_{\OSS}$. Let $S_{a\inv}=S_a^*$ for every
  $a\in\al$, let $S_e=1$ and let for $g\in \fg\al$ written in the
  reduced form $b_1b_2\dotsm b_k$, $S_g=S_{b_1}S_{b_2}\dotsm S_{b_k}$.
  Then $(S_g)_{g\in\fg{\al}}$ is a family of elements from $\O_{\OSS}$ which
  clearly satisfies \eqref{eq:oxa} and \eqref{eq:oxb} of Definition
  \ref{def:ox}. We will show that $(S_g)_{g\in\fg{\al}}$ also
  satisfies \eqref{eq:ox}, \eqref{eq:oxc} and
  \eqref{eq:oxd} of Definition \ref{def:ox}.
  
  It is easy to check that $\ioss\bigl(\lambda_a(f)\bigr)=S_a^*\ioss(f)S_a$
  and $\ioss\bigl(\phi_a(f)\bigr)=S_a\ioss(f)S_{a\inv}$ for $a\in\al$ and
  $f\in \Dx$. Thus if we for $b\in\al\cup\al\inv$ define $\omega_b$ by
  \begin{equation*}
    \omega_b=
    \begin{cases}
      \phi_b&\text{if }b\in\al,\\
      \lambda_{b\inv}&\text{if }b\inv\in\al,
    \end{cases}
  \end{equation*}
  then we have that
  \begin{equation*}
    \begin{split}
      \ioss(1_{D_g})
      &=\ioss\Bigl(1_{\theta_{b_1}\circ\theta_{b_2}\circ\dotsb
        \circ\theta_{b_k}(\OSS)}\Bigr)\\
      &=\ioss\bigl(\omega_{b_1}\circ\omega_{b_2}\circ
      \dotsb\circ\omega_{b_k}(1)\bigr)\\
      &=S_{b_1}S_{b_2}\dotsm S_{b_k}S_{b_k}^*\dotsm S_{b_1}^*\\
      &=S_gS_g^*
    \end{split}
  \end{equation*}
  for every $g\in\fg{\al}$ written in the reduced form $b_1b_2\dotsm b_k$.
  This shows that $(S_g)_{g\in\fg{\al}}$ satisfies \eqref{eq:ox} and
  \eqref{eq:oxd} of Definition \ref{def:ox}.
  
  Let $h,i\in \fg\al$, and let $b_1b_2\dotsm b_k$ be the reduced
  form of $h$ and $b_1'b_2'\dotsm b_{k'}'$ be the reduced form of
  $i$. Consider those $l\in\{1,2,\dotsc,k\}$ for which
  $k+1-l\le k'$ and $b_1'=b_k\inv,
  b_2'=b_{k-1}\inv,\dotsc,b_{k+1-l}'=b_l\inv$, and for which $b_{k+2-l}'\ne
  b_{l-1}'$ if $l\ne 1$ and $k+1-l\ne k'$. Notice that if
  such an $l$ exists, then it is necessarily unique. In this case, we
  let $j=b_lb_{l+1}\dotsm b_k$. If no such $l$ exists,
  then we let $j$ be equal to the neutral element.
  
  Let $h'=hj\inv$ and $i'=ji$. We then have that the reduced form of
  $h'$ is $b_1b_2\dotsm b_{l-1}$, the reduced form of $i'$ is
  $b_{k+2-l}'b_{k+3-l}'\dotsm b_{k'}'$ and the reduced form of
  $hi=h'i'$ is $b_1b_2\dotsm b_{l-1}b_{k+2-l}'b_{k+3-l}'\dotsm
  b_{k'}'$ (if $j=e$, then $h'=h$, $i'=i$ and the reduced form of
  $hi=h'i'$ is $b_1b_2\dotsm b_kb_1'b_2'\dotsm b_{k'}'$). Thus
  $S_{h'}S_{i'}=S_{h'i'}=S_{hi}$, and we have that
  \begin{equation*}
    \begin{split}
      S_hS_i&= S_{h'}S_jS_j^*S_{i'}\\
      &= S_{h'}S_{h'}^*S_{h'}S_jS_j^*S_{i'}\\
      &= S_{h'}S_jS_j^*S_{h'}^*S_{h'}S_{i'}\\
      &= S_hS_h^*S_{h'i'}\\
      &= S_hS_h^*S_{hi}.
    \end{split}
  \end{equation*}
  This shows that $(S_g)_{g\in\fg{\al}}$ satisfies \eqref{eq:oxc}
  of Definition \ref{def:ox}.
  
  Hence it follows from the universal property of
  $\csp{\OSS,\theta_{\OSS},\fg{\al}}$ that there is a $*$-homomorphism
  $\psi$ from $\csp{\OSS,\theta_{\OSS},\fg{\al}}$ to $\O_{\OSS}$ such that
  $\psi(s_g)= S_{b_1}S_{b_2}\dotsm S_{b_k}$ for every $g\in \fg\al$,
  where $g\in\fg\al$ is written in the reduced form $b_1b_2\dotsm b_k$.
  
  We have that $\psi\bigl(\tilde{\gamma}(S_a)\bigr)=\psi(s_a)=S_a$ for every
  $a\in\al$, and since $\O_{\OSS}$ is generated by
  $(S_a)_{a\in\al}$, this shows that
  $\psi\circ\tilde{\gamma}=\Id_{\O_{\OSS}}$.
  
  We also have by Lemma \ref{lemma:onestruc} that
  \begin{equation*}
    \tilde{\gamma}\bigl(\psi(s_g)\bigr)= 
    \tilde{\gamma}\bigl(S_{b_1}S_{b_2}\dotsm S_{b_k}\bigr)= 
    s_{b_1}s_{b_2}\dotsm s_{b_k}=s_g
  \end{equation*}
  for every $g\in \fg\al$,
  where $g\in \fg\al$ is written in the reduced form $b_1b_2\dotsm b_k$,
  and since $\csp{\OSS,\theta_{\OSS},\fg{\al}}$ is generated by
  $(s_g)_{g\in\fg\al}$, this shows that
  $\tilde{\gamma}\circ\psi=\Id_{\csp{\OSS,\theta_{\OSS},\fg{\al}}}$.
    
  Thus $\tilde{\gamma}$ is a $*$-isomorphism from $\O_{\OSS}$ to
  $\csp{\OSS,\theta_{\OSS},\fg{\al}}$ which maps $\ioss(1_{D_g})$ to
  $s_gs_g^*$ for every $g\in \fg\al$ and $S_a$ to $s_a$ for every
  $a\in\al$..
\end{proof}

\subsection{Cuntz-Krieger algebras} \label{subsec:CK}

Let $\ind$ be an arbitrary index set and let
$A=(A(i,j))_{i,j\in\ind}$ be a matrix with entries in $\{0,1\}$ and
having no zero rows. Exel and Laca have in
\cite{MR2000i:46064} defined a $\cs$-algebra $\O_A$ associated with
$A$. It is the universal $\cs$-algebra generated by a family
$(S_i)_{i\in\ind}$ of partial isometries satisfying the following 4
conditions:     
\begin{align}
  &\label{eq:cka}\\
  &\raisebox{0pt}[38pt][47pt]{\parbox{0.8\textwidth}{for each pair of
      finite subsets $X$ and $Y$ of $\ind$ such that the number
      $A(X,Y,j)$ defined by 
      \begin{equation*}
        A(X,Y,j)=\prod_{x\in X}A(x,j)\prod_{y\in Y}(1-A(y,j))
      \end{equation*}
      is zero for all but a finite number of $j$'s in $\ind$, the
      following equation holds:
      \begin{equation*}
        \prod_{x\in X}S_x^*S_x\prod_{y\in Y}(1-S_y^*S_y) =\sum_{j\in
          \ind}A(X,Y,j)S_jS_j^*,
      \end{equation*}}}\\
  &S_i^*S_i \text{ and } S_j^*S_j \text{ commute for all
  }i,j\in\ind,\label{eq:ckb}\\
  &S_i^*S_j=0, \text{ if }i\ne j\in \ind,\label{eq:ckc}\\
  &S_i^*S_iS_j=A(i,j)S_j, \text{ for all }i,j\in \ind\label{eq:ckd}.
\end{align}
We will call the family $(S_i)_{i\in\ind}$ the generators of $\O_A$.

If $A$ is a finite matrix, then the conditions
\eqref{eq:cka}--\eqref{eq:ckd} reduce to the ordinary Cuntz-Krieger
relations, and so $\O_A$ is isomorphic to the usual
Cuntz-Krieger algebra of $A$ \cite{MR561974}, or the universal
Cuntz-Krieger algebra considered in \cite{MR1452183}
(cf. \cite{MR2000i:46064}*{Examples 8.9}). 

Let $(\XA,\msh)$ be the topological Markov chain associated to
$A$. That is: $\XA$ is the set defined by 
\begin{equation*}
  \XA=\bigl\{(x_n)_{n\in\No}\in\ind^{\No}\bigm| \forall
  k\in\No:A(x_k,x_{k+1})=1\bigr\},
\end{equation*}
and $\msh:\XA\to\XA$ is the map
\begin{equation*}
  x_0x_1x_2\dotsm\mapsto x_1x_2\dotsm .
\end{equation*}

Since $(\XA,\msh)$ is a one-sided symbolic dynamical system, we can
associate to it the discrete partial dynamical system
$(\XA,\theta_\XA,\fg{\ind})$ of Definition \ref{onesideddef}.
We then have the following theorem:

\begin{theorem} \label{theorem:CK}
  Let $\ind$ be an arbitrary index set, let
  $A=(A(i,j))_{i,j\in\ind}$ be a matrix with entries in $\{0,1\}$ and
  having no zero rows, and let $(\XA,\theta_\XA,\fg{\ind})$
  be the discrete partial dynamical system defined above. Then
  $\csp{\XA,\theta_\XA,\fg{\ind}}$ is isomorphic to the unitization
  $\widetilde{\O}_A$ of $\O_A$.
 
  More precisely: if $(s_g)_{g\in\fg\ind}$ denotes the generators of
  $\csp{\XA,\theta_\XA,\fg{\ind}}$ and $(S_i)_{i\in\ind}$ denotes the
  generators of $\O_A$, then 
  $\csp{\XA,\theta_\XA,\fg{\ind}}$ is generated by its unit and
  $\{s_i\mid i\in\ind\}$, and there exists a unital
  $*$-isomorphism from $\csp{\XA,\theta_\XA,\fg{\ind}}$ to
  $\widetilde{\O}_A$ which maps $s_i$ to $S_i$ for every $i\in\ind$.
\end{theorem}

\begin{proof}
  Let $(D_g)_{g\in\fg\ind}$ and $(\theta_g)_{g\in\fg\ind}$ denote the
  domains and partial one-to-one maps of $\theta_\XA$.

  Using the facts that $s_is_i^*s_i=s_i$ for every $i\in\ind$ and
  that the map $D_g\mapsto s_gs_g^*$ extends to a Boolean
  homomorphism from the Boolean algebra $\BA$
  to the Boolean algebra of projections in the unital abelian $\cs$-algebra
  \begin{equation*}
    \cs\Bigl(\bigl\{s_gs_g^*\in\csp{\XA,\theta_\XA,\fg{\ind}}\bigm|
    g\in\fg\ind\bigr\}\Bigr),
  \end{equation*}
  it is easy to check that the family $(s_i)_{i\in\ind}$ of elements
  of $\csp{\XA,\theta_\XA,\fg{\ind}}$ satisfies condition \eqref{eq:cka},
  \eqref{eq:ckb}, \eqref{eq:ckc} and
  \eqref{eq:ckd} above. Thus it follows from the universal property of $\O_A$
  that there is a $*$-homomorphism $\phi$ from $\O_A$ to
  $\csp{\XA,\theta_\XA,\fg{\ind}}$ such that $\phi(S_i)=s_i$ for all
  $i\in\ind$. 

  If $\O_A$ is not unital, then $\phi$ extends to a unital $*$-homomorphism
  $\tilde{\phi}$ from $\widetilde{\O}_A$ to
  $\csp{\XA,\theta_\XA,\fg{\ind}}$. If $\O_A$ is unital (in which case
  $\widetilde{\O}_A=\O_A$), then it follows from
  \cite{MR2000i:46064}*{Proposition 8.5} that there exist finite
  subsets $X$ and $Y$ of $\ind$ such that the equation
  \begin{equation*}
    \XA=\bigcup_{x\in X}D_x\cup\bigcup_{y\in Y}D_{y\inv}
  \end{equation*}
  holds.
  That means that the unit of $\csp{\XA,\theta_\XA,\fg{\ind}}$ is
  contained in the $\cs$-subalgebra generated by $\{s_i\mid
  i\in\ind\}$ and thus is in the image of $\phi$, and so $\phi$ is
  unital. 

  So we have in both cases that there exists a unital $*$-homomorphism
  $\tilde{\phi}$ from $\widetilde{\O}_A$ to
  $\csp{\XA,\theta_\XA,\fg{\ind}}$ which maps $S_i$ to $s_i$ for all $i\in\ind$.

  Let us now turn towards $\widetilde{\O}_A$. We let for every
  $i\in\ind$, $S_{i\inv}=S_i^*$ and we let $S_e=1$ and
  $S_g=S_{b_1}S_{b_2}\dotsm S_{b_k}$, where $g=b_1b_2\dotsm b_k\in
  \fg\ind$ is written in reduced form. It then follows from
  condition \eqref{eq:ckb}, \eqref{eq:ckc} and \eqref{eq:ckd} that the map
  \begin{equation*}
    g\mapsto S_g,\quad g\in \fg\ind
  \end{equation*}
  is a partial representation of $\fg\ind$ (see
  \cite{MR2000i:46064}*{Proposition 3.2} for a proof of this) and
  thus that $(S_g)_{g\in\fg\ind}$ is a family of partial isometries which
  satisfies condition \eqref{eq:ox}--\eqref{eq:oxc} of Definition
  \ref{def:ox}.
  We will now show that $(S_g)_{g\in \fg\ind}$ also satisfies
  \eqref{eq:oxd} of Definition \ref{def:ox}.
  
  It follows from Lemma \ref{lemma:onestruc} that
  $s_g=s_{b_1}s_{b_2}\dotsm s_{b_k}$ for every 
  $g\in \fg\ind$ written in the reduced form $b_1b_2\dotsm b_k$.
  Thus the $*$-homomorphism $\tilde{\phi}:\widetilde{\O}_A\to
  \csp{\XA,\theta_\XA,\fg{\ind}}$ 
  mentioned above maps $S_g$ to $s_g$ for every $g\in \fg\ind$, and
  thus induces a Boolean homomorphism from the Boolean algebra of
  projections in the unital abelian $\cs$-algebra
  \begin{equation*}
    \cs\Bigl(\bigl\{S_g S_g^*\in\widetilde{\O}_A\bigm| g\in\fg\ind\bigr\}\Bigr)
  \end{equation*}
  to the Boolean algebra of projections in the unital abelian $\cs$-algebra
  \begin{equation*}
    \cs\Bigl(\bigl\{s_gs_g^*\in\csp{\XA,\theta_\XA,\fg{\ind}}\bigm|
    g\in\fg\ind\bigr\}\Bigr)
  \end{equation*}
  which maps $S_gS_g^*$ to $s_gs_g^*$. Let us denote the Boolean
  subalgebra of the Boolean algebra of projections in the unital
  abelian $\cs$-algebra 
  \begin{equation*}
    \cs\Bigl(\bigl\{S_g S_g^*\in\widetilde{\O}_A\bigm| g\in\fg\ind\bigr\}\Bigr)
  \end{equation*}
  generated by $\{S_gS_g^*\mid g\in \fg\ind\}$ by $\BO$. Remember
  (cf. Definition 
  \ref{remarkboomap} and Corollary \ref{boomapinj}) that
  $\boomap[(\XA,\theta_\XA,\fg{\ind})]$
  is an injective Boolean homomorphism from $\BA$ to the Boolean algebra
  of projections in the unital abelian $\cs$-algebra
  \begin{equation*}
    \cs\Bigl(\bigl\{s_gs_g^*\in\csp{\XA,\theta_\XA,\fg{\ind}}\bigm|
      g\in\fg\ind\bigr\}\Bigr)
  \end{equation*} 
  which maps $D_g$ to $s_gs_g^*$. Thus there is a Boolean
  homomorphism $\eta$ from $\BO$ to $\BA$ which maps
  $S_gS_g^*$ to $D_g$. We claim that $\eta$ is injective.
  
  To prove this we will use that the family $(S_g)_{g\in\fg\ind}$ has the
  following properties (remember that we by $\ind^*$ denote the set
  of finite words with letters from $\ind$, and that we identify
  $\ind^*$ with the unital sub-semigroup of $\fg\ind$ generated by
  $\ind$, cf. Section \ref{sec:notat}):
  \begin{eqnarray}
    &&S_h^*S_hS_i^*S_i=S_i^*S_iS_h^*S_h \text{ for all
    }h,i\in\fg\ind,\label{eq:ckpa}\\
    &&S_\mu S_\nu=S_{\mu\nu} \text{ for all }\mu,\nu\in\ind^*,\label{eq:ckpb}\\
    &&S_\mu^*S_\nu=0 \text{ if }\mu,\nu\in\ind^*,\
    \abs{\mu}=\abs{\nu}\text{ and }\mu\ne\nu,\label{eq:ckpc}\\
    &&S_\mu^*S_\mu=A(\mu_1,\mu_2)\dotsm
    A(\mu_{\abs{\mu}-1},\mu_{\abs{\mu}})
    S_{\mu_{\abs{\mu}}}^*S_{\mu_{\abs{\mu}}}\text{ for all
    }\mu\in\ind^*,\label{eq:ckpd}\\
    &&\text{if }S_g\ne 0\text{ then }g=\mu\nu\inv \text{ and }S_g=S_\mu
    S_\nu^*\text{ for some
    }\mu,\nu\in\ind^*.\label{eq:ckpe} 
  \end{eqnarray}
  Property \eqref{eq:ckpa} follows from the fact that $g\mapsto S_g$
  is a partial representation of $\fg\ind$, \eqref{eq:ckpb}
  follows from the definition of $S_g$, \eqref{eq:ckpc} is
  \emph{Claim 2} and \eqref{eq:ckpd} is \emph{Claim 1} in the proof
  of \cite{MR2000i:46064}*{Proposition 3.2}. To see \eqref{eq:ckpe}
  notice that it follows from \cite{MR2000i:46064}*{Proposition 3.1}
  that if $S_g\ne 0$, then $g=\mu\nu\inv$ for some $\mu,\nu\in\ind^*$,
  and if we choose $\mu$ and $\nu$ such that the last symbol of $\mu$
  is different from the first letter of $\nu$ (or $\mu$ or $\nu$ is
  the empty word), then $S_g=S_\mu S_\nu^*$ by definition.    
  
  Let for every $\mu\in\ind^*$ and every pair
  $(I,J)$ of (possible empty) finite subsets of $\ind$, $C(\mu,I,J)$
  be the subset of $\XA$ defined by
  \begin{equation*}
    \begin{split}
      C(\mu,I,J)&=\theta_\mu\biggl(\Bigl(\bigcap_{i\in
        I}D_{i\inv}\Bigr)\cap\Bigl(\bigcap_{j\in J}\XA\setminus
      D_{j\inv}\Bigr)\biggr)\\
      &= \bigl\{\mu x\in \XA\bigm| \forall i\in I: ix\in\XA, \forall
      j\in J:jx\notin \XA\bigr\}.
    \end{split}
  \end{equation*}
  Notice that we have that 
  \begin{equation*}
    \eta\Biggl(S_u\prod_{i\in I}S_i^*S_i\prod_{j\in
        J}(1-S_j^*S_j)S_u^*\Biggr)= C(\mu,I,J)
  \end{equation*}
  for finite subsets $I,J$ of $\ind$ and $\mu\in\ind^*$.

  We will now prove that every element of $\BA$ can be written as a
  finite union of elements of the form 
  \begin{equation*}
    S_\mu\prod_{i\in I}S_i^*S_i \prod_{j\in J}(1-S_j^*S_j)S_\mu^* \left(
      \prod_{k=1}^n \Biggl(1-S_{\mu^k} \prod_{i\in I_k}S_i^*S_i
        \prod_{j\in J_k}(1-S_j^*S_j)S_{\mu^k}^* \Biggr) \right)
  \end{equation*}
  where $\mu,\mu^1,\mu^2,\dots \mu^n\in \ind^*$ and
  $I,J,I_1,J_1,I_2, J_2, \dotsc ,I_n,J_n$ all are finite (possible
  empty) subsets of $\ind$.

  Notice first that if $\mu\in\ind^*$ and $i\in\ind$, then we have that
  \begin{equation*}
    S_\mu S_i^*S_i S_\mu^*=
    \begin{cases}
      S_{i\inv} S_{i\inv}^*&\text{if }u=\emptyword,\\
      S_{\mu i\inv}S_{\mu i\inv}^*&\text{if }\mu_{|\mu|}\ne i,\\
      S_\mu S_\mu^*&\text{if }\mu_{|\mu|}=i,
    \end{cases}
  \end{equation*}
  and thus that $S_\mu S_i^*S_i S_\mu^*$ belongs to $\BO$. Since $\BO$
  is closed under intersection and complement, and the complement of
  an element $p$ of $\BO$ is defined to be $1-p$ and the intersection
  of two elements $p,q\in\BO$ is $pq$, we also have that 
  \begin{equation*}
    S_\mu(1-S_j^*S_j)S_\mu^*=S_\mu S_\mu^*(1-S_\mu S_j^*S_j S_\mu^*)\in\BO
  \end{equation*}
  for every $\mu\in\ind^*$ and every $j\in\ind$.

  Let $\mu,\tilde{\mu}\in\ind^*$ and let $I,J,\tilde{I},\tilde{J}$
  be subsets of $I$. It then follows from \eqref{eq:ckd} and
  \eqref{eq:ckpa}--\eqref{eq:ckpc} that the element
  \begin{equation*}
    S_\mu \prod_{i\in I}S_i^*S_i \prod_{j\in J}
    \left(1-S_j^*S_j\right) S_\mu^*
    S_{\tilde{\mu}}\prod_{i\in\tilde{I}}S_i^*S_i \prod_{j\in
      \tilde{J}}\left(1-S_j^*S_j\right) S_{\tilde{\mu}}^*
  \end{equation*}
  is equal to 
  \begin{equation*}
    S_\mu\smashoperator{\prod_{i\in I}}
    S_i^*S_i\smashoperator{\prod_{j\in J}}  
    \left(1-S_j^*S_j\right) S_\mu^*
  \end{equation*} 
  if $\mu=\tilde{\mu}\nu$ for some $\nu=\nu_1\nu_2\dotsm \nu_{|\nu|}
  \in\ind^*$ which satisfies that $A(\tilde{I},\tilde{J},\nu_1)=1$, to
  \begin{equation*}
    S_{\tilde{\mu}}\smashoperator{\prod_{i\in\tilde{I}}}S_i^*S_i
    \smashoperator{\prod_{j\in \tilde{J}}}\left(1-S_j^*S_j\right)
    S_{\tilde{\mu}}^*
  \end{equation*}
  if $\tilde{\mu}=\mu\nu$ for some
  $\nu=\nu_1\nu_2\dotsm \nu_{|\nu|} \in\ind^*$ which satisfies that
  $A(I,J,\nu_1)=1$, to 
  \begin{equation*}
    S_\mu \smashoperator{\prod_{i\in I\cup\tilde{I}}}
    S_i^*S_i \smashoperator{\prod_{j\in J\cup\tilde{J}}} 
    \left(1-S_j^*S_j\right) S_\mu^*
  \end{equation*} 
  if $\mu=\tilde{\mu}$, and to 0 otherwise. 

  Thus if we let $\mathcal{Z}$ be the set
  \begin{equation*}
    \Biggl\{S_\mu\prod_{i\in I}S_i^*S_i\prod_{j\in J}
    \left(1-S_j^*S_j\right) S_\mu^*\Bigm| \mu\in\ind^*,\ I,J
    \text{ are finite subset of }\ind^*\Biggr\},
  \end{equation*}
  then $\mathcal{Z}$ is a subset of $\BO$, and it is closed under intersection. 
  So if we let  
  \begin{equation*}
    \mathcal{Y}=\left\{Z_0\bigcap_{Z\in \find}(1-Z)\Bigm| Z_0\in
    \mathcal{Z},\ \find\text{ is a finite subset of }\mathcal{Z} \right\},
  \end{equation*}
  then $\mathcal{Y}$ is also closed under intersection.
  Let $\mathcal{X}$ be the set of
  elements in $\BO$ which can be written as a finite union of 
  elements from $\mathcal{Y}$. Then $\mathcal{X}$ is closed under
  union and intersection, and since we also have that $1-Y$ belongs to
  $\mathcal{X}$ if $Y$ belongs to $\mathcal{Y}$, we also have that
  $\mathcal{X}$ is closed under complement and thus is a Boolean
  subalgebra of $\BO$.
  It follows from \eqref{eq:ckpb}, \eqref{eq:ckpd} and \eqref{eq:ckpe}
  that if $g\in\fg\ind$ and $S_g\ne 0$, then $g=\mu\nu\inv$ for
  suitable $\mu,\nu\in\ind^*$ and that
  \begin{equation*}
    S_g^*S_g=S_\nu S_\mu^*S_\mu S_\nu^*=S_\nu S_i^*S_iS_\nu^*\in\mathcal{X}
  \end{equation*}
  where $i$ is the last letter of $\mu$ (if $\mu$ is the empty word,
  then $S_g^*S_g=S_\nu S_\nu^*\in \mathcal{X}$).
  Thus every element in $\BO$ belongs to $\mathcal{X}$, and is thus a finite
  union of elements of the form
  \begin{equation*}
    S_\mu\prod_{i\in I}S_i^*S_i \prod_{j\in J}(1-S_j^*S_j)S_\mu^* \left(
      \prod_{k=1}^n \Biggl(1-S_{\mu^k} \prod_{i\in I_k}S_i^*S_i
        \prod_{j\in J_k}(1-S_j^*S_j)S_{\mu^k}^* \Biggr) \right)
  \end{equation*}
  where $\mu,\mu^1,\mu^2,\dots \mu^n\in \ind^*$ and
  $I,J,I_1,J_1,I_2, J_2, \dotsc ,I_n,J_n$ all are finite (possible
  empty) subsets of $\ind$.  

  In order to prove that $\eta$ is injective, it is therefore enough
  to show that if $\mu,\mu^1,\mu^2,\dots \mu^n\in \ind^*$ and
  $I,J,I_1,J_1,I_2, J_2, \dotsc ,I_n,J_n$ all are finite (possible
  empty) subsets of $\ind$ and 
  \begin{equation*}
    C(\mu,I,J)\cap\left(\bigcap_{k=1}^n \XA\setminus
      C(\mu^k,I_k,J_k)\right)=\emptyset, 
  \end{equation*}
  then we have that
  \begin{equation*}
    S_\mu\prod_{i\in I}S_i^*S_i\prod_{j\in J}(1-S_j^*S_j)S_\mu^*
    \left(\prod_{k=1}^n \Biggl(1-S_{\mu^k}\prod_{i\in I_k}
      S_i^*S_i\prod_{j\in J_k}(1-S_j^*S_j)S_{\mu^k}^*\Biggr)\right)=0.
  \end{equation*}
  
  Let $\mu,\nu\in\ind^*$ with $|\mu|>|\nu|$ (remember that $|\mu|$
  and $|\nu|$ denote the length of $\mu$ and $\nu$, respectively, cf. Section
  \ref{sec:notat}) and let $I,J,I',J'$ be finite subsets of
  $\ind$. Then either the equality 
  \begin{equation*}
    C(\mu,I,J)\cap \XA\setminus C(\nu,I',J')=C(\mu,I,J)
  \end{equation*}
  holds, or $\mu=\nu\alpha$ for some $\alpha=\alpha_1\alpha_2\dotsm
  \alpha_{|\mu|-|\nu|}\in\ind^*$ which satisfies that $A(I',J',\alpha_1)=1$. 
  In the latter case, we have that
  \begin{equation*}
    S_\mu S_\mu^*S_\nu\prod_{i\in I_k}S_i^*S_i
    \prod_{j\in J_k}(1-S_j^*S_j)S_\nu^* 
    = S_\mu S_\mu^*, 
  \end{equation*}
  and thus that
  \begin{equation*}
    S_\mu\prod_{i\in I}S_i^*S_i\prod_{j\in J}(1-S_j^*S_j)S_\mu^*
    \Biggl(1-S_\nu\prod_{i\in I'}S_i^*S_i\prod_{j\in J'}
    (1-S_j^*S_j)S_\nu^*\Biggr)=0.
  \end{equation*}
  So in order to prove that the condition
  \begin{equation*}
    C(\mu,I,J)\cap\left(\bigcap_{k=1}^n \XA\setminus
      C(\mu^k,I_k,J_k)\right)=\emptyset,
  \end{equation*}
  implies that
  \begin{equation*}
    S_\mu\prod_{i\in I}S_i^*S_i\prod_{j\in J}(1-S_j^*S_j)S_\mu^*
    \left(\prod_{k=1}^n \Biggl(1-S_{\mu^k}\prod_{i\in
          I_k}S_i^*S_i\prod_{j\in
          J_k}(1-S_j^*S_j)S_{\mu^k}^*\Biggr)\right)=0,
  \end{equation*}
  we may assume that $|\mu^k|\ge|\mu|$ for every $k\in\{1,2,\dotsc ,n\}$.
 
  So let $\mu,\mu^1,\mu^2,\dots \mu^n\in \ind^*$ with
  $|\mu^k|\ge|\mu|$ for every $k\in\{1,2,\dotsc ,n\}$ and let
  $I,J,I_1,J_1,I_2, J_2, \dotsc ,I_n,J_n$ be finite subsets of $\ind$
  such that
  \begin{equation*}
    C(\mu,I,J)\cap\left(\bigcap_{k=1}^n \XA\setminus
      C(\mu^k,I_k,J_k)\right)=\emptyset.
  \end{equation*}
  
  \begin{clai}
    We claim that there for every $m>|\mu|$ exists a finite subset
    $\fw_m$ of the set
    \begin{multline*}
      \biggl\{\nu\in\ind^*\Bigm| |\nu|=m,\ \forall x\in
      \XA:A(\nu_m,x_0)=1\Longrightarrow \\
      \nu x\in C(\mu,I,J)\cap
      \Bigl(\smashoperator[r]{\bigcap_{\substack{k=1\\
            |\mu^k|<m}}^n} 
      \XA\setminus C(\mu^k,I_k,J_k)\Bigr)\biggr\}
    \end{multline*}
    such that the following equality holds:
    \begin{multline*}
      S_\mu\prod_{i\in I}S_i^*S_i\prod_{j\in J}(1-S_j^*S_j)S_\mu^*
      \Biggl(\prod_{k=1}^n \biggl(1-S_{\mu^k}\prod_{i\in I_k}
      S_i^*S_i\prod_{j\in J_k}(1-S_j^*S_j)S_{\mu^k}^*\biggr)\Biggr)\\
      =\sum_{\nu\in\fw_m}S_\nu S_\nu^* \smashoperator[r]{
        \prod_{\substack{k=1\\
            |\mu^k|\ge m}}^n}
      \biggl(1-S_{\mu^k} \prod_{i\in I_k}S_i^*S_i\prod_{j\in J_k}
      (1-S_j^*S_j)S_{\mu^k}^*\biggr).
    \end{multline*}   
  \end{clai}
  
  It follows from this claim that the equality
  \begin{equation*}
    S_\mu\prod_{i\in I}S_i^*S_i\prod_{j\in J}(1-S_j^*S_j)S_\mu^*
    \left(\prod_{k=1}^n \Biggl(1-S_{\mu^k}\prod_{i\in
        I_k}S_i^*S_i\prod_{j\in
        J_k}(1-S_j^*S_j)S_{\mu^k}^*\Biggr)\right)=0,
  \end{equation*}
  holds, because since $A$ has no identically zero row, there is for every
  $\nu=\nu_1\nu_2\dotsm \nu_m\in\ind^*$, an $x\in\XA$ such that
  $A(\nu_m,x_0)=1$, and it follows from this that the set
  \begin{multline*}
    \biggl\{\nu\in\ind^*\Bigm| |\nu|=m,\ \forall x\in
    \XA:A(\nu_m,x_0)=1\Longrightarrow \\
    \nu x\in C(\mu,I,J)\cap
    \Bigl(\smashoperator[r]{\bigcap_{\substack{k=1\\ |\mu^k|< m}}^n} 
    \XA\setminus C(\mu^k,I_k,J_k)\Bigr)\biggr\}
  \end{multline*}
  is empty for $m>\max\{|\mu^k|\mid k\in \{1,2,\dotsc ,n\}\}$.
  
  To prove the claim, we need a little lemma:

  \begin{lemma} \label{modsat}
    Let $R$ be a unital ring, let $k\in\N$ and let 
    $x_1,x_2,\dotsc x_n,\in R$. Then we have that
    \begin{equation*}
      1-\prod_{i=1}^nx_i=\sum_{E\in\cset_n}
      \smashoperator[r]{\prod_{i\in\{1,2,\dotsc n\}\setminus E}}
      x_i\hspace{10pt}\prod_{j\in E}(1-x_j)
    \end{equation*}
    where $\cset_n$ is the set of non-empty subsets of 
    $\{1,2,\dots ,n\}$.
  \end{lemma}
  
  \begin{proof}
    We will prove the lemma by induction. The lemma obviously holds
    for $n=1$. Assume now that the lemma holds for $n=m$. Then we have
    that 
    \begin{equation*}
      \begin{split}
        1-\prod_{i=1}^{m+1}x_i=& \left(1-\prod_{i=1}^mx_i\right)x_{m+1} +
        \prod_{i=1}^mx_i(1-x_{m+1}) +
        \left(1-\prod_{i=1}^mx_i\right)(1-x_{m+1}) \\
        =& \begin{multlined}[t]
          \sum_{E\in\cset_m}\smashoperator[r]{\prod_{i\in\{1,2,\dotsc
              m\}\setminus E}}x_i\hspace{11pt}\prod_{j\in E}(1-x_j)x_{m+1} +
          \prod_{i=1}^mx_i(1-x_{m+1})\\ +
          \sum_{E\in\cset_m}\smashoperator[r]{\prod_{i\in\{1,2,\dotsc
            m\}\setminus E}}x_i\hspace{11pt}\prod_{j\in
          E}(1-x_j)(1-x_{m+1}) 
        \end{multlined}\\
        =& \sum_{E\in\cset_{m+1}}\smashoperator[r]{\prod_{i\in\{1,2,\dotsc
          m+1\}\setminus E}}x_i\hspace{16pt}\prod_{j\in E}(1-x_j).
      \end{split}
    \end{equation*}
  \end{proof}

  \begin{proof}[Proof of the claim.]
    We will now prove the claim by induction. First let
    $m=|\mu|+1$. Let $K_\mu$ be the subset of $\{1,2,\dotsc ,n\}$
    defined by
    \begin{equation*}
      K_\mu=\bigl\{k\in \{1,2,\dotsc ,n\}\bigm| \mu^k=\mu\bigr\},
    \end{equation*}
    and let for each $k\in K_\mu$, $\cset_k$ be the set of non-empty subsets of
    $I_k\cup J_k$. We let $\pnes$ be the set of families $(E_k)_{k\in
      K_\mu}$ where for each $k\in K_\mu$, $E_k\in\cset_k$, and we 
    let for $E=(E_k)_{k\in K_\mu}\in \pnes$, $I_E$ be the finite
    subset of $\ind$ defined by
    \begin{equation*}
      I_E=I\cup\{\mu_{|\mu|}\}\cup\bigcup_{k\in K_\mu}(J_k\cap
      E_k)\cup (I_k\setminus E_k),
    \end{equation*}
    and we let $J_E$ be the finite subset of $\ind$ defined by
    \begin{equation*}
      J_E=J\cup\bigcup_{k\in K_\mu}(I_k\cap E_k)\cup 
      (J_k\setminus E_k).
    \end{equation*}
    We then have that if $y\in \XA$ and $A(I_E,J_E,y_0)=1$, then
    \begin{equation*}
      \mu y \in C(\mu,I,J)\cap
      \Bigg(\smashoperator[r]{\bigcap_{\substack{k=1\\
            |\mu^k|=|\mu|}}^n} 
      \XA\setminus C(\mu^k,I_k,J_k)\Bigg)\subseteq
      \smashoperator{\bigcup_{\substack{k=1\\ |\mu^k|>|\mu|}}^n}
      C(\mu^k,I_k,J_k), 
    \end{equation*}
    so if we let
    $\fw_m=\Bigl\{\mu j\Bigm| j\in\ind,\ \exists E\in
    \smashoperator{\prod_{k\in B_\mu}} \cset_k:A(I_E,J_E,j)=
    1\Bigr\}$, then $\fw_m$ is a subset of the set 
    \begin{multline*}
      \biggl\{\nu\in\ind^*\Bigm| |\nu|=m,\ \forall x\in
      \XA:A(\nu_m,x_0)=1\Longrightarrow \\
      \nu x\in C(\mu,I,J)\cap
      \Bigl(\smashoperator[r]{\bigcap_{\substack{k=1\\
            |\mu^k|< m}}^n} \XA\setminus C(\mu^k,I_k,J_k)\Bigr)\biggr\}.
    \end{multline*}
    Since $A$ has no zero row, there is for every $j\in \ind$ an $y\in
    \XA$ such that $y_0=j$. Now if $\mu j\in \fw_m$, then we have, as
    mentioned above, that
    \begin{equation*}
      \mu y\in \smashoperator{\bigcup_{\substack{k=1\\
            |\mu^k|>|\mu|}}^n} C(\mu^k,I_k,J_k),
    \end{equation*}
    which means that 
    \begin{equation*}
      j=y_0\in \bigl\{(\mu^k)_m\bigm| k\in \{1,2,\dotsc,n\},\
      |\mu^k|>|\mu|\bigr\}.
    \end{equation*}
    So $\fw_m$ is finite.

    If $|\mu^k|=|\mu|$, but $\mu^k\notin K_\mu$, then
    $S_\mu^*S_{\mu^k}=0$ by \eqref{eq:ckpc}, and so we have that
    \begin{multline*}
      S_\mu\prod_{i\in I}S_i^*S_i\prod_{j\in J}(1-S_j^*S_j)S_\mu^*
      \biggl(1-S_{\mu^k}\prod_{i\in I_k}S_i^*S_i\prod_{j\in J_k}
      (1-S_j^*S_j)S_{\mu^k}^*\biggr) \\
      = S_\mu\prod_{i\in I}S_i^*S_i\prod_{j\in J}(1-S_j^*S_j)S_\mu^*.
    \end{multline*}
    Thus we have that
    \begin{equation*}
      \begin{split}
        \begin{multlined}
          \mathmakebox[9.82cm][l]{S_\mu \prod_{i\in I} S_i^*S_i \prod_{j\in J}
            (1-S_j^*S_j)S_\mu^*}\\ 
          \left( \smashoperator[r]{\prod_{\substack{k=1\\
                  |\mu^k|=|\mu|}}^n} 
            \Biggl( 1-S_{\mu^k} \prod_{i\in I_k} S_i^*S_i \prod_{j\in
              J_k} 
            (1-S_j^*S_j) S_{\mu^k}^* \Biggr) \right)
          =\end{multlined}&\\
        \begin{multlined}        
          \mathmakebox[9.82cm][l]{S_\mu \prod_{i\in I} S_i^*S_i
            \prod_{j\in J} (1-S_j^*S_j) S_\mu^*}
          \\\left( \smashoperator[r]{\prod_{k\in K_\mu}}
            \Biggl(1-S_{\mu^k} \prod_{i\in I_k}
            S_i^*S_i \prod_{j\in J_k} (1-S_j^*S_j) S_{\mu^k}^* \Biggr)
          \right)
          =\end{multlined}&\\
        \begin{multlined}
          \mathmakebox[9.82cm][l]{S_\mu \prod_{i\in I} S_i^*S_i
            \prod_{j\in J} (1-S_j^*S_j) S_\mu^*S_\mu} 
          \\\left( \smashoperator[r]{\prod_{k\in K_\mu}}
            \Biggl(1-\prod_{i\in I_k} S_i^*S_i \prod_{j\in J_k}
            (1-S_j^*S_j) \Biggr) \right) S_\mu^*
          =\end{multlined}&\\
        \begin{multlined}
          \mathmakebox[9.82cm][l]{S_\mu \prod_{i\in I} S_i^*S_i
            \prod_{j\in J} (1-S_j^*S_j) S_\mu^*S_\mu} 
          \\\left( \smashoperator[r]{\prod_{k\in K_\mu}}
            \Biggl( \smashoperator[r]{\sum_{E_k\in \cset_k}} \biggl(
            \smashoperator[r]{\prod_{\substack{i\in (J_k\cap E_k)\\ \cup
                  (I_k\setminus E_k)}}} S_i^*S_i \hspace{-0.7em}
            \smashoperator[r]{\prod_{\substack{j\in (I_k\cap E_k)\\ \cup
                  (J_k\setminus E_k)}}} (1-S_j^*S_j) \biggr) \Biggr) \right)
          S_\mu^*
          =\end{multlined}&\\
        \begin{multlined}
          \mathmakebox[9.82cm][l]{S_\mu \prod_{i\in I} S_i^*S_i
            \prod_{j\in J}(1-S_j^*S_j) S_{\mu_{|\mu|}}^*S_{\mu_{|\mu|}}} \\
          \left( \smashoperator[r]{\prod_{k\in K_\mu}} \Biggl(
            \sum_{E_k\in \cset_k} \biggl(
            \smashoperator[r]{\prod_{\substack{i\in (J_k\cap E_k)\\ \cup
                  (I_k\setminus E_k)}}} S_i^*S_i \hspace{-0.7em}
            \smashoperator[r]{\prod_{\substack{j\in (I_k\cap E_k)\\ \cup
                  (J_k\setminus E_k)}}}(1-S_j^*S_j) \biggr) \Biggr) \right)
        S_\mu^*
        =\end{multlined}&\\
      \begin{multlined}
        \mathmakebox[9.82cm][l]{S_\mu\sum_{E\in \pnes} 
          \left(\prod_{i\in I} S_i^*S_i \prod_{j\in J} (1-S_j^*S_j)
            S_{\mu_{|\mu|}}^*S_{\mu_{|\mu|}}\right)}\\
        \left(
          \prod_{k\in K_\mu}
          \smashoperator[r]{\prod_{
              \substack{i\in (J_k\cap E_k)\\ 
                \cup(I_k\setminus E_k)}}}
          S_i^*S_i \hspace{-0.7em}
          \smashoperator[r]{\prod_{
              \substack{j\in (I_k\cap E_k)\\ 
                \cup(J_k\setminus E_k)}}}
          (1-S_j^*S_j) \right) S_\mu^*
        =\end{multlined}&\\
      S_\mu \sum_{E\in \pnes} 
      \left(\smashoperator[r]{\prod_{i\in I_E}} S_i^*S_i
        \smashoperator[r]{\prod_{j\in J_E}}(1-S_j^*S_j) \right) S_\mu^*
      =S_\mu\sum_{E\in \pnes}
      \smashoperator{\sum_{\substack{j\in\ind \\
            A(I_E,J_E,j) = 1}}}S_jS_j^*S_\mu^*=&
      \smashoperator{\sum_{\nu\in \fw_m}}S_\nu S_\nu^*,
      \end{split}
    \end{equation*}
    where the second equality follows from \eqref{eq:ckpa} and the
    facts that $\mu^k=\mu$ for $k\in K_\mu$ and that $S_\mu$ is a
    partial isometry, the third follows from Lemma \ref{modsat}, the
    fourth from the \eqref{eq:ckpd} (if $A(\mu_i,\mu_{i+1})=0$ for some
    $i=\{1,2,\dotsc,|\mu|-1\}$, then $S_\mu=0$ according to
    (\ref{eq:ckd}) and (\ref{eq:ckpb}), and the equality still
    holds), the fifth from the distribute law,
    the sixth from the definition of $I_E$ and $J_E$,
    the seventh from \eqref{eq:ckb} (that $ A(I_E,J_E,j)$ vanishes for
    all but a finite number of $j$'s in $\ind$ follows from the fact
    that $\fw_m$ is finite), and the eight from
    \eqref{eq:ckpb} and the definition of $\fw_m$. Thus we have
    proved the claim in the case where $m=|\mu|+1$.
    
    Assume now that $m>|\mu|$ and that
    \begin{multline*}
      S_\mu\prod_{i\in I}S_i^*S_i\prod_{j\in J}(1-S_j^*S_j)S_\mu^*
      \left(\prod_{k=1}^n \Biggl(1-S_{\mu^k}
        \prod_{i\in I_k}S_i^*S_i
        \prod_{j\in J_k}(1-S_j^*S_j)S_{\mu^k}^*\Biggr)\right)\\
      =\smashoperator{\sum_{\nu\in \fw_m}}S_\nu S_\nu^* 
      \smashoperator{\prod_{\substack{k=1\\ |\mu^k|\ge m}}^n}
      \Biggl(1-S_{\mu^k} \prod_{i\in I_k}
      S_i^*S_i\prod_{j\in J_k}(1-S_j^*S_j)S_{\mu^k}^*\Biggr)
    \end{multline*}
    for some finite subset $\fw_m$ of the set
   \begin{multline*}
     \Biggl\{\nu\in\ind^*\Bigm| |\nu|=m,\ \forall x\in
     \XA:A(\nu_m,x_0)=1\Longrightarrow \\
     \nu x\in C(\mu,I,J)\cap
     \Bigl(\smashoperator[r]{\bigcap_{\substack{k=1\\ |\mu^k|<m}}^n} 
     \XA\setminus C(\mu^k,I_k,J_k)\Bigr)\Biggr\}.
   \end{multline*}
   Let for every $\gamma\in \fw_m$, $K_\gamma$ be the subset of
   $\{1,2,\dotsc ,n\}$ defined by
   \begin{equation*}
     K_\gamma=\bigl\{k\in \{1,2,\dotsc ,n\}\bigm| \mu^k=\gamma\bigr\},
   \end{equation*}
   and let for
   each $k\in K_\gamma$, $\cset_k$ be the set of non-empty subset of
   $I_k\cup J_k$.  We 
   let $\pnes$ be the set of families $(E_k)_{k\in K_\gamma}$ where
   $E_k\in\cset_k$ for
   each $k\in K_\gamma$, and we 
   let for $E=(E_k)_{k\in K_\gamma}\in \pnes$, $I_E$ be the finite
   subset of $\ind$ defined by
   \begin{equation*} 
     I_E=\{\gamma_m\}\cup\bigcup_{k\in K_\gamma}(J_k\cap E_k)\cup
     (I_k\setminus E_k)
   \end{equation*}
   and $J_E$ be the finite subset of $\ind$ defined by
   \begin{equation*}
     J_E=\bigcup_{k\in K_\gamma}(I_k\cap E_k)\cup (J_k\setminus
     E_k).
   \end{equation*}
   We then have that if $y\in \XA$ and $A(I_E,J_E,y_0)=1$,
   then
   \begin{equation*}
     \gamma y \in \smashoperator{\bigcap_{\substack{k=1\\
         |\mu^k|= m}}^n} \XA\setminus C(\mu^k,I_k,J_k),
   \end{equation*}
   and since we have that
    \begin{multline*}
      \gamma\in \fw_m \subseteq\Biggl\{\nu\in\ind^*\Bigm| |\nu|=m,\ \forall
      x\in \XA:A(\nu_m,x_0)=1\Longrightarrow \\
      \nu x\in C(\mu,I,J)\cap
      \Bigl(\smashoperator[r]{\bigcap_{\substack{k=1\\ |\mu^k|< 
          m}}^n} \XA\setminus C(\mu^k,I_k,J_k)\Bigr)\Biggr\},
    \end{multline*}
    it follows that
    \begin{equation*}
      \gamma y\in
      C(\mu,I,J)\cap\Bigg(\smashoperator[r]{\bigcap_{\substack{k=1\\
            |\mu^k|\le m}}^n} 
      \XA\setminus C(\mu^k,I_k,J_k)\Bigg)\subseteq
      \smashoperator{\bigcup_{\substack{k=1\\ |\mu^k|>m}}^n} C(\mu^k,I_k,J_k).
    \end{equation*}
    Let $\fw_\gamma$ be the subset of $\ind$ defined by
    \begin{equation*}
      \fw_\gamma=\biggl\{\gamma j\Bigm| j\in\ind,\ \exists E\in
      \smashoperator{\prod_{k\in K_\gamma}}
      \cset_k:A(X_E,Y_E,j)=1\biggr\}.
    \end{equation*}
    Since $A$ has no zero row, there is for every $j\in \ind$ an $y\in
    \XA$ such that $y_0=j$. Now if $\gamma j\in \fw_\gamma$, then we
    have, as mentioned above, that  
    \begin{equation*}
      \gamma y\in \smashoperator{\bigcup_{\substack{k=1\\
            |\mu^k|>m}}^n} C(\mu^k,I_k,J_k), 
    \end{equation*}
    which means that 
    \begin{equation*}
      j=y_0\in \bigl\{(\mu^k)_m\bigm| k\in \{1,2,\dotsc,n\},\ |\mu^k|>m\bigr\}.
    \end{equation*}
    Thus $\fw_\gamma$ is a finite set.
    Hence if we let $\fw_{m+1}=\bigcup_{\gamma\in \fw_m}\fw_\gamma$, then
    $\fw_{m+1}$ is a finite subset of the set
    \begin{multline*}
      \biggl\{\nu\in\ind^*\Bigm| |\nu|=m+1,\ \forall x\in
      \XA:A(\nu_{m+1},x_0)=1\Longrightarrow \\
      \nu x\in C(\mu,I,J)\cap
      \Bigl(\smashoperator[r]{\bigcap_{\substack{k=1\\
            |\mu^k|< m+1}}^n} 
      \XA\setminus C(\mu^k,I_k,J_k)\Bigr)\biggr\}.
    \end{multline*}

    Let $\gamma\in\fw_m$. If $|\mu^k|=m$, but $\mu^k\notin K_\gamma$, then
    $S_\gamma^*S_{\mu^k}=0$ by \eqref{eq:ckpc}, so the equality
    \begin{equation*}
      S_\gamma S_\gamma^*
      \left(1-S_{\mu^k}\prod_{i\in I_k}S_i^*S_i
        \prod_{j\in J_k}(1-S_j^*S_j)S_{\mu^k}^*\right) 
      =   S_\gamma S_\gamma^* 
    \end{equation*}
    holds.

    Thus we have that
    \begin{equation*}
      \begin{split}
        S_\gamma S_\gamma^* \Biggl( \smashoperator[r]{\prod_{\substack{k=1\\
              |\mu^k|=m}}^n} \Biggr( 1-S_{\mu^k} \prod_{i\in I_k}
            S_i^*S_i \prod_{j\in J_k} (1-S_j^*S_j) S_{\mu^k}^* \biggr)
        \Biggr)=&\\
        S_\gamma S_\gamma^* \Biggl( \smashoperator[r]{\prod_{k\in
            K_\gamma}} \biggl(1-S_{\mu^k} 
            \prod_{i\in I_k} S_i^*S_i \prod_{j\in J_k} (1-S_j^*S_j)
            S_{\mu^k}^* \biggr) \Biggr)=&\\
        S_\gamma S_\gamma^* S_\gamma \Biggl(
        \smashoperator[r]{\prod_{k\in K_\gamma}} 
        \biggl(1-\prod_{i\in I_k} S_i^*S_i \prod_{j\in J_k}
        (1-S_j^*S_j) \biggr) \Biggr) S_\gamma^*=& \\
        S_\gamma S_\gamma^* S_\gamma \Biggl(
        \smashoperator[r]{\prod_{k\in K_\gamma}} \biggl(
        \smashoperator[r]{\sum_{E_k\in \cset_k}} \Bigl(
        \smashoperator[r]{\prod_{\substack{i\in (J_k\cap E_k)\\ 
              \cup (I_k\setminus E_k)}}}
        S_i^*S_i\hspace{-0.7em} 
        \smashoperator[r]{\prod_{\substack{j\in (I_k\cap E_k)\\ 
              \cup (J_k\setminus E_k)}}}
        (1-S_j^*S_j) \Bigr) \biggr) \Biggr)
        S_\gamma^*=&\\
        S_\gamma S_{\gamma_m}^*S_{\gamma_m} \Biggl( 
        \smashoperator[r]{\prod_{k\in K_\gamma}} 
        \biggl( \sum_{E_k\in \cset_k} \Bigl(
        \smashoperator[r]{\prod_{\substack{i\in (J_k\cap E_k)\\ 
              \cup (I_k\setminus E_k)}}}
        S_i^*S_i \hspace{-0.7em}
        \smashoperator[r]{\prod_{\substack{j\in (I_k\cap E_k)\\ 
              \cup (J_k\setminus E_k)}}}
        (1-S_j^*S_j) \Bigr) \biggr) \Biggr)
        S_\gamma^*=&\\
        S_\gamma
        \Biggl(\smashoperator[r]{\sum_{E\in\pnes}}S_{\gamma_m}^*S_{\gamma_m}
        \smashoperator[l]{\prod_{k\in K_\gamma}} 
        \smashoperator[r]{\prod_{\substack{i\in (J_k\cap E_k)\\ 
              \cup (I_k\setminus E_k)}}}
        S_i^*S_i\hspace{-0.7em}
        \smashoperator[r]{\prod_{\substack{j\in (I_k\cap E_k)\\ 
              \cup (J_k\setminus E_k)}}} 
        (1-S_j^*S_j) \Biggr)
        S_\gamma^*=&\\
        S_\gamma\sum_{E\in\pnes} 
        \Biggl(\prod_{i\in I_E} S_i^*S_i \prod_{j\in J_E} 
        (1-S_j^*S_j) \Biggr) S_\gamma^*=&\\
        S_\gamma \sum_{E\in \pnes}
        \smashoperator{\sum_{\substack{j\in\ind \\
              A(I_E,J_E,j) = 1}}}
        S_jS_j^*S_\gamma^* =&
        \smashoperator{\sum_{\eta\in \fw_\gamma}}S_\eta S_\eta^*,
      \end{split}
    \end{equation*}
    where the second equality follows from \eqref{eq:ckpa}, the
    fact that $\mu^k=\gamma$ 
    for $k\in K_\gamma$ and that $S_\gamma$ is a partial isometry, the third
    follows from Lemma \ref{modsat}, the fourth from \eqref{eq:ckpd}
    (if $A(\gamma_i,\gamma_{i+1})=0$ for some
    $i=\{1,2,\dotsc,|\gamma|-1\}$, then $S_\gamma=0$ according to
    (\ref{eq:ckd}) and (\ref{eq:ckpb}), and the equality
    still holds), the fifth from the distribute law,
    the sixth from the definition of $I_E$ and $J_E$, the seventh
    from \eqref{eq:ckb} (that $ A(I_E,J_E,j)$ vanishes for
    all but a finite number of $j$'s in $\ind$ follows from the fact
    that $\fw_\gamma$ is finite), and the eight from
    \eqref{eq:ckpb} and the definition of $\fw_\gamma$.

   Thus we have that
   \begin{equation*}
     \begin{split}
        S_\mu\prod_{i\in I}S_i^*S_i\prod_{j\in J}&(1-S_j^*S_j)S_\mu^*
        \Biggl(\prod_{k=1}^n \biggl(1-S_{\mu^k}\prod_{i\in I_k}
        S_i^*S_i\prod_{j\in J_k}(1-S_j^*S_j)S_{\mu^k}^*\biggr)\Biggr)\\
        =&\smashoperator[r]{\sum_{\gamma\in \fw_m}} 
        S_\gamma S_\gamma^*\hspace{.7ex} \smashoperator{\prod_{\substack{k=1\\ 
              |\mu^k|\ge m}}^n}
        \Biggl(1-S_{\mu^k}\prod_{i\in I_k}S_i^*S_i
        \prod_{j\in J_k}(1-S_j^*S_j)S_{\mu^k}^*\Biggr)\\
        =&\smashoperator[r]{\sum_{\eta\in \fw_{m+1}}} 
        S_\eta S_\eta^*\hspace{.7ex} 
        \smashoperator{\prod_{\substack{k=1\\ |\mu^k|\ge
              m+1}}^n}\hspace{.7ex}
        \Biggl(1-S_{\mu^k} \prod_{i\in I_k}
        S_i^*S_i\prod_{j\in J_k}(1-S_j^*S_j)S_{\mu^k}^*\Biggr),
      \end{split}
    \end{equation*}
    which finalizes the induction proof of the claim.
  \end{proof}
  
  Hence the Boolean homomorphism $\eta$ from $\BO$ to $\BA$ is
  injective, and since it obviously also is surjective, it is
  invertible. Thus there is a Boolean homomorphism from $\BA$ to
  the Boolean algebra of projections in the unital abelian
  $\cs$-algebra 
  \begin{equation*}
    \cs\Bigl(\bigl\{S_gS_g^*\in\widetilde{\O}_A\bigm|
    g\in\fg\ind\bigr\}\Bigr),
  \end{equation*}
  sending $D_g$ to $S_gS_g^*$ for every $g\in\fg\ind$. Hence
  $(S_g)_{g\in\fg{\ind}}$ also satisfies \eqref{eq:oxd} of 
  Definition \ref{def:ox}. So it follows from the universal property of
  $\csp{\XA,\theta_\XA,\fg{\ind}}$ that there is a unital $*$-homomorphism
  $\psi$ from $\csp{\XA,\theta_\XA,\fg{\ind}}$ to $\widetilde{\O}_A$, such that
  $\psi(s_g)= S_{b_1}S_{b_2}\dotsm S_{b_k}$ for $g\in \fg\ind$
  written in the reduced form $b_1b_2\dotsm b_k$.
  
  We have that
  $\psi\bigl(\tilde{\phi}(S_i)\bigr)=\psi(s_i)=S_i$ for every $i\in\ind$ and
  that $\psi\bigl(\tilde{\phi}(1)\bigr)=1$, and since $\widetilde{\O}_A$ is
  generated by $\{S_i\mid i\in\ind\}\cup\{1\}$, this shows that
  $\psi\circ\tilde{\phi}=\Id_{\widetilde{\O}_A}$.
  
  According to Lemma \ref{lemma:onestruc}, $s_g=s_{b_1}s_{b_2}\dotsm
  s_{b_k}$ for every $g\in
    \fg\ind$ written in the reduced form $b_1b_2\dotsm b_k$. Thus we
    have that
  \begin{equation*}
    \tilde{\phi}\bigl(\psi(s_g)\bigr)=\tilde{\phi}(S_{b_1}S_{b_2}\dotsm
    S_{b_k})=s_{b_1}s_{b_2}\dotsm s_{b_k}=s_g
  \end{equation*} for
  every $g\in \fg\ind$ written in the reduced form $b_1b_2\dotsm b_k$,
  and since $\csp{\XA,\theta_\XA,\fg{\ind}}$ is generated by
  $(s_g)_{g\in\fg\ind}$, this shows that
  $\tilde{\phi}\circ\psi=\Id_{\csp{\XA,\theta_\XA,\fg{\ind}}}$.  

  Thus $\psi$ is a unital
  $*$-isomorphism from $\csp{\XA,\theta_\XA,\fg{\ind}}$ to
  $\widetilde{\O}_A$ which maps $s_i$ to $S_i$ for every $i\in\ind$,
  and since $\widetilde{\O}_A$ is generated by its unit and
  $\{S_i\mid i\in\ind\}$, $\csp{\XA,\theta_\XA,\fg{\ind}}$ is generated by
  its unit and $\{s_i\mid i\in\ind\}$ .
\end{proof}

\section{The ideal structure of $\csp{X,\theta,G}$} \label{sec:ideal}
One of the advantages of having a unified construction of the
$\cs$-algebras associated to one-sided shift spaces, crossed
product of two-sided shift spaces and Cuntz-Krieger algebras is
that it is easy to obtain results which holds for all of these
$\cs$-algebras (and of course other $\cs$-algebras which can be constructed as
$\cs$-algebras of discrete partial dynamical systems) and results which relate
these different kind of $\cs$-algebras to each other.

We will in this section for every discrete partial dynamical system
$(X,\theta,G)$ show how a $\theta$-invariant (see the definition of
$\theta$-invariant 
below) subset of $X$ gives raise to an ideal in
$\csp{X,\theta,G}$, and from this construct an injective order
preserving map between certain $\theta$-invariant subsets of $X$ and
ideals of $\csp{X,\theta,G}$. This will enable us to recover some
well-known result about the ideal structure of crossed
product of two-sided shift spaces and Cuntz-Krieger algebras, and will
shade new light on the ideal structure of $\cs$-algebra associated to
one-sided shift spaces.

We will also obtain a result which relates the $\cs$-algebras of two
different partial dynamical systems and use this to show that for
two-sided shift spaces having a certain property, the crossed product
of the two-sided shift space is a quotient of the $\cs$-algebra
associated to the corresponding one-sided shift space. This lays the
ground for a description of the $K$-theory of the $\cs$-algebra
associated to the one-sided shift space which is explained in
\cite{tmcseii} and \cite{tmcseiv}. 

\begin{definition}
   Let $(X,\theta,G)$ be a discrete partial dynamical system and
   $(\theta_g)_{g\in G}$ the partial one-to-one maps of $\theta$. Then we
   say that a subset $Y$ of $X$ is \emph{$\theta$-invariant} if
   $\theta_g(Y)\subseteq Y$ for all $g\in G$.
\end{definition}
  
Let $(X,\theta,G)$ be a discrete partial dynamical system,
$(D_g)_{g\in G}$ the domains and $(\theta_g)_{g\in G}$ the partial
one-to-one maps of $\theta$, and let $Y$ be an $\theta$-invariant
subset of $X$. Then we also have that $\theta_g(X\setminus Y)\subseteq
X\setminus Y$ for all $g\in G$, so if we for every $g\in G$ by
$\theta_{g|X\setminus Y}$ denote the restriction of $\theta_g$ to
$D_{g\inv}\cap X\setminus Y$, then the triple 
\begin{equation*}
  \left((D_g\cap X\setminus Y)_{g\in G},
    (\theta_{g|X\setminus Y})_{g\in G}\right)
\end{equation*}
is a partial action of $G$ on $X\setminus Y$ which we will denote by 
$\theta_{|X\setminus Y}$. 
Thus $(X\setminus Y,G,\theta_{|X\setminus Y})$ 
is a discrete partial dynamical system.

\begin{proposition} \label{shortexact}
  Let $(X,\theta,G)$ be a discrete partial dynamical system and
  let $Y$ be a $\theta$-invariant subset of $X$.
  
  Let $\ideal(Y)$ be the ideal of $\csp{X,\theta,G}$ generated
  by the set
  \begin{equation*}
    \{\boomap(A)\mid A\in \B,\ A\subseteq Y\}.
  \end{equation*}
  Then the quotient $\csp{X,\theta,G}/\ideal(Y)$ is
  isomorphic to $\csp{X\setminus Y,\theta_{|X\setminus Y},G}$.

  More precisely: if $\bigl(s_g^X\bigr)_{g\in G}$ denotes the generators of
  $\csp{X,\theta,G}$, and $\bigl(s_g^{X\setminus Y}\bigr)_{g\in G}$ denotes
  the generators of 
  $\csp{X\setminus Y,\theta_{|X\setminus Y},G}$, then the map
  \begin{equation*}
    s_g^X+\ideal(Y)\mapsto s_g^{X\setminus Y}
  \end{equation*}
  extends to a $*$-isomorphism from 
  $\csp{X,\theta,G}/\ideal(Y)$ to $\csp{X\setminus
    Y,\theta_{|X\setminus Y},G}$ which maps
  $\boomap(A)+\ideal(Y)$ to 
  $\boomap[(X\setminus Y,\theta_{|X\setminus Y},G)]
  (A\cap X\setminus Y)$
  for every $A\in \B$.
\end{proposition}

\begin{proof}
  Let $(D_g)_{g\in G}$ denote the domains of $\theta$ and let $\BXY$ be
  the Boolean algebra on $X\setminus Y$ generated by
  $\{D_g\cap X\setminus Y\mid g\in G\}$. We claim that the following
  identity holds:
  \begin{equation} \label{ww}
    \BXY=\{A\cap X\setminus Y\mid A\in \B\}.
  \end{equation}
  To see this, notice first that $\{A\cap X\setminus Y\mid A\in \B\}$
  is a Boolean algebra on $X\setminus Y$. Then let $g\in G$. Since
  $D_g\in\B$, we have that 
  \begin{equation*}
    D_g\cap X\setminus Y\in \{A\cap X\setminus Y\mid A\in \B\},
  \end{equation*}
  and therefore that $\BXY\subseteq \{A\cap X\setminus Y\mid A\in \B\}$.

  It is easy to check that $\{A\subseteq X\mid A\cap X\setminus Y\in
  \BXY\}$ is a Boolean algebra on $X$, and since we have that 
  \begin{equation*}
    D_g\in\{A\subseteq
    X\mid A\cap X\setminus Y\in \BXY\}
  \end{equation*}
  for every $g\in G$, it follows that $\B\in\{A\subseteq X\mid A\cap
  X\setminus Y\in \BXY\}$,  
  which shows that \eqref{ww} holds.
 
  It follows from (\ref{ww}) that the map 
  \begin{equation*}
    \pi:A\mapsto A\cap X\setminus Y
  \end{equation*}
  is a Boolean homomorphism from $\B$ to $\BXY$ which maps $D_g$ to
  $D_g\cap X\setminus Y$ for every $g\in G$. Thus
  $\boomap[(X\setminus Y,\theta_{|X\setminus Y},G)]\circ\pi$ is
  a Boolean homomorphism from $\B$ to the Boolean algebra of
  projections in the unital abelian $\cs$-algebra
  $\cs\Bigl(\bigl\{s_g^{X\setminus Y}s_g^{X\setminus Y*}\bigm| g\in G\bigr\}\Bigr)$
  which maps $D_g$ to $s_g^{X\setminus Y}s_g^{X\setminus Y*}$ for
  every $g\in G$. So it follows from the universal
  property of $\csp{X,\theta,G}$ that there exists a 
  $*$-homomorphism $\psi$ from $\csp{X,\theta,G}$ to 
  $\csp{X\setminus Y,\theta_{|X\setminus Y},G}$ which maps $s_g^X$ to
  $s_g^{X\setminus Y}$ for every $g\in G$, and $\boomap(A)$ to 
  $\boomap[(X\setminus Y,\theta_{|X\setminus Y},G)]\circ\pi(A)$ 
  for every $A\in\B$. We are done with the proof when we have shown that
  $\ker\psi=\ideal(Y)$. 

  If $A\subseteq Y$, then $\pi(A)=\emptyset$ and
  therefore $\psi\bigl(\boomap(A)\bigr)=0$. Thus $\ideal(Y)\subseteq \ker
  \psi$.

  The quotient map from $\csp{X,\theta,G}$ to
  $\csp{X,\theta,G}/\ideal(Y)$ induces a Boolean map
  $\widetilde{\pi}$ from the Boolean algebra of projections in the
  unital abelian $\cs$-algebra
  \begin{equation*}
    \cs\Bigl(\bigl\{s_g^Xs_g^{X*}\bigm| g\in G\bigr\}\Bigr)
  \end{equation*}
  to the Boolean algebra of projections in the
  unital abelian $\cs$-algebra
  \begin{equation*}
    \cs\Bigl(\bigl\{s_g^Xs_g^{X*}+\ideal(Y)\in
    \csp{X,\theta,G}/\ideal(Y)\bigm| g\in G\bigr\}\Bigr) 
  \end{equation*}
  which maps $p$ to
  $p+\ideal(Y)$ for every projection $p$ in
  $\cs\Bigl(\bigl\{s_g^Xs_g^{X*}\mid g\in G\bigr\}\Bigr)$. The map
  $\widetilde{\pi}\circ\boomap$ 
  is then a Boolean homomorphism from $\B$ to the Boolean algebra of
  projections in the
  unital abelian $\cs$-algebra
  \begin{equation*}
    \cs\Bigl(\bigl\{s_g^Xs_g^{X*}+\ideal(Y)\in
    \csp{X,\theta,G}/\ideal(Y)\bigm| g\in G\bigr\}\Bigr)
  \end{equation*}
  which maps $D_g$ to $s_g^Xs_g^{X*}+ \ideal(Y)$ for every
  $g\in G$. Since $\widetilde{\pi}\circ\boomap (A)=0$ if $A\subseteq
  Y$, we have that
  $\widetilde{\pi}\circ\boomap$ induces a Boolean homomorphism from
  $\BXY$ to the Boolean algebra of projections in the
  unital abelian $\cs$-algebra
  \begin{equation*}
    \cs\Bigl(\bigl\{s_g^Xs_g^{X*}+\ideal(Y)\in
    \csp{X,\theta,G}/\ideal(Y)\bigm| g\in G\bigr\}\Bigr)
  \end{equation*}
  which maps $D_g\cap
  X\setminus Y$ to $s_g^Xs_g^{X*}+ \ideal(Y)$ for every $g\in
  G$. Thus it follows from the universal property of 
  $\csp{X\setminus Y,\theta_{|X\setminus Y},G}$ that there exists
  a $*$-homomorphism $\tau$ from  
  $\csp{X\setminus Y,\theta_{|X\setminus Y},G}$ to
  $\csp{X,\theta,G}/\ideal(Y)$ which maps $s_g^X$ to
  $s_g^X+\ideal(Y)$ for every $g\in G$.

  The $*$-homomorphism $\psi$ from $\csp{X,\theta,G}$ to 
  $\csp{X\setminus Y,\theta_{|X\setminus Y},G}$ induces a 
  $*$-homomorphism $\widetilde{\psi}$ from 
  $\csp{X,\theta,G}/\ker(\psi)$ to 
  $\csp{X\setminus Y,\theta_{|X\setminus Y},G}$ which maps 
  $s_g^X+\ker(\psi)$ to $s_g^{X\setminus Y}$. 
  Thus $\tau\circ\widetilde{\psi}$ is a $*$-homomorphism from \\
  $\csp{X,\theta,G}/\ker(\psi)$ to
  $\csp{X,\theta,G}/\ideal(Y)$ which maps $s_g^X+\ker(\psi)$ to
  $s_g^X+\ideal(Y)$. This shows that $\ker(\psi)\subseteq \ideal(Y)$. 
\end{proof}

Let $(X,\theta,G)$ be a discrete partial dynamical system and $Y$ a
$\theta$-invariant subset of $X$. Clearly, the ideal $\ideal(Y)$ from
Proposition \ref{shortexact} only depends of the set 
$\{A\in \B\mid A\subseteq Y\}$.
We will for a subset $Y$ of $X$ call the set
\begin{equation*}
  \bigcup \{A\in \B\mid A\subseteq Y\}
\end{equation*}
for the \emph{$\theta$-admissible core} of $Y$, and we will call $Y$
\emph{$\theta$-admissible} if it is equal to its $\theta$-admissible
core. We then have for every $\theta$-invariant subset $Y$ of $X$ that
the ideal $\ideal(Y)$ from Proposition \ref{shortexact} its identical
to the ideal $\ideal(Y^\circ)$ where $Y^\circ$ denotes the
$\theta$-admissible core of $Y$.

Notice that a subset $Y$ of $X$ is $\theta$-admissible if and only if
there for every $x\in Y$ exists a $A\in\B$ such that $x\in A\subseteq Y$.

We have just seen how an $\theta$-invariant $\theta$-admissible subset
of $X$ gives raise to an ideal of $\csp{X,\theta,G}$. We will now go
the other way and from an ideal of $\csp{X,\theta,G}$ construct a
$\theta$-invariant $\theta$-admissible subset of $X$.

\begin{proposition} \label{prop:invariant}
  Let $(X,\theta,G)$ be a discrete partial dynamical system and let
  $\ideal$ be an ideal in $\csp{X,\theta,G}$. Then the set 
  \begin{equation*}
    \uset=\bigcup\{A\in \B\mid \boomap(A)\in\ideal\}
  \end{equation*}
  is a $\theta$-invariant $\theta$-admissible subset of $X$.
\end{proposition}

\begin{proof}
  The set $\uset$ is clearly a $\theta$-admissible subset of $X$, and
  it follows from Lemma \ref{booinvariant} that it is $\theta$-invariant.
\end{proof}

\begin{proposition} \label{prop:ideal}
  Let $(X,\theta,G)$ be a discrete partial dynamical system and let
  $Y$ be a $\theta$-invariant $\theta$-admissible subset of $X$. Then
  we have that
  \begin{equation*}
    Y=U\bigl(\ideal(Y)\bigr),
  \end{equation*}
  where $\ideal(Y)$ is as in Proposition \ref{shortexact} and
  $U\bigl(\ideal(Y)\bigr)$ is as in Proposition \ref{prop:invariant}.
\end{proposition}

\begin{proof}
  Since $Y$ and $U\bigl(\ideal(Y)\bigr)$ both are $\theta$-admissible
  subsets of $X$, it is enough to prove that the equivalence
  \begin{equation*}
    A\subseteq Y\iff \boomap(A)\in \ideal(Y)
  \end{equation*}
  holds for all $A\in \B$. It is clear that $\boomap (A)\in \ideal(Y)$ if
  $A\subseteq Y$. 

  Assume that $A$ is not a subset of $Y$. Then $A\cap X\setminus Y$ is
  non-empty, and since 
  $\boomap[(X\setminus U,\theta_{|X\setminus Y},G)]$ is injective, 
  $\boomap[(X\setminus Y,\theta_{|X\setminus Y},G)]
  (A\cap X\setminus Y)$ is a non-zero element of 
  $\csp{X\setminus Y,\theta_{|X\setminus Y},G}$. 
  By Proposition \ref{shortexact}, this
  implies that $\boomap(A)+\ideal(Y)$ is non-zero,
  and thus that $\boomap(A)$ is not in $\ideal(Y)$. 
\end{proof}

From Proposition \ref{prop:ideal} now directly follows the promised
theorem about the existence of  an injective order
preserving map between certain $\theta$-invariant subsets of $X$ and
ideals of $\csp{X,\theta,G}$: 

\begin{theorem} \label{theorem:ideal}
  Let $(X,\theta,G)$ be a discrete partial dynamical system and let
  for every  $\theta$-invariant $\theta$-admissible subset $Y$ of
  $X$, $\ideal(Y)$ be the ideal of $\csp{X,\theta,G}$ generated
  by the set 
  \begin{equation*}
    \{\boomap(A)\mid A\in \B,\ A\subseteq Y\}.
  \end{equation*}
  Then the map
  \begin{equation*}
    Y\mapsto \ideal(Y)
  \end{equation*}
  is an injective order preserving (i.e., $Y\subseteq Z\Rightarrow
  \ideal(Y)\subseteq \ideal(Z)$) map from the set of $\theta$-invariant
  $\theta$-admissible subsets of $X$ to the set of ideals of
\end{theorem}

We will now present a result which relate the $\cs$-algebras of two
partial dynamical systems to each other:

\begin{theorem} \label{theorem:surmap}
  Let $G$ be a group and let for $i\in\{1,2\}$, $X_i$ be a set,
  $\theta_i$ a partial action of $G$ on $X_i$ and
  $(\theta_g^{X_i})_{g\in G}$ the partial one-to-one maps of
  $\theta_i$. If there exists a 
  Boolean homomorphism $\boof$ from $\Bx(X_1,\theta_1,G)$ to
  $\Bx(X_2,\theta_2,G)$ such that the identity
  \begin{equation*}
    \boof\bigl(\theta_g^{X_1}(A)\bigr)=\theta_g^{X_2}\bigl(\boof(A)\bigr)
  \end{equation*}
  holds for all $A\in \Bx(X_1,\theta_1,G)$ and $g\in G$, then
  $\csp{X_2,\theta_2,G}$ is a quotient of $\csp{X_1,\theta_1,G}$.
 
  More precisely: if for $i\in\{0,1\}$, $(s_g^{X_i})_{g\in G}$ denotes the
  generators of $\csp{X_i,\theta_i,G}$,
  then there exists a surjective $*$-homomorphism from
  $\csp{X_1,\theta_1,G}$ to $\csp{X_2,\theta_2,G}$ which maps $s_g^{X_1}$ to
  $s_g^{X_2}$ for every $g\in G$, and $\boomap[(X_1,\theta_1,G)](A)$ to
  $\boomap[(X_2,\theta_2,G)]\bigl(\boof(A)\bigr)$ for every
  $A\in\Bx(X_1,\theta_1,G)$. 

  The kernel of this $*$-homomorphism is
  the ideal generated by the set
  \begin{equation*}
    \{\boomap[(X_1,\theta_1,G)](A)\mid A\in\Bx(X_1,\theta_1,G),\
    \boof(A)=\emptyset\}.
  \end{equation*}
\end{theorem}

\begin{proof}
  Notice that $Y=\bigcup\{A\in\Bx(X_1,\theta_1,G)\mid
  \boof(A)=\emptyset\}$ is a $\theta_1$-invariant
  $\theta_1$-admissible subset of $X_1$ and that 
  if we let $\ideal(Y)$ be as in Proposition \ref{shortexact}, then
  $\ideal(Y)$ is the ideal generated by the set
  \begin{equation*}
    \{\boomap[(X_1,\theta_1,G)](A)\mid A\in\Bx(X_1,\theta_1,G),\
    \boof(A)=\emptyset\}. 
  \end{equation*}

  It directly follows from the universal property of
  $\csp{X_1,\theta_1,G}$ that there exists a $*$-homomorphism $\psi$ from
  $\csp{X_1,\theta_1,G}$ to $\csp{X_2,\theta_2,G}$ which maps $s^{X_1}_g$ to
  $s^{X_2}_g$ for every $g\in G$ and $\boomap[(X_1,\theta_1,G)](A)$ to
  $\boomap[(X_2,\theta_2,G)](\boof(A))$ for every
  $A\in\Bx(X_1,\theta_1,G)$, and it is clear that the set
  \begin{equation*}
    \{\boomap[(X_1,\theta_1,G)](A)\mid A\in\Bx(X_1,\theta_1,G),\
    \boof(A)=\emptyset\} 
  \end{equation*}
  is contained in the kernel of $\psi$. Since $\psi$ maps $s^{X_1}_g$ to
  $s^{X_2}_g$, $\psi$ is surjective, so we only have to show that
  the kernel of $\psi$ is contained in $\ideal(Y)$. We will do that by
  proving that there exists a $*$-homomorphism 
  from $\csp{X_1,\theta_1,G}/\ker\psi$ to $\csp{X_1,\theta_1,G}/\ideal(Y)$
  which maps $s^{X_1}_g+\ker\psi$ to $s^{X_1}_g+\ideal(Y)$ for every $g\in G$.
  
  Let us denote the generators of
  $\csp{X_1\setminus Y,\theta_{1|X_1\setminus Y},G}$ by
  $(s^{X_1\setminus Y}_g)_{g\in G}$.  
  Since $\psi$ is surjective, it induces a $*$-isomorphism from $
  \csp{X_1,\theta_1,G}/\ker\psi$ to $ \csp{X_2,\theta_2,G}$ which maps
  $s^{X_1}_g+\ker\psi$ to $s^{X_2}_g$ and it follows from Proposition
  \ref{shortexact} 
  that there exists a $*$-isomorphism from $
  \csp{X_1,\theta_1,G}/\ideal(Y)$ to $\csp{X_1\setminus
    Y,\theta_{1|X_1\setminus Y},G}$ which maps $s^{X_1}_g+\ideal(Y)$
  to $s^{X_1\setminus Y}_g$
  for every $g\in G$. So all we have to do is to show that there
  exists a $*$-homomorphism from $\csp{X_2,\theta_2,G}$ to
  $\csp{X_1\setminus Y,\theta_{1|X_1\setminus Y},G}$ which maps
  $s^{X_2}_g$ to $s^{X_1\setminus Y}_g$ for every $g\in G$.

  Let $(D_g^{X_1})_{g\in G}$ denote the domains of $\theta_1$ and
  $(D_g^{X_2})_{g\in G}$ denote the domains of $\theta_2$. We then have
  that
  \begin{equation*}
    \boof\bigl(D_g^{X_1}\bigr)=\boof\bigl(\theta_g^{X_1}(X_1)\bigr)=
    \theta_g^{X_2}\bigl(\boof(X_1)\bigr) 
    =\theta_g^{X_2}\bigl(X_2\bigr)= D_g^{X_2}
  \end{equation*}
  for every $g\in G$, $\boof$ is
  surjective, and if  $A,B\in\Bx(X_1,\theta_1,G)$ and
  $\boof(A)=\boof(B)$, then $A\cap X_1\setminus Y=B\cap X_1\setminus
  Y$. So the map  
  \begin{equation*}
    \boof(A)\mapsto A\cap X_1\setminus Y
  \end{equation*}
  is a well-defined map from $\Bx(X_2,\theta_2,G)$ to
  $\Bx(X_1\setminus Y,\theta_{1|X_1\setminus Y},G)$ which maps $D_g$
  to $D_g\cap X_1\setminus Y$ for every $g\in G$. It is easy to check
  that it is a Boolean homomorphism, so it follows from the universal
  property of $\csp{X_2,\theta_2,G}$ that there  exists a
  $*$-homomorphism from $\csp{X_2,\theta_2,G}$ to 
  $\csp{X_1\setminus Y,\theta_{1|X_1\setminus Y},G}$ which maps
  $s^{X_2}_g$ to $s^{X_1\setminus Y}_g$ for every $g\in G$.
\end{proof}

\subsection{The ideal structure of $C(\TSS)\rtimes_{\tsh^\star}\Z$} 
Let $(\TSS,\tsh)$ be a two-sided shift space over the finite
alphabet $\al$. As it is proved in the proof of
Theorem \ref{kryds}, the Boolean algebra
$\Bx(\TSS,\theta_\TSS,\fg\al)$ is equal the Boolean algebra of clopen subsets
of $\TSS$, and since the clopen subsets generate the topology of
$\TSS$, a subset of $\TSS$ is $\theta_\TSS$-admissible if and only it is
open. It is easy to check that a subset $Y$ of $\TSS$ is
$\theta_\TSS$-invariant if and only if $\tsh(Y)=Y$. 

Thus we recover from Theorem \ref{kryds} and \ref{theorem:ideal} the
well know fact 
that there exists an injective order preserving map from the set of
open $\tsh$-invariant subsets of $\TSS$ to the set of ideals of
$C(\TSS)\rtimes_{\tsh^\star}\Z$. If $(\TSS,\tsh)$ is free (meaning that
$\tsh^n(x)\ne x$ for every $x\in\TSS$ and every $n\in\Z\setminus\{0\}$)
then this map
 is bijective, cf. \cite{MR0241994}*{Proposition 5.10 and
  Th\'eor\`eme 5.15}.

\subsection{The ideal structure of $\widetilde{\O}_A$}
\label{subsec:idealkryds} 
Let $\ind$ be an arbitrary index set and let
$A=(A(i,j))_{i,j\in\ind}$ be a matrix with entries in $\{0,1\}$ and
having no zero rows. There is, as mentioned in Remark \ref{remark:exellaca},
a homeomorphism from $\h{\mathsf{X}}_A^+$ to
$\widetilde{\Omega}_A$. Thus there is a bijective correspondence
between elements of $\BA$ and clopen subsets of
$\widetilde{\Omega}_A$. This correspondence extend to a correspondence
between $\theta_\XA$-admissible subsets of $\XA$ and open subsets of
$\widetilde{\Omega}_A$, and this correspondence takes
$\theta_\XA$-invariant subsets to invariant subsets of
$\widetilde{\Omega}_A$. 

Thus we get from Theorem \ref{theorem:CK} and \ref{theorem:ideal}
an injective order preserving map from the set of open invariant
subsets of $\widetilde{\Omega}_A$ to the set of ideals of
$\widetilde{\O}_A$. Exel and Laca have in
\cite{MR2000i:46064}*{Theorem 15.1} proved that this map is bijective
if the directed graph $Gr(A)$ of $A$ has no transitory circuits
(cf. \cite{MR2000i:46064}*{Section 12}).

\subsection{The ideal structure of $\O_{\OSS}$} 
Let $(\OSS,\osh)$ be a one-sided shift space over the finite alphabet
$\al$, and let $\laOSS$ be the language of $\OSS$
(cf. \cite{MR97a:58050}*{\S1.3}), that is 
\begin{equation*}
  \laOSS=\{\mu\in\al^*\mid \exists x\in\OSS,0\le k\le
  l:x_kx_{k+1}\dotsm x_l=\mu\}.
\end{equation*}
We let $(\OSS,\theta_{\OSS},\fg{\al})$ be the discrete partial
dynamical system of Theorem \ref{theorem:one}, and 
$(D_g)_{g\in\fg\al}$ be the domains and $(\theta_g)_{g\in\fg\al}$ the
partial one-to-one maps of $\theta_{\OSS}$.

Following \cite{MR1691469}, we let for every $x\in \OSS$ and every
$k\in \No$, $\Past_k(x)$ be the subset of $\laOSS$ defined by
\begin{equation*}
  \Past_k(x)=\{\mu\in \laOSS\mid \mu x\in \OSS, |\mu|=k\},
\end{equation*}
and define for every $l\in \No$ an equivalence relation $\sim_l$ on
$\OSS$ by
\begin{equation*}
  x\sim_l x' \iff \Past_l(x)=\Past_l(x').
\end{equation*}

We then define for every $(k,l)\in \indd$ an equivalence relation
$\eq$ on $\OSS$ by
\begin{equation*}
  x\eq y\iff x_{[0,k[}=y_{[0,k[}\land
  \Past_l(x_{[k,\infty[})=\Past_l(y_{[k,\infty[}), 
\end{equation*}
and we let $\ec{k}{x}{l}$ denote the equivalence class of $x$ under
this equivalence relation. We then define an order $\lee$ on $\indd$ by
\begin{equation*}
  (k_1,l_1)\lee (k_2,l_2) \iff k_1\le k_2 \land l_1-k_1\le l_2-k_2,
\end{equation*}
and notice that if $ (k_1,l_1)\lee (k_2,l_2)$, then
$\ec{k_2}{x}{l_2}\subseteq\ec{k_1}{x}{l_1}$ for all $x\in\OSS$. Notice also that
$\ec{k}{x}{l}\in\Bx(\OSS,\theta_{\OSS},\fg\al)$ for all $x\in\OSS$ and
all $(k,l)\in\indd$ because we have that
\begin{equation*}
  \ec{k}{x}{l}=\theta_{x_{[0,k[}}
  \left(\Biggl(\smashoperator[r]{\bigcap_{\mu\in\Past_l(x_{[k,\infty[})}}
    (D_{\mu\inv})\Biggr) \cap \Biggl(
  \smashoperator[r]{\bigcap_{\mu\in\neg\Past_l(x_{[k,\infty[})}}
  \OSS\setminus D_{\mu\inv}\Biggr)\right),
\end{equation*}
where $\neg\Past_l(x_{[k,\infty[})=\{\mu\in\laOSS\mid |\mu|=l,\
\mu\notin\Past_l(x_{[k,\infty[})\}$.

\begin{lemma} \label{lemma:clopen}
  There exists for every $A\in\Bx(\OSS,\theta_{\OSS},\fg\al)$ a
  $(k,l)\in\indd$ such that the implication
  \begin{equation*}
    x\in A\Rightarrow \ec{k}{x}{l}\subseteq A
  \end{equation*}
  holds for all $x\in\OSS$.
\end{lemma}

\begin{proof}
  Let $\Aa$ be the subset of $\Bx$ defined by
  \begin{equation*}
    \Aa=\bigl\{A\in\Bx(\OSS,\theta_{\OSS},\fg\al)\bigm| \exists
    (k,l)\in\indd\forall x\in A:\ec{k}{x}{l}\subseteq A\bigr\}. 
  \end{equation*}
  We must then show that $\Aa=\Bx(\OSS,\theta_{\OSS},\fg\al)$. Clearly
  $\OSS\in \Aa$. 
  Assume that $A,B\in \Aa$ and choose $(k_a,l_a),(k_b,l_b)\in \indd$
  such that the two implications
  \begin{equation*}
    x\in A\Rightarrow \ec{k_a}{x}{l_a}\subseteq A
  \end{equation*}
  and
  \begin{equation*}
    x\in B\Rightarrow \ec{k_b}{x}{l_b}\subseteq B
  \end{equation*}
  hold for all $x\in\OSS$.
  
  Let $k=\max\{k_a,k_b\}$ and $l=\max\{l_a-k_a,l_b-k_b\}+k$. Then we
  have that $(k_a,l_a),(k_b,l_b)\lee (k,l)$, from which it follows
  that the implication 
  \begin{equation*}
    x\in A\cap B\Rightarrow
    \ec{k}{x}{l}\subseteq\ec{k_a}{x}{l_a}\cap\
    \ec{k_b}{x}{l_b}\subseteq A\cap B,
  \end{equation*}
  holds for all $x\in\OSS$. This shows that $A\cap B\in \Aa$. We
  also have that the three implications
  \begin{gather*}
    x\in \OSS\setminus A\Rightarrow \ec{k_a}{x}{l_a}\subseteq 
    \OSS\setminus A,\\
    x\in \theta_a(A)\Rightarrow \ec{k_a+1}{x}{l_a}\subseteq\theta_a(A)
  \end{gather*}
  and
  \begin{equation*}
    x\in \theta_{a\inv}(A)\Rightarrow \ec{k_a}{x}{l_a+1}\subseteq
    \theta_{a\inv}(A) 
  \end{equation*}
  hold for all $x\in\OSS$. Thus the sets $\OSS\setminus A,\
  \theta_a(A)$, and $\theta_{a\inv}(A)$ all belong to $\Aa$. Hence
  $\Aa$ is a Boolean algebra containing $D_g$ for every 
  $g\in \fg\al$, which means that $\Aa=\Bx(\OSS,\theta_{\OSS},\fg\al)$.
\end{proof}

\begin{lemma} \label{lemma:open}
  A subset $Y$ of $\OSS$ is $\theta_{\OSS}$-admissible if and only if there
  for every $x\in Y$ exists a $(k,l)\in\indd$ such that
  $\ec{k}{x}{l}\subseteq Y$.  
\end{lemma}

\begin{proof}
  Let $\Aa$ be the subset of $\OSS$ defined by 
  \begin{equation*}
    \Aa=\bigl\{Y\subseteq \OSS\bigm| \forall x\in Y\exists
    (k,l)\in\indd: \ec{k}{x}{l}\subseteq Y\bigr\}.
  \end{equation*}
  We will show that a subset $Y$ of $\OSS$ belongs to $\Aa$ if and
  only if it is $\theta_{\OSS}$-admissible.

  It is clear that if $(Y_i)_{i\in\ind}$ is a family of elements of
  $\Aa$, then $\bigcup_{i\in\ind}Y_i\in \Aa$, and since it follows
  from Lemma \ref{lemma:clopen} that 
  $\Bx(\OSS,\theta_{\OSS},\fg\al)\subseteq\Aa$, we have that every
  $\theta_{\OSS}$-admissible subset of $\OSS$ belongs to $\Aa$.

  Now let $Y\in\Aa$. It is in order to show that $Y$ is
  $\theta_{\OSS}$-admissible enough to show that there for every
  $x\in Y$ exists a $A\in\Bx(\OSS,\theta_{\OSS},\fg\al)$ such that
  $x\in A\subseteq Y$, but we have that $x\in\ec{k}{x}{l}\subseteq Y$,
  and $\ec{k}{x}{l}\in\Bx(\OSS,\theta_{\OSS},\fg\al)$ for some 
  $(k,l)\in\indd$, so we are done. 
\end{proof}

It is not difficult to see that a subset $Y$ of $\OSS$ is
$\theta_{\OSS}$-invariant if and only if $\osh(Y)\subseteq Y$ and
$\osh\inv(Y)\subseteq Y$. Combining this with Theorem
\ref{theorem:one}, Theorem \ref{theorem:ideal} and Lemma
\ref{lemma:open} we get: 

\begin{theorem} \label{theorem:ossideal}
  Let $(\OSS,\osh)$ be a one-sided shift space over the finite alphabet
  $\al$, let $\mathcal{H}(\OSS)$ denote the set of subsets $Y$ of
  $\OSS$ which satisfies the following 3 conditions: 
\begin{enumerate}
  \item $\osh(Y)\subseteq Y$,
  \item $\osh\inv(Y)\subseteq Y$,
  \item $\forall x\in Y\exists (k,l)\in\indd: \ec{k}{x}{l}\subseteq
    Y$,
\end{enumerate}
and let for every subset $Y$ in $\mathcal{H}(\OSS)$, $\ideal(Y)$ be the
ideal of $\csp{\OSS,\theta_{\OSS},\fg\al}$ generated
by the set
\begin{equation*}
  \{\boomap[(\OSS,\theta_{\OSS},\fg\al)](A)\mid A\in
  \Bx(\OSS,\theta_{\OSS},\fg\al),\ A\subseteq Y\}.
\end{equation*}  
Then the map
\begin{equation*}
  Y\mapsto \ideal(Y)
\end{equation*}
is an injective order preserving (i.e., $ Y\subseteq Z\Rightarrow
\ideal(Y)\subseteq \ideal(Z)$) map from $\mathcal{H}(\OSS)$ to the
set of ideals of $\O_{\OSS}$.
\end{theorem}

If $Y$ is an open subset of $\OSS$ then there exists for every $x\in Y$
a $(k,l)\in\indd$ such that $\ec{k}{x}{l}\subseteq Y$ (in fact, one
can choose $l$ to be $0$). There may on the
other hand be subsets $Y$ of $\OSS$ which are not open, but with the
property that $\forall x\in Y\exists (k,l)\in\indd:
\ec{k}{x}{l}\subseteq Y$. If however $\OSS$ is a shift of finite
type (cf. \cite{MR97a:58050}*{\S2.1}), then this can not happen (in
fact according to \cite{MR0346134}*{Theorem 1}, $\OSS$ is of finite
type if and only if $\osh$ is an open map, and it is not difficult to
show that this is equivalent to $\ec{k}{x}{l}$ being clopen for all
$x\in\OSS$ and all $(k,l)\in\indd$). Thus if $\OSS$ is a shift of finite
type, then the set $\mathcal{H}(\OSS)$ from Theorem
\ref{theorem:ossideal} is the set of all open subsets $Y$ of
$\OSS$ with the property that $\osh(Y)\subseteq Y$ and 
$\osh\inv(Y)\subseteq Y$. Hence we
get the following corollary:

\begin{corollary}
  Let $(\OSS,\osh)$ be a one-sided shift of finite type and let for
  every subset $Y$ of $\OSS$, $\ideal(Y)$ be the 
  ideal of $\csp{\OSS,\theta_{\OSS},\fg\al}$ generated
  by the set
  \begin{equation*}
    \{\boomap[(\OSS,\theta_{\OSS},\fg\al)](A)\mid A\in
    \Bx(\OSS,\theta_{\OSS},\fg\al),\ A\subseteq Y\}.
  \end{equation*}  
  Then the map
  \begin{equation*}
    Y\mapsto \ideal(Y)
  \end{equation*}
  is an injective order preserving (i.e., $ Y\subseteq Z\Rightarrow
  \ideal(Y)\subseteq \ideal(Z)$) map from the set of open subsets $Y$ of
  $\OSS$ with the property that $\osh(Y)\subseteq Y$ and 
  $\osh\inv(Y)\subseteq Y$ to the set of ideals of $\O_{\OSS}$.
\end{corollary}

\subsection{Connections between the $\cs$-algebras of one- and
  two-sided shift spaces}
We will in this section let $\TSS$ be a two-sided shift space, and
\begin{equation} \label{eq:twoone}
  \OSS=\{(z_n)_{n\in \No}\mid (z_n)_{n\in \Z}\in \TSS\}
\end{equation}
be the corresponding one-sided shift space. We then have that map
\begin{equation*}
  \spro:(z_n)_{n\in \Z}\mapsto (z_n)_{n\in \No}
\end{equation*}
is a surjective continuous map from $\TSS$ to $\OSS$ and
$\spro\circ\tsh=\osh\circ\spro$, where $\tsh$ and $\osh$ are the maps
defined by \eqref{eq:twoshift} and \eqref{eq:oneshift}, respectively.

We will throughout this section let $(\TSS,\theta_\TSS,\fg{\al})$ be
the discrete partial dynamical system associated to $(\TSS,\tsh)$ as
done in Theorem \ref{kryds} and $(\OSS,\theta_{\OSS},\fg{\al})$ be the
discrete partial dynamical system associated to $(\OSS,\osh)$ as done in
Theorem \ref{theorem:one}. We let $(D_g^{\TSS})_{g\in\fg\al}$ and
$(\theta_g^{\TSS})_{g\in\fg\al}$ denote the domains and partial one-to-one
maps of $\theta_\TSS$, and $(D_g^{\OSS})_{g\in\fg\al}$ and
$(\theta_g^{\OSS})_{g\in\fg\al}$ denote the domains and partial one-to-one
maps of $\theta_{\OSS}$.

\begin{remark}
  If $\TSS$ only contains finitely many elements, then every element of
  $\TSS$ is periodic (meaning that there exists an $n\in\N$ such that
  $\tsh^n(z)=z$), so $\spro$ is bijective, and since it also satisfies
  that $\theta_g\circ\spro=\spro\circ\theta_g$ for every $g\in\fg\al$, it
  easily follows from Theorem \ref{kryds} and \ref{theorem:one} that
  there exists a $*$-isomorphism from $\O_{\OSS}$ to
  $C(\TSS)\rtimes_{\tsh^\star}\Z$ which maps $\sum_{a\in\al}S_a$ to $U$
  and $\ioss\big(1_{D_g^{\OSS}}\big)$ to $1_{D_g^{\TSS}}$  
  for every $g\in\fg\al$, where $U$ is as in
  Section \ref{sec:cross}, and $(S_a)_{a\in\al}$ and $\ioss$ are as in
  Section \ref{sec:onesided}.
\end{remark}

If $\TSS$ contains infinitely many elements, then
$C(\TSS)\rtimes_{\tsh^\star}\Z$ and $\O_{\OSS}$ are, as we will see
below, in general quit different. We will however show later that if
$\TSS$ has a certain property, then $C(\TSS)\rtimes_{\tsh^\star}\Z$ is
a quotient of $\O_{\OSS}$.

\begin{remark}
  If $Y$ is an open subset of $\OSS$ with the property that
  $\osh(Y)\subseteq Y$ and $\osh\inv(Y)\subseteq Y$, then
  $\spro\inv(Y)$ is an open subset of $\TSS$ and
  $\tsh(\spro\inv(Y))=\spro\inv(Y)$. This explains why
  $C(\TSS)\rtimes_{\tsh^\star}\Z$ in general have more ideals than
  $\O_{\OSS}$ if $\TSS$ (and thus $\OSS$) is of finite type.

  It is in fact easy to construct an example of a shift of finite type
  $\TSS$ such that $\O_{\OSS}$ is simple, but $\TSS$ contains
  infinitely many open $\tsh$-invariant subset and thus 
  $C(\TSS)\rtimes_{\tsh^\star}\Z$ infinitely many ideals according to Section
  \ref{subsec:idealkryds}. The full
  two-shift $\{0,1\}^\Z$ will for example do the trick (in this case
  $\O_{\OSS}$ will be the Cuntz-algebra $\O_2$, cf. \cite{MR0467330}).
\end{remark}

If $\TSS$ and $\OSS$ are not of finite type, then it might happen
that even though $\OSS$ is minimal (meaning that the only closed
subsets $Y$ of $\OSS$ such that $\osh(Y)\subseteq Y$ are $\OSS$ and
$\emptyset$, cf. \cite{MR97a:58050}*{\S13.7}), 
there exists a subset $Y$ of $\OSS$ which is neither equal to $\OSS$
nor $\emptyset$ and with the property that
$\osh(Y)\subseteq Y$, 
$\osh\inv(Y)\subseteq Y$ and $\forall x\in Y\exists
(k,l)\in\ind:\ec{k}{x}{l}\subseteq Y$, and thus
that $\O_{\OSS}$ is not simple. We will in Example \ref{exam:sub} see an
example of this phenomenon.

We will now describe a class of shift spaces for which
$C(\TSS)\rtimes_{\tsh^\star}\Z$ is a quotient of $\O_{\OSS}$, and we
will then describe a subclass of this class for which we can show that
the ideal of this quotient is a direct sum of a finite number of the
compact operators $\K$

This result has been used in \cite{tmcseii} and \cite{tmcseiv} to
compute the $K$-theory of $\O_{\OSS}$ and relate it to the $K$-theory of
$C(\TSS)\rtimes_{\tsh^\star}\Z$ for these classes of shift spaces.

\begin{definition}
  We say that a shift space $\TSS$ has
  \emph{property $\coni$} if for every $\mu\in \laOSS$ there exists an
  $x\in \OSS$ such that $\Past_{|\mu|}(x)=\{\mu\}$.
\end{definition}

An element $z\in\TSS$ is called \emph{left special} if there exists
$z'\in\TSS$ such that $\spro(z)=\spro(z')$, but $z_{-1}\ne z'_{-1}$.

Notice that $z\in\TSS$ is left special if and only if
$\Past_1(\spro(z))$ consists of a least two elements.

\begin{definition}
  We say that a shift space $\TSS$ has property $(**)$ if it has
  property $(*)$ and the number of left special elements of $\TSS$
  is finite, and no such left special word is periodic.
\end{definition}

It has in \cite{tmcseii} been proved that every finite shift space and
every minimal shift space with a finite number of left special
elements (for example every shift space of a primitive substitution
and every shift space of a Sturmian sequence) have property $(**)$, and
that the shift space of a non-regular Toeplitz sequence has property
$(*)$, but not necessarily property $(**)$.

We say that $z,z'\in \TSS$ are \emph{right shift tail equivalent} if
there exist $m,N\in\Z$ such that $z_n=z'_{n+m}$ for all $n>N$. 

\begin{lemma} \label{lemma:faa}
   Let $\TSS$ be a shift space which has property $(**)$ and let
   $(k,l)\in\indd$. Then the set
   \begin{equation*}
     \{x\in\OSS\mid \Past_l(x_{[k,\infty[})\text{ contains more than
       one element}\}
   \end{equation*}
   is finite, and every element of it is of the form $\spro(z)$, where
   $z$ is an element of $\TSS$ which is right shift tail equivalent to
   a left special element of $\TSS$. 
\end{lemma}

\begin{proof}
  Assume that $u,v,w\in\al^*$ are such that the last letter of
  $u$ is different from the last letter of $v$ and
  $uw,vw\in\Past_l(x_{[k,\infty[})$ ($w$ might be the empty
  word). Then $\Past_1(wx_{[k,\infty[})$ consists of a least two
  elements so $wx_{[k,\infty[}$ is equal to $\spro(z')$ for some left
  special element $z'$. 

  Thus the set
  \begin{equation*}
     \{x\in\OSS\mid \Past_l(x_{[k,\infty[})\text{ contains more than
       one element}\}
   \end{equation*}
   is finite, and every element of it is of the form $\spro(z)$, where
   $z$ is an element of $\TSS$ which is right shift tail equivalent to
   a left special element of $\TSS$. 
\end{proof}

\begin{lemma} \label{lemma:nonperiodic}
  Let $\TSS$ be a shift space which has property $(**)$, and let $z\in
  X$ be right shift tail equivalent to a left special element. Then
  $\spro(z)$ is not eventually periodic, meaning that there exist no
  $k,N\in\N$ such that $z_n=z_{n+k}$ for all $n>N$. 
\end{lemma}

\begin{proof}
  Assume that $m,N_1\in\Z$ and $z_n=z'_{n+m}$ for all $n>N_1$ with $z'$
  being a left special element, and that $k,N_2\in\N$ are such that
  $z_{n}=z_{n+k}$ for all $n>N_2$. We then have that $z'_{n+m}=z'_{n+m+k}$ for
  $n>\max\{N_1,N_2\}$, so the set $\{n\in\Z\mid z'_n\ne z'_{n+k}\}$
  is bounded above, and since $z'$ is left special and thus not
  periodic, this set is not empty. So we can define $n_0$ by the equation
  \begin{equation*}
    n_0=\max\{n\in\Z\mid z'_n\ne z'_{n+k}\}+1.
  \end{equation*}
  We then have that $z'_{n_0-1}\ne z'_{n_0-1+k}$ and that 
  $z'_{n}= z'_{n+k}$ for $n\ge n_0$. 

  We define a sequence $z''$ by letting $z''_{n_0+i+lk}=z'_{n_0+i}$ for
  $i\in\{0,1,\dotsc,k-1\}$ and $l\in\Z$. We then have for all $n\in\Z$ that 
  $z''_{[n,\infty[}=z'_{[n+lk,\infty[}$, where $l\in\Z$ is chosen such that
  $n+lk\ge n_0$, and thus that $z''\in\TSS$. We also have that $z''$ is
  periodic, that $z''_{n_0-1}=z'_{n_0-1+k}\ne z'_{n_0-1}$ and that
  $z''_{[n_0,\infty[}=z'_{[n_0,\infty[}$. Thus $\osh^{n_0}(z'')$ is an
  element of $\TSS$ which is both periodic and left special, but that
  contradicts the assumption that $\TSS$ has property $(**)$.
\end{proof}

\begin{lemma} \label{lemma:iso}
  Let $\TSS$ be a shift space which has property $(**)$. Then we have
  that 
  \begin{equation*}
    \{\spro(z)\}\in\Bx(\OSS,\theta_{\OSS},\fg\al)
  \end{equation*}
  for every $z\in\TSS$ which is right shift equivalent to a left
  special element. 
\end{lemma}

\begin{proof}
  Let $z$ be an element of $\TSS$ which is right shift equivalent to a left
  special element. Then there exist a left special element $z'$ and
  $m,N\in\Z$ such that $z_n=z'_{n+m}$ for all $n>N$. Since $z'$ is left
  special, there exist $a,b\in\al$ such
  that $a\ne b$ and $a,b\in\Past_1(\spro(z'))$. That means that 
  \begin{equation*}
    \spro(z')\in D_{a\inv}^{\OSS}\cap D_{b\inv}^{\OSS},
  \end{equation*}
  and since $\TSS$ only contains a finite number of special elements,
  $D_{a\inv}^{\OSS}\cap D_{b\inv}^{\OSS}$ is finite. Thus there
  exists a $k>0$ such that
  \begin{equation*}
    \{z'_{[k,\infty[}\}=D_{(au)\inv}^{\OSS}\cap D_{(bu)\inv}^{\OSS},
  \end{equation*}
  where $u=z'_{[0,k[}$. 

  Let us denote $\max\{m+N,0,m-k\}$ by $n_0$ and let $v=z'_{[k,k+n_0[}$  and
  $w=z_{[0,n_0[}$. We then have that $\spro(z)=wy$ where
  $y$ is the unique element in $\OSS$ such that
  $vy=z'_{[k,\infty[}$. Thus it follows from Lemma \ref{lemma:onefund}
  that 
  \begin{equation*}
    \{\spro(z)\}=\theta_{wv\inv}^{\OSS}\left(D_{(au)\inv}^{\OSS}\cap
      D_{(bu)\inv }^{\OSS}\right) 
    \in\Bx(\OSS,\theta_{\OSS},\fg\al). 
  \end{equation*}
\end{proof}

Remember that there exist an inclusion $\ioss$ of $\DX$ into
$\O_{\OSS}$ (cf. Section \ref{sec:onesided})
and a that $C(\TSS)$ sits inside
$C(\TSS)\rtimes_{\tsh^\star}\Z$ (cf. Section \ref{sec:cross}). We will denote
the inclusion of 
$C(\TSS)$ into $C(\TSS)\rtimes_{\tsh^\star}\Z$ by $\crosmap$. We then
have the follow result:

\begin{theorem} \label{theorem:toke}
  Let $\TSS$ be a two-sided shift space which has property $\coni$
  and let $\OSS$ be the corresponding one-sided shift space defined by
  \eqref{eq:twoone}.
  Then there are surjective $*$-homomorphisms $\kappa:\DX\to C(\TSS)$
  and $\rho:\O_{\OSS}\to C(\TSS)\rtimes_{\tsh^\star}\Z$ making the
  diagram
  \begin{equation*}
    \xymatrix{{\DX}\ar[d]_ - {\eta_{\mathcal{O}}}\ar[rr]^ -
      {\kappa}&&{C(\TSS)}\ar[d]^ - {\eta_{\rtimes}}&\\
      {\O_{\OSS}}\ar[rr]_-{\rho}&&{C(\TSS)\rtimes_{\tsh^\star}\Z}}
  \end{equation*}
  commute. We furthermore have that
  \begin{equation*} 
    \kappa(1_{\cy\mu\nu})=
    \begin{cases}
      1_{\left\{x\in\TSS\mid x_{[0,|v|[}=v\right\}}&\text{if }\exists
      w\in\al^*:v=wu,\\ 
      1_{\left\{x\in\TSS\mid x_{[|v|-|u|,|v|[}=u\right\}}&\text{if
      }\exists w\in\al^*:u=wv,\\ 
      0&\text{else,}
    \end{cases}
  \end{equation*}
  for every $\mu,\nu\in\al^*$,
  and that $\rho(S_a)=\eta_{\rtimes}(1_{D_a^{\TSS}})U$ for every
  $a\in\al$ where $U$ 
  is as in Section \ref{sec:cross} and $(S_a)_{a\in\al}$ and
  $\cy\mu\nu$ are as in Section \ref{sec:onesided}.

  If $\TSS$ also has property $(**)$, then the kernel of $\rho$ is
  isomorphic $\K^{\nTSS}$ where $\K$ is the $\cs$-algebra of compact
  operator on an infinite dimensional separable Hilbert space and
  $\nTSS$ is the number of right shift tail equivalence classes of $\TSS$
  containing a left special element.
\end{theorem}

\begin{proof} 
  Let $(\TSS,\theta_\TSS,\fg{\al})$ be the discrete partial dynamical
  system associated to $(\TSS,\tsh)$ as done in Theorem \ref{kryds}
  and $(\OSS,\theta_{\OSS},\fg{\al})$ be the discrete partial
  dynamical system associated to $(\OSS,\osh)$ as done in Theorem
  \ref{theorem:one}. We let $(D_g^{\TSS})_{g\in\fg\al}$ and
  $(\theta_g^{\TSS})_{g\in\fg\al}$ be the domains and partial
  one-to-one maps of $\theta_\TSS$, and $(D_g^{\OSS})_{g\in\fg\al}$
  and $(\theta_g^{\OSS})_{g\in\fg\al}$ the domains and partial
  one-to-one maps of $\theta_{\OSS}$.

  Let for every $A\in \Bx(\OSS,\theta_{\OSS},\fg\al)$, $\psi(A)$ be
  the subset of $\TSS$ defined by
  \begin{equation*}
    \psi(A)=\Bigl\{z\in \TSS\bigm| \forall (k,l)\in\ind\exists x\in
    A:x_{[0,k[}=z_{[0,k[}\land
    \Past_l(x_{[k,\infty[})=\{z_{[k-l,k[}\}\Bigr\}. 
  \end{equation*}
  We claim that the map $\psi$ is a Boolean homomorphism from
  $\Bx(\OSS,\theta_{\OSS},\fg\al)$ to $\Bx(\TSS,\theta_\TSS,\fg\al)$
  which satisfies that
  \begin{equation*}
    \psi\bigl(\theta_g^{\OSS}(A)\bigr)=\theta_g^{\TSS}\bigl(\psi(A)\bigr)
  \end{equation*}
  for all $A\in \Bx(\OSS,\theta_{\OSS},\fg\al)$ and $g\in \fg\al$.

  We will prove that by establishing a sequence of claims.

  \begin{claim}
    The map $\psi$ is a Boolean homomorphism.
  \end{claim}

  \begin{proof}
    We will prove that $\psi$ is a Boolean homomorphism by showing
    that $\psi(A\cap B)=\psi(A)\cap \psi(B)$ and $\psi(\OSS\setminus
    A)=\TSS\setminus \psi(A)$ for all
    $A,B\in\Bx(\OSS,\theta_{\OSS},\fg\al)$.
  
    Let $A,B\in \Bx(\OSS,\theta_{\OSS},\fg\al)$. It is obvious that
    $\psi(A\cap B)\subseteq \psi(A)\cap \psi(B)$.  Assume that $z\in
    \psi(A)\cap \psi(B)$. There exist by Lemma \ref{lemma:clopen}
    $(k_a,l_a),(k_b,l_b)\in \indd$ such that the two implications
    \begin{equation*}
      x\in A\Rightarrow \ec{k_a}{x}{l_a}\subseteq A
    \end{equation*}
    and
    \begin{equation*}
      x\in B\Rightarrow \ec{k_b}{x}{l_b}\subseteq B,
    \end{equation*}
    hold for all $x\in\OSS$.
  
    Let $(k_0,l_0)\in\indd$ and choose $(k,l)\in\indd$ such that we
    have that
    \begin{equation*}
      (k,l)\gee (k_a,l_a),(k_b,l_b), (k_0,l_0),
    \end{equation*}
    and choose $x^A\in A$ and $x^B\in B$ such that we have that
    \begin{equation*}
      z_{[0,k[}=x^A_{[0,k[}=x^B_{[0,k[} \text{ and }\Past_l(x^A_{[k,\infty[})=
      \Past_l(x^B_{[k,\infty[})=\{z_{[k-l,k[}\}.
    \end{equation*}
    Then $x^A\eq x^B$, so $x^A\in A\cap B$, and since
    $x^A_{[0,k_0[}=z_{[0,k_0[}$ and
    $\Past_{l_0}((x^A_{[k_0,\infty[})=\{z_{[k_0-l_0,k[}\}$, this shows
    that $z\in \psi(A\cap B)$. Thus $\psi(A\cap B)=\psi(A)\cap
    \psi(B)$.
  
    Assume now that $z\in \psi(\OSS\setminus A)$. We then have that if
    $x\in\OSS$, $z_{[0,k_a[}=x_{[0,k_a[}$ and
    $\Past_{l_a}(x_{[k_a,\infty[}) = \{z_{[k_a-l_a,k_a[}\}$, then
    $x\in\OSS\setminus A$. Thus $z\in\TSS\setminus\psi(A)$.
  
    If $z\in\TSS\setminus\psi(A)$, then there is a $(k_z,l_z)\in\indd$
    such that we for every $x\in A$ have that
    \begin{equation*}
      z_{[0,k_z[}\ne x_{[0,k_z[}\text{ or
      }\Past_{l_z}(x_{[k_z,\infty[})\ne \{z_{[k_z-l_z,k_z[}\}. 
    \end{equation*}
    Let $(k_0,l_0)\in\indd$ and choose $(k,l)\gee
    (k_z,l_z),(k_0.l_0)$. Since $\TSS$ has property $\coni$, there
    exists an $x\in \OSS$ such that $z_{[0,k[}=x_{[0,k[}$ and
    $\Past_l(x_{[k,\infty[}) = \{z_{[k-l,k[}\}$. This $x$ must belong
    to $\OSS\setminus A$, and since $z_{[0,k_0[}=x_{[0,k_0[}$ and
    $\Past_{l_0}(x_{[k_0,\infty[}) = \{z_{[k_0-l_0,k_0[}\}$, this
    shows that $z\in \psi(\OSS\setminus A)$. Thus $\psi(\OSS\setminus
    A)=\TSS\setminus\psi(A)$.
  \end{proof}

  \begin{claim} \label{et} We have that
    $\psi\bigl(\theta_a^{\OSS}(A)\bigr)=\theta_a^{\TSS}\bigl(\psi(A)\bigr)$
    for every $A\in \Bx(\OSS,\theta_{\OSS},\fg\al)$ and every
    $a\in\al$.
  \end{claim}

  \begin{proof}
    Let $z\in \psi\bigl(\theta_a^{\OSS}(A)\bigr)$ and let
    $(k,l)\in\indd$. Then there is an $x\in\theta_a^{\OSS}(A)$ such
    that $x_{[0,k+1[}=z_{[0,k+1[}$ and
    $\Past_l(x_{[k+1,\infty[})=\{z_{[k+1-l,k+1[}\}$. That means that
    there is a $y\in A\cap D_{a\inv}^{\OSS}$ such that
    $x=\theta_a^{\OSS}(y)$. We then have that $ z_0=x_0=a$, that
    $y_{[0,k[}=x_{[1,k+1[}=z_{[1,k+1[}$ and that
    \begin{equation*}
      \Past_l(y_{[k,\infty[})=\Past_l(x_{[k+1,\infty[})=\{z_{[k+1-l,k+1[}\}.
    \end{equation*}
    Thus $\tsh(z)\in \psi(A)$, and $z_0=a$, which shows that $z\in
    \theta_a^{\TSS}\bigl(\psi(A)\bigr)$.
  
    Now let $z\in \theta_a^{\OSS}\bigl(\psi(A)\bigr)$. Then there is a
    $z'\in \psi(A)$ such that $z=\theta_a^{\TSS}(z')$. That means that
    there for every $(k,l)\in\indd$, is an $x\in A$ such that
    \begin{equation*}
      z'_{[0,k[}= x_{[0,k[}\text{ and }\Past_l(x_{[k,\infty[})=
      \{z'_{[k-l,k[}\}. 
    \end{equation*}
    We especially have that $\Past_1(x)=\{z'_{-1}\}=\{a\}$, so $x\in
    D_{a\inv}^{\OSS}$, and we furthermore have that
    \begin{equation*}
      z_{[0,k+1[}=az'_{[0,k[}=ax_{[0,k[}=\bigl(\theta_a^{\OSS}(x)\bigr)_{[0,k+1[}
    \end{equation*}
    and that
    \begin{equation*}
      \Past_l\bigl((\theta_a^{\OSS}(x)\bigr)_{[k+1,\infty[})=
      \Past_l(x_{[k,\infty[})= \{z'_{[k-l,k[}\}= \{z_{[k+1-l,k+1[}\},
    \end{equation*}
    which shows that $z\in\psi(\theta_a^{\OSS}(A))$.
  \end{proof}
  
  \begin{claim} \label{to} We have that
    $\psi\bigl(\theta_{a\inv}^{\OSS}(A)\bigr)=
    \theta_{a\inv}^{\TSS}\bigl(\psi(A)\bigr)$ for every $a\in\al$ and
    every $A\in \Bx(\OSS,\theta_{\OSS},\fg\al)$.
  \end{claim}
  
  \begin{proof}
    Let $z\in \psi\bigl(\theta_{a\inv}^{\OSS}(A)\bigr)$ and let
    $(k,l)\in\indd$. Then there is an $x\in\theta_{a\inv}^{\OSS}(A)$
    such that $x_{[0,k+1[}=z_{[0,k+1[}$ and
    $\Past_{l+k+2}(x_{[k+1,\infty[})=\{z_{[-l-1,k+1[}\}$. That means
    that there is a $y\in A\cap D_{a}^{\OSS}$ such that
    $x=\theta_{a\inv}^{\OSS}(y)$, and we then have that
    \begin{equation*}
      a=y_0\in\Past_1(x)=\{z_{-1}\},
    \end{equation*}
    and thus that
    \begin{equation*}
      y_{[0,k[}=ax_{[0,k-1[}=z_{[-1,k-1[}\text{ and }
      \Past_l(y_{[k,\infty[})=\Past_l(x_{[k-1,\infty[})=\{z_{[k-1-l,k-1[}\}
    \end{equation*}
    if $k>0$, and
    \begin{equation*}
      \Past_l(y_{[k,\infty[})=\Past_l(ax)=\{z_{[k-1-l,k-1[}\}
    \end{equation*}
    if $k=0$.  Thus $\tsh\inv(z)\in \psi(A)$ and $z_{-1}=a$, which shows
    that
    \begin{equation*}
      z=\theta_{a\inv}^{\TSS}\bigl(\tsh(z)\bigr)
      \in\theta_{a\inv}^{\TSS}\bigl(\psi(A)\bigr).
    \end{equation*}
  
    Now let $z\in \theta_{a\inv}^{\TSS}\bigl(\psi(A)\bigr)$. We then
    have that there is a $z'\in \psi(A)\cap D_a^{\TSS}$ such that
    $z=\theta_{a\inv}^{\TSS}(z')$. That means that there for every
    $(k,l)\in\indd$ is an $x\in A$ such that 
    $z'_{[0,k+1[}=x_{[0,k+1[}$ and 
    $\Past_l(x_{[k+1,\infty[})=\{z'_{[k+1-l,k+1[}\}$. 
    We especially have that $x_0=z'_0=a$, so
    $x\in A\cap D_a$, and we also have that
    \begin{equation*}
      z_{[0,k[}=z'_{[1,k+1[}=x_{[1,k+1[}=\bigl(\theta_{a\inv}^{\OSS}(x)\bigr)_{[0,k[}
    \end{equation*}
    and that
    \begin{equation*}
      \Past_l\Bigl(\bigl(\theta_{a\inv}^{\OSS}(x)\bigr)_{[k,\infty[}\Bigr)=
      \Past_l(x_{[k+1,\infty[})= 
      \{z'_{[k+1-l,k+1[}\}=\{z_{[k-l,l[}\}.
    \end{equation*}
    This shows that $z\in\psi\bigl(\theta_{a\inv}^{\OSS}(A)\bigr)$.
  \end{proof}
  
  It follows from Claim \ref{et} and Claim \ref{to} and the definition
  of $\theta_g^{\OSS}$ and $\theta_g^{\TSS}$ that
  \begin{equation*}
    \psi\bigl(\theta_g^{\OSS}(A)\bigr)=\theta_g^{\TSS}\bigl(\psi(A)\bigr)
  \end{equation*}
  for every $A\in\Bx(\OSS,\theta_{\OSS},\fg\al)$ and every
  $g\in\fg\al$.

  Thus if we let $(s_g^{\OSS})_{g\in\fg\al}$ denote the generators of
  $\csp{\OSS,\theta_{\OSS},\fg\al}$ and $(s_g^{\TSS})_{g\in\fg\al}$
  the generators of $\csp{\TSS,\theta_{\TSS},\fg\al}$, then it follows
  from Theorem \ref{theorem:surmap} that there exists a surjective
  $*$-homomorphism $\phi$ from $\csp{\OSS,\theta_{\OSS},\fg\al}$ to
  $\csp{\TSS,\theta_\TSS,\fg\al}$ which maps $s_g^{\OSS}$ to
  $s_g^{\TSS}$ for every $g\in\fg\al$, and thus from Theorem
  \ref{kryds} and \ref{theorem:one} that there is a surjective
  $*$-homomorphism $\rho:\O_{\OSS}\to C(\TSS)\rtimes_{\tsh^\star}\Z$
  which maps $\sum_{a\in\al}S_a$ to $U$ and $\ioss(1_{D_g})$ to
  $\eta_{\rtimes}(1_{D_g})$ for every $g\in\fg\al$. We then have that
  \begin{equation*}
    \rho(S_a)=\rho\left(\ioss\biggl(1_{D_a}\sum_{a'\in\al}S_{a'}\biggr)\right)=
    \eta_{\rtimes}(1_{D_a})U
  \end{equation*}
  for every $a\in\al$, and that
  \begin{equation*}
    \begin{split}
      \rho\bigl(\ioss(1_{\cy\mu\nu})\bigr)&=
      \rho\bigl(\ioss(1_{D_\nu})\bigr)
      \rho\bigl(\ioss(1_{D_{\nu\mu\inv}})\bigr)\\
      &= \eta_{\rtimes}(1_{D_\nu})\eta_{\rtimes}(1_{D_{\nu\mu\inv}})\\
      &= \begin{cases} \eta_{\rtimes}\Bigl(1_{\left\{x\in\TSS\mid
            x_{[0,|v|[}=v\right\}}\Bigr)&
        \text{if }\exists w\in\al^*:v=wu,\\
        \eta_{\rtimes}\Bigl(1_{\left\{x\in\TSS\mid
            x_{[|v|-|u|,|v|[}=u\right\}}\Bigr)& \text{if }\exists w\in\al^*:u=wv,\\
        0&\text{else,}
      \end{cases}
    \end{split}
  \end{equation*}
  for $u,v\in\al^*$, according to Lemma \ref{lemma:onefund} and
  \ref{lemma:twofund}.  Since $\Dx$ is generated by
  $\{1_{\cy\mu\nu}\mid u,v\in\al^*\}$ and $C(X)$ is generated by the
  set
  \begin{equation*}
    \Bigl\{1_{\{x\in \TSS\mid x_{[k-|u|,k[}=u\}}\bigm| u\in\al^*,\
    k\in\Z\Bigr\},
  \end{equation*}
  it follows that there exists a surjective $*$-homomorphism
  $\kappa:\DX\to C(\TSS)$ with the desired properties.

  Assume now that $\TSS$ has property $(**)$. Let $\rte$ be the set of
  right shift tail equivalence classes of $\TSS$ which contains a left
  special element. We will show that the kernel of $\phi$, and thus
  the kernel of $\rho$, is isomorphic to $\nTSS$ copies of $\K$ by
  constructing a family
  \begin{equation*}
    (e^\jj_{x,y})_{\bigl(\jj\in\rte,\ x,y\in\spro(\jj)\bigr)}
  \end{equation*}
  of non-zero elements of $\csp{\OSS,\theta_{\OSS},\fg\al}$ such that
  the equations
  \begin{equation*}
    (e^\jj_{x,y})^*=e^\jj_{y,x}\text{ and }
    e^\jj_{x,y}e^{\jj'}_{x',y'}=
    \begin{cases}
      e^\jj_{x,y'}&\text{if }\jj=\jj'\text{ and }y=x',\\
      0&\text{else,}
    \end{cases}
  \end{equation*}
  hold for all $\jj,\jj'\in\rte$ and all $x,y\in\spro(\jj)$ and
  $x',y'\in\spro(\jj')$, and such that we have that
  \begin{equation*}
    \spc\bigr\{e^\jj_{x,y}\bigm|\jj\in\rte,\ x,y\in\spro(\jj)\bigr\}=\ker\phi.
  \end{equation*}

  So let $\jj\in\rte$ and $x,y\in\spro(\jj)$. Then there exist
  $n,m\in\N$ such that $x_{[n,\infty[}=y_{[m,\infty[}$. It follows
  from Lemma \ref{lemma:iso} that $\{x\}$ and $\{y\}$ belong to
  $\Bx(\OSS,\theta_{\OSS},\fg\al)$ so we can define an element
  $e^\jj_{x,y}$ in $\csp{\OSS,\theta_{\OSS},\fg\al}$ by the equation
  \begin{equation*}
    e^\jj_{x,y}= \boomap[(\OSS,\theta_{\OSS},\fg\al)]\bigl(\{x\}\bigr)\
    s_{x_{[0,n[}y\inv_{[0,m[}}^{\OSS}\
    \boomap[(\OSS,\theta_{\OSS},\fg\al)] \bigl(\{y\}\bigr). 
  \end{equation*}
  Notice that $e^\jj_{x,y}$ does not depend on the choice of $n$ and
  $m$, because if $k,l\in\N$ also satisfy that
  $x_{[k,\infty[}=y_{[l,\infty[}$, then we have that
  \begin{equation*}
    x_{n+l+i}=y_{m+l+i}=x_{m+k+i}
  \end{equation*}
  for all $i\in\N$, and since $x$ is not eventually periodic according
  to Lemma \ref{lemma:nonperiodic} that means that $n+l=m+k$ and thus
  that
  \begin{equation*}
    x_{[0,n[}y\inv_{[0,m[}= x_{[0,n+l[}y\inv_{[0,m+l[}=
    x_{[0,m+k[}y\inv_{[0,m+l[}= x_{[0,k[}y\inv_{[0,l[}. 
  \end{equation*}

  It is clear that $(e^\jj_{x,y})^*=e^\jj_{y,x}$ and that
  $e^\jj_{x,y}e^{\jj'}_{x',y'}=0$ if $y\ne x'$. Assume now that
  $\jj\in\rte$, that $x,y,z\in\spro(\jj)$ and that $k,l,m,n\in\N$ are
  such that $x_{[n,\infty[}=y_{[m,\infty[}$ and
  $y_{[k,\infty[}=z_{[l,\infty[}$. We then have that
  \begin{equation*}
    y\in D_{\bigl(x_{[0,n[}y\inv_{[0,m[}\bigr)\inv}^{\OSS}\text{ and }
    \theta_{x_{[0,n[}y\inv_{[0,m[}}^{\OSS}(y)=x.
  \end{equation*} 
  Thus it follows from Lemma \ref{booinvariant} that we have that
  \begin{equation*}
    \begin{split}
      s_{x_{[0,n[}y\inv_{[0,m[}}^{\OSS}\
      \boomap[(\OSS,\theta_{\OSS},\fg\al)]\bigl(\{y\}\bigr)&=
      s_{x_{[0,n[}y\inv_{[0,m[}}^{\OSS}
      \Bigl(s^{\OSS}_{x_{[0,n[}y\inv_{[0,m[}}\Bigr)^*
      s_{x_{[0,n[}y\inv_{[0,m[}}^{\OSS}\
      \boomap[(\OSS,\theta_{\OSS},\fg\al)]\bigl(\{y\}\bigr)\\
      &= s_{x_{[0,n[}y\inv_{[0,m[}}^{\OSS}\
      \boomap[(\OSS,\theta_{\OSS},\fg\al)]\bigl(\{y\}\bigr)\
      \Bigl(s_{x_{[0,n[}y\inv_{[0,m[}}^{\OSS}\Bigr)^*
      s_{x_{[0,n[}y\inv_{[0,m[}}^{\OSS}\\
      &=\boomap[(\OSS,\theta_{\OSS},\fg\al)]\bigl(\{x\}\bigr)\
      s_{x_{[0,n[}y\inv_{[0,m[}}^{\OSS}.
    \end{split}
  \end{equation*}
  We also have that $x\in D_{x_{[0,n[}y\inv_{[0,m[}}^{\OSS}$, and thus
  that
  \begin{equation*}
    \begin{split}
      &\boomap[(\OSS,\theta_{\OSS},\fg\al)]\bigl(\{x\}\bigr)\
      s_{x_{[0,n[}y\inv_{[0,m[}}^{\OSS} s_{y_{[0,k[}z\inv_{[0,l[}}^{\OSS} \\
      &\qquad =\boomap[(\OSS,\theta_{\OSS},\fg\al)]\bigl(\{x\}\bigr)\
      \boomap[(\OSS,\theta_{\OSS},\fg\al)]
      \Bigl(D_{x_{[0,n[}y\inv_{[0,m[}}^{\OSS}\Bigr)\
      s_{x_{[0,n[}y\inv_{[0,m[}y_{[0,k[}z\inv_{[0,l[}}^{\OSS}\\
      &\qquad = \boomap[(\OSS,\theta_{\OSS},\fg\al)]\bigl(\{x\}\bigr)\
      s_{x_{[0,n+k[ z\inv_{[0,m+l[}}}^{\OSS}
    \end{split}
  \end{equation*}
  according to (\ref{eq:oxc}) and (\ref{eq:oxd}). Thus we have that
  \begin{equation*}
    \begin{split}
      e^\jj_{x,y}e^\jj_{y,z}&=
      \begin{multlined}[t]
        \boomap[(\OSS,\theta_{\OSS},\fg\al)]\bigl(\{x\}\bigr)\
        s_{x_{[0,n[}y\inv_{[0,m[}}^{\OSS}\
        \boomap[(\OSS,\theta_{\OSS},\fg\al)]\bigl(\{y\}\bigr)\\
        \boomap[(\OSS,\theta_{\OSS},\fg\al)]\bigl(\{y\}\bigr)\
        s_{y_{[0,k[}z\inv_{[0,l[}}^{\OSS}\
        \boomap[(\OSS,\theta_{\OSS},\fg\al)]\bigl(\{z\}\bigr)
      \end{multlined}\\
      &= \boomap[(\OSS,\theta_{\OSS},\fg\al)]\bigl(\{x\}\bigr)\
      s_{x_{[0,n[}y\inv_{[0,m[}}^{\OSS}
      s_{y_{[0,k[}z\inv_{[0,l[}}^{\OSS}\
      \boomap[(\OSS,\theta_{\OSS},\fg\al)]\bigl(\{z\}\bigr)\\
      &= \boomap[(\OSS,\theta_{\OSS},\fg\al)]\bigl(\{x\}\bigr)\
      s_{x_{[0,n+k[}z\inv_{[0,m+l[}}^{\OSS}\
      \boomap[(\OSS,\theta_{\OSS},\fg\al)]\bigl(\{z\}\bigr)\\
      &=e^\jj_{x,z}.
    \end{split}
  \end{equation*}
  Finally we notice that it follows from Corollary \ref{boomapinj}
  that $e^\jj_{x,y}$ is non-zero, because we have that
  \begin{equation*}
    e^\jj_{x,y}\bigl(e^\jj_{x,y}\bigr)^*=e^\jj_{x,x}=
    \boomap[(\OSS,\theta_{\OSS},\fg\al)]\bigl(\{x\}\bigr). 
  \end{equation*}

  It follows from Theorem \ref{theorem:surmap} that the kernel of
  $\phi$ is generated by the set
  \begin{equation*}
    \{\boomap[\OSS,\theta_{\OSS},\fg\al](A)\mid \psi(A)=\emptyset\},
  \end{equation*}
  and since
  $\psi\Bigl(\boomap[(\OSS,\theta_{\OSS},\fg\al)]\bigl(\{x\}\bigr)\Bigr)=
  \emptyset$ if $x\in\spro(\jj)$ for some $\jj\in\rte$, we have that
  the set
  \begin{equation*}
    \spc\bigl\{e^\jj_{x,y}\mid
    \jj\in\rte,\ x,y\in\spro(\jj)\bigr\}
  \end{equation*}
  is contained in the kernel of $\phi$.

  Assume that $A\in\Bx(\OSS,\theta_{\OSS},\fg\al)$ and that
  $\psi(A)=\emptyset$. Choose by Lemma \ref{lemma:clopen}
  $(k,l)\in\indd$ such that the implication
  \begin{equation*}
    x\in A\Rightarrow \ec{k}{x}{l}\subseteq A
  \end{equation*}
  holds for all $x\in\OSS$. Then we have that
  \begin{equation*}
    A=\bigcup_{x\in A}\ec{k}{x}{l}.
  \end{equation*}
  Assume that $x\in A$ and that the number of elements of
  $\Past_l(x_{[k,\infty[})$ is one. Then it follows from Lemma
  \ref{lemma:twofund} that
  \begin{equation*}
    \begin{split}
      \psi(\ec{k}{x}{l}) &=
      \psi\left(\theta_{x_{[0,k[}}^{\OSS}\Biggl(D_{\nu\inv}^{\OSS}\cap
        \biggl(
        \smashoperator[r]{\bigcap_{\mu\in\neg\Past_l(\osh^k(x))}}
        \OSS\setminus D_{\mu\inv}^{\OSS}\biggr)\Biggr)\right)\\
      &= \theta_{x_{[0,k[}}^{\TSS}\left(D_{\nu\inv}^{\TSS}\cap \Biggl(
        \smashoperator[r]{\bigcap_{\mu\in\neg\Past_l(\osh^k(x))}}
        \OSS\setminus D_{\mu\inv}^{\TSS}\Biggr)\right)\\
      &= \theta_{x_{[0,k[}}^{\TSS}\bigl(D_{\nu\inv}^{\TSS}\bigr)\ne
      \emptyset,
    \end{split}
  \end{equation*}
  where $\nu$ is the unique element of $\Past_l(x_{[k,\infty[})$ and
  $\neg\Past_l(x_{[k,\infty[})= \bigl\{\mu\in\laOSS\mid |\mu|=l,\
  \mu\notin\Past_l(x_{[k,\infty[})\bigr\}$, but this contradict our
  assumption that $\psi(A)=\emptyset$.

  Thus $\Past_l(x_{[k,\infty[})$ consists of at least two elements for
  every $x\in A$, so it follows from Lemma \ref{lemma:faa} that $A$ is
  finite and that every element of $A$ is of the form $\spro(z)$ for
  some $z\in\TSS$ which is right shift tail equivalent to a left
  special element. Hence we get by Lemma \ref{lemma:iso} that
  \begin{multline*}
    \boomap[(\OSS,\theta_{\OSS},\fg\al)](A)= 
    \sum_{x\in A}\boomap[(\OSS,\theta_{\OSS},\fg\al)](\{x\})=\\ 
    \sum_{x\in A}e^\xx_{x,x}\in \spc\bigl\{e^\jj_{x,y}\bigm|\jj\in\rte,\
    x,y\in\spro(\jj)\bigr\},
  \end{multline*}
  where for every $x\in A$, $\xx$ denotes the right shift tail
  equivalence class which contains an element $z$ such that
  $\spro(z)=x$.

  So to show that the equality
  \begin{equation*}
    \spc\bigl\{e^\jj_{x,y}\bigm|\jj\in\rte,\ x,y\in\spro(\jj)\bigr\}=\ker\phi,
  \end{equation*}
  holds, we just have to show that the subset
  \begin{equation*}
    \spc\bigl\{e^\jj_{x,y}\bigm|\jj\in\rte,\ x,y\in\spro(\jj)\bigr\}
  \end{equation*}
  is an ideal of $\csp{\OSS,\theta_{\OSS},\fg\al}$, and since
  $\csp{\OSS,\theta_{\OSS},\fg\al}$ is generated by the set
  $\bigl\{s_g^{\OSS}\bigm| g\in\fg\al\bigr\}$, it is enough to prove
  that $s_g^{\OSS}e^\jj_{x,y}$ and $e^\jj_{x,y}s_g^{\OSS}$ belong to
  $\spc\bigl\{e^\jj_{x',y'}\bigm|\jj\in\rte,\
  x',y'\in\spro(\jj)\bigr\}$ for every $g\in\fg\al$, $\jj\in\rte$ and
  $x,y\in\spro(\jj)$.

  So let $g\in\fg\al$, $\jj\in\rte$ and $x,y\in\spro(\jj)$, and let
  $n,m\in\N$ such that $x_{[n,\infty[}=y_{[m,\infty[}$. It follows
  from Lemma \ref{lemma:onestruc} that if $s_g^{\OSS}\ne 0$, then
  there exist $\mu,\nu\in\al^*$ such that $g=\mu\nu\inv$. If $y$ does
  not belong to $D_{\mu\nu\inv}^{\OSS}$ then we have that
  \begin{equation*}
    \boomap[(\OSS,\theta_{\OSS},\fg\al)]\bigl(\{y\}\bigr)\ s_{\mu\nu\inv}^{\OSS} =
    \boomap[(\OSS,\theta_{\OSS},\fg\al)]\bigl(\{y\}\bigr)\
    \boomap[(\OSS,\theta_{\OSS},\fg\al)]\bigl(D_{\mu\nu\inv}^{\OSS}\bigr)\
    s_{\mu\nu\inv}^{\OSS}
    =\emptyset,
  \end{equation*}
  from which it follows that
  \begin{equation*}
    e^\jj_{x,y} s_g^{\OSS}= \boomap[(\OSS,\theta_{\OSS},\fg\al)]\bigl(\{x\}\bigr)\
    s_{x_{[0,n[}y\inv_{[0,m[}}^{\OSS}\
    \boomap[(\OSS,\theta_{\OSS},\fg\al)]\bigl(\{y\}\bigr)\ 
    s_{\mu\nu\inv}^{\OSS} =0. 
  \end{equation*}
  If $y\in D_{\mu\nu\inv}^{\OSS}$, then we have according to Lemma
  \ref{lemma:onefund} and \ref{booinvariant} that
  $\theta_{\nu\mu\inv}^{\OSS}(y)\in\spro(\jj)$,
  $\theta_{\nu\mu\inv}^{\OSS}(y)_{[m+|\nu|,\infty[}=
  x_{[n+|\mu|,\infty[}$ and that
  \begin{equation*}
    \begin{split}
      \boomap[(\OSS,\theta_{\OSS},\fg\al)]\bigl(\{y\}\bigr)\
      s_{\mu\nu\inv}^{\OSS} &=
      \boomap[(\OSS,\theta_{\OSS},\fg\al)]\bigl(\{y\}\bigr)\
      s_{\mu\nu\inv}^{\OSS}
      \bigl(s_{\mu\nu\inv}^{\OSS}\bigr)^*s_{\mu\nu\inv}^{\OSS}\\
      &= s_{\mu\nu\inv}^{\OSS} \bigl(s_{\mu\nu\inv}^{\OSS}\bigr)^*
      \boomap[(\OSS,\theta_{\OSS},\fg\al)]\bigl(\{y\}\bigr)\
      s_{\mu\nu\inv}^{\OSS}\\
      &= s_{\mu\nu\inv}^{\OSS}\ \boomap[(\OSS,\theta_{\OSS},\fg\al)]
      \bigl(\{\theta_{\nu\mu\inv}^{\OSS}(y)\}\bigr).
    \end{split}
  \end{equation*}
  We further more have that
  \begin{equation*}
    \begin{split}
      y\inv_{[0,m+|\mu|[}uv\inv &= \bigl(uy_{[|\mu|,m+|\mu|[}\bigr)\inv uv\inv\\
      &= y_{[|\mu|,m+|\mu|[}\inv v\inv\\
      &= \left(\theta_{vu\inv}^{\OSS}(y)_{[0,m+|\nu|[}\right)\inv,
    \end{split}
  \end{equation*}
  and hence that
  \begin{multline*}
    s_{x_{[0,n+|\mu|[}y\inv_{[0,m+|\mu|[}}^{\OSS} s_{\mu\nu\inv}^{\OSS} =\\
    \boomap[(\OSS,\theta_{\OSS},\fg\al)]
    \Bigl(D_{x_{[0,n+|\mu|[}y\inv_{[0,m+|\mu|[}}^{\OSS}\Bigr)
    s_{x_{[0,n+|\mu|[}\bigl(\theta_{vu\inv}^{\OSS}(y)_{[0,m+|\nu|[}\bigr)\inv}^{\OSS}
  \end{multline*}
  according to \eqref{eq:oxc} and \eqref{eq:oxd} of Definition
  \ref{def:ox}.

  Thus we have that
  \begin{equation*}
    \begin{split}
      e^\jj_{x,y} s_g^{\OSS} &=
      \boomap[(\OSS,\theta_{\OSS},\fg\al)]\bigl(\{x\}\bigr)\
      s_{x_{[0,n[}y\inv_{[0,m[}}^{\OSS}\
      \boomap[(\OSS,\theta_{\OSS},\fg\al)]\bigl(\{y\}\bigr)\
      s_{\mu\nu\inv}^{\OSS}\\
      &= \boomap[(\OSS,\theta_{\OSS},\fg\al)]\bigl(\{x\}\bigr)\
      s_{x_{[0,n[}y\inv_{[0,m[}}^{\OSS} s_{\mu\nu\inv}^{\OSS}\
      \boomap[(\OSS,\theta_{\OSS},\fg\al)]
      \bigl(\{\theta_{\nu\mu\inv}^{\OSS}(y)\}\bigr)
      \\
      &= \begin{multlined}[t]
        \boomap[(\OSS,\theta_{\OSS},\fg\al)]\bigl(\{x\}\bigr)\
        \boomap[(\OSS,\theta_{\OSS},\fg\al)]\
        \bigl(D_{x_{[0,n+|\mu|[}y\inv_{[0,m+|\mu|[}}^{\OSS}\bigr)
        \\
        s_{x_{[0,n+|\mu|[}\left(\theta_{vu\inv}(y)_{[0,m+|\nu|[}\right)\inv}^{\OSS}\
        \boomap[(\OSS,\theta_{\OSS},\fg\al)]
        \bigl(\{\theta_{\nu\mu\inv}^{\OSS}(y)\}\bigr)
      \end{multlined}\\
      &= \begin{multlined}[t]
        \boomap[(\OSS,\theta_{\OSS},\fg\al)]\bigl(\{x\}\bigr)\
        s_{x_{[0,n+|\mu|[}(\theta_{vu\inv}(y)_{[0,m+|\nu|[})\inv}^{\OSS}\\
        \boomap[(\OSS,\theta_{\OSS},\fg\al)]
        \bigl(\{\theta_{\nu\mu\inv}^{\OSS}(y)\}\bigr)
      \end{multlined}\\
      &= e^\jj_{x,\theta_{\nu\mu\inv}^{\OSS}(y)}.
    \end{split}
  \end{equation*}
  So we have in all cases that $e^\jj_{x,y} s_g^{\OSS}\in
  \spc\bigl\{e^\jj_{x',y'}\bigm|\jj\in\rte,\
  x',y'\in\spro(\jj)\bigr\}$.

  One can in a similar way prove that we have that
  \begin{equation*}
    s_g^{\OSS} e^\jj_{x,y}\in
    \spc\bigl\{e^\jj_{x',y'}\mid\jj\in\rte,\
    x',y'\in\spro(\jj)\bigr\}. 
  \end{equation*}

  Hence $\spc\bigl\{e^\jj_{x,y}\mid\jj\in\rte,\
  x,y\in\spro(\jj)\bigr\}$ is an ideal of
  $\csp{\OSS,\theta_{\OSS},\fg\al}$ and thus equal to the kernel of
  $\phi$.
\end{proof}

\begin{example} \label{exam:sub}
  Let $\TSS_\eta$ be the shift space of an aperiodic proper
  substitution $\eta$. It then follows from \cite{MR924156}*{page 90
    and 107} and \cite{MR1355295}*{Theorem 3.9} that $\TSS_\eta$ is
  minimal (and thus also $\TSS^+_\eta$) and contains a finite, but
  nonzero, number of left special elements. Thus according to
  \cite{tmcseii}*{Example 3.6}, $\TSS_\eta$ 
  has property $(**)$, and since ${\nTSS}_\eta$ is nonzero, it follows
  from Theorem \ref{theorem:toke} that $\O_{\TSS^+_\eta}$ is not
  simple. On the other hand $C(\TSS_\eta)\rtimes_{\tsh^\star}\Z$ is
  simple because $\TSS_\eta$ is minimal.
\end{example}

\appendix
\section{Partial representations of groups} \label{sec:parrep}
A \emph{partial representation} of a group $G$ was defined in
\cite{MR1469405} to be a map $u$ from $G$ to a $\cs$-algebra $\Aa$
such that the following conditions hold:
\begin{subequations}    
  \begin{align}
    &u(e)=1,\label{eq:prxa}\\
    &u(g\inv)=u(g)^* \text{ for every }g\in G\label{eq:prxb},\\
    &u(h)u(i)u(i\inv)=u(hi)u(i\inv) \text{ for every }h,i\in G.\label{eq:prxc}
  \end{align}
\end{subequations}
In \cite{MR99a:46121} another definition of a partial representation
of a group was given, namely as a map $u$ from $G$ to a $\cs$-algebra
$\Aa$ such that the following conditions hold: 
\begin{subequations}
  \begin{align}
    &\parbox[t]{0.8\textwidth}{$\bigl(u(g)\bigr)_{g\in G}$ is a family of partial
      isometries with commuting\\ range projections,}\label{eq:prqra} \\ 
    &u(e)u(e)^*=1,\label{eq:prqrb}\\
    &u(g)^*u(g)=u(g\inv)u(g\inv)^* \text{ for every } g\in G,\label{eq:prqrc}\\
    &u(h)u(i)u(i)^*u(h)^*=u(h)u(i)u(hi)^* \text{ for every } h,i\in
    G.\label{eq:prqrd} 
  \end{align}
\end{subequations}
It was in \cite{MR99a:46121}*{Lemma 1.8 and Remark 1.9} also noticed
that the conditions 
\eqref{eq:prqra}--\eqref{eq:prqrd} are equivalent to the conditions:
\begin{subequations}
  \begin{align}
    &\parbox[t]{0.8\textwidth}{$\bigl(u(g)\bigr)_{g\in G}$ is a family of partial
      isometries with commuting\\ range projections,}\label{eq:pra} \\ 
    &u(e)=1,\label{eq:prb}\\
    &u(g)^*=u(g\inv) \text{ for every } g\in G,\label{eq:prc}\\
    &u(h)u(i)=u(h)u(h)^*u(hi) \text{ for every } h,i\in
    G.\label{eq:prd} 
  \end{align}
\end{subequations}
As mentioned in \cite{MR2003f:46108}, the conditions
\eqref{eq:prxa}--\eqref{eq:prxc} are equivalent to the conditions
\eqref{eq:prqra}--\eqref{eq:prqrd} and thus to the conditions
\eqref{eq:pra}--\eqref{eq:prd}, but since we have not been able to
find a complete proof of this, we will now give one:
\begin{proposition}
  Let $u$ be a map from a group $G$ to a $\cs$-algebra $\Aa$. Then the
  following are equivalent:
  \begin{enumerate}
  \item $u$ satisfies the conditions \eqref{eq:prxa}--\eqref{eq:prxc},
  \item $u$ satisfies the conditions \eqref{eq:prqra}--\eqref{eq:prqrd},
  \item $u$ satisfies the conditions \eqref{eq:pra}--\eqref{eq:prd}.
  \end{enumerate}
\end{proposition}

\begin{proof}
  As mentioned above, the equivalence of $(2)$ and $(3)$ follows from
  \cite{MR99a:46121}*{Lemma 1.8 and Remark 1.9}. 

  Assume that  $u$ satisfies the conditions
  \eqref{eq:prxa}--\eqref{eq:prxc}. Then we have that
  \begin{equation*}
    u(g)u(g)^*u(g)=u(g)u(g\inv)u(g)=u(e)u(g)=u(g),
  \end{equation*}
  so $u(g)$ is a partial isometry for all $g\in G$.

  If $h,i\in G$, then we have that
  \begin{equation*}
    \begin{split}
      u(h)u(h)^*u(i)u(i)^*&=u(h)u(h\inv i)u(i)^*\\
      &=u(i)u(i\inv h)u(h\inv i)u(i)^*\\
      &=\bigl(u(h\inv i)u(i\inv)\bigr)^*\bigl(u(h\inv i)u(i\inv)\bigr)\\
      &=u(i)u(i)^*u(h)u(h)^*u(i)u(i)^*,
    \end{split}
  \end{equation*}
  from which it follows that
  \begin{equation*}
    \begin{split}
      u(i)u(i)^*u(h)u(h)^*&=\bigl(u(h)u(h)^*u(i)u(i)^*\bigr)^*\\
      &=\bigl(u(i)u(i)^*u(h)u(h)^*u(i)u(i)^*\bigr)^*\\
      &=u(i)u(i)^*u(h)u(h)^*u(i)u(i)^*\\
      &=u(h)u(h)^*u(i)u(i)^*.
    \end{split}
  \end{equation*}
  This shows that $u$ satisfies condition
  \eqref{eq:prqra}.

  It is clear that $u$ satisfies condition
  \eqref{eq:prqrb} and \eqref{eq:prqrc}. If $h,i\in G$, then we have that
  \begin{equation*}
    \begin{split}
      u(h)u(i)u(i)^*u(h)^*&=u(h)\bigl(u(h)u(i)u(i)^*\bigr)^*\\
      &=u(h)\bigl(u(hi)u(i)^*\bigr)^*\\
      &=u(h)u(i)u(hi)^*,
    \end{split}
  \end{equation*}
  so $u$ also satisfies condition
  \eqref{eq:prqrd}. Thus $(1)$ implies $(2)$.

  Assume now that $u$ satisfies the conditions
  \eqref{eq:pra}--\eqref{eq:prd}, and that $h,i\in G$. Then we have that
  \begin{equation*}
    \begin{split}
      u(h)u(i)u(i)^*&=\bigl(u(i)u(i)^*u(h\inv)\bigr)^*\\
      &=\bigl(u(i)u(i\inv h\inv)\bigr)^*\\
      &=u(hi)u(i\inv),
    \end{split}
  \end{equation*}
  which shows that $(3)$ implies $(1)$.
\end{proof}
 
\section{Boolean Algebras} \label{sec:boo}
We recommend \cite{MR29:4713} for a very nice introduction to Boolean
algebras. A \emph{Boolean algebra} is a set $\Bx$ with two
distinct elements $0,1\in\Bx$ which act like the empty set and the
whole set, respectively, and three
operations $\lor:\Bx\times\Bx\to\Bx$, $\land:\Bx\times\Bx\to\Bx$ and
$\neg:\Bx\to\Bx$ which act like union, intersection and complement,
respectively. To be precise, they satisfy the following axioms:
\setlength\arraycolsep{2pt}
\begin{eqnarray}
  \neg 0=1&& \neg 1=0\\
  A\land 0=0&& A\lor 1=1\\
  A\land 1=A&& A\lor 0=0\label{eq:1}\\
  A\land\neg A=0&& A\lor\neg A=1\label{eq:2}\\
  \neg\neg A&=&A\\
  A\land A=A&&A\lor A=A\\
  \neg(A\land B)=\neg A\lor\neg B&& \neg(A\lor B)=\neg A\land\neg B\\
  A\land B=B\land A&& A\lor B=B\lor A\label{eq:3}\\
  A\land(B\land C)=(A\land B)\land C&& A\lor(B\lor C)=(A\lor B)\lor C\\
  A\land(B\lor C)&=&(A\land B)\lor(A\land C)\label{eq:4}\\A\lor(B\land
  C)&=&(A\lor B)\land(A\lor C)\label{eq:5} 
\end{eqnarray}
This set of axioms is not the shortest one possible. In fact one could
do with for example just axiom (\ref{eq:1}), (\ref{eq:2}),
(\ref{eq:3}), (\ref{eq:4}) and (\ref{eq:5}). We will call
$A\lor B$ for the \emph{intersection of $A$ and $B$}, $A\land B$ for the
\emph{union of $A$ and $B$}, and $\neg A$ the \emph{complement of $A$}.

The generic example of a Boolean algebra is of course the power set of
a set $X$, where $0=\emptyset$, $1=X$, $\lor=\cup$, $\land=\cap$ and
$\neg A=X\setminus A$. 

If $\Aa$ is a subset of a Boolean algebra $\Bx$, then we call it a
\emph{Boolean subalgebra} of $\Bx$ if $0,1\in\Aa$ and $A\lor B, A\land B,
\neg A\in \Aa$ for every $A,B\in \Aa$. In this case $\Aa$ is of course
itself a Boolean algebra with operations inherited from $\Bx$. When
$X$ is a set, then we will by a \emph{Boolean algebra on $X$} mean a Boolean
subalgebra of the power set of $X$. 

An example of this which we will use in this paper, is if $X$ is a
topological space. Then the set of clopen subsets of $X$ is a Boolean
subalgebra of the power set of $X$ and thus a Boolean algebra on $X$. 

If $\Aa$ is some subset of a Boolean
algebra $\Bx$, then we will by the \emph{Boolean algebra generated by
  $\Aa$} mean the
smallest Boolean subalgebra of $\Bx$ containing $\Aa$. Notice that if a
subset $\Aa$ of a Boolean algebra $\Bx$ is closed under intersection
(respectively union) and complement, then it also closed under union
(respectively intersection) and thus is a Boolean subalgebra.

Other examples of Boolean algebras
which we will meet in this paper is the set $\{0,1\}$ where 
\begin{eqnarray*}
  0\lor 0=0\land0=0\land 1=1\land 0=\neg 1=0,\\
  1\land 1=1\lor 1=1\lor 0=0\lor 1=\neg 0=1,
\end{eqnarray*}
and the set of projections in a unital abelian $\cs$-algebra, where
\begin{align*}
  p\lor q=pq&& p\land q=p+q-pq&& \neg p=1-p.
\end{align*}

A map $\phi$ between two Boolean algebras $\Bx$ and $\Bx'$ is called 
\emph{a Boolean homomorphism} if $\phi(A\lor B)=\phi(A)\lor\phi(B)$,
$\phi(A\land B)=\phi(A)\land\phi(B)$ and $\phi(\neg A)=\neg\phi(A)$ for
every $A,B\in \Bx$. In fact, the first (respectively the second)
together with the last equality 
imply the second (respectively the first), so in order to verify
that $\phi$ is a Boolean homomorphism, it is enough to check these two
equalities. Notice that when $\phi$ is a Boolean homomorphism, then
$\phi(1)=1$ and $\phi(0)=0$.

If $\Aa$ is a unital abelian $\cs$-algebra and $\Bx$ is the Boolean algebra of
projections of $\Aa$, then it is easy to check that $\spa(\Bx)$ is a
$*$-subalgebra of $\Aa$ and thus that $\spc(\Bx)$ is a
$\cs$-subalgebra of $\Aa$. Hence $\spc{\Bx}=\csp{\Bx}$.

\begin{lemma} \label{lemma:homo}
  If for $i\in\{1,2\}$, $\Aa_i$ is a unital abelian $\cs$-algebra 
  and $\Bx_i$ is the Boolean algebra of projections of $\Aa_i$, and
  $\phi:\Bx_1\to\Bx_2$ is a Boolean homomorphism, then there is a
  uniquely determined $*$-homomorphism from $\spc{\Bx_1}$ to
  $\spc{\Bx_2}$ which maps $A$ to $\phi(A)$ for $A\in\Bx_1$.
\end{lemma}

\begin{proof}
  Since $\Bx_1$ generates $\spc{\Bx_1}$, there can at most be one
  $*$-homomorphism from $\spc{\Bx_1}$ to 
  $\spc{\Bx_2}$ which maps $A$ to $\phi(A)$ for $A\in\Bx_1$.

  Let $\cset$ denote the set of finite subset of $\Bx_1$. We then have
  that $\spc{\Bx_1}$ is the closure of the set
  \begin{equation*}
    \bigcup_{C\in\cset}\spa{C},
  \end{equation*}
  and since $\spa{C}$ is the $\cs$-subalgebra of $\Aa_1$ generated by
  $C$, it is enough to show that there for every $C\in \cset$ exists a
  $*$-homomorphism from $\spa{C}$ to
  $\spc{\Bx_2}$ which maps $A$ to $\phi(A)$ for $A\in C$.

  So let $C\in\cset$. There then exists a finite family
  $p_1,p_2,\dotsc, p_n$ of mutually orthogonal projections in $C$ such
  that every element of $\spa{C}$ uniquely can be written as
  \begin{equation*}
    \sum_{i_1}^n\lambda_ip_i
  \end{equation*}
  with $\lambda_1,\lambda_2,\dotsc,\lambda_n\in \C$. The map 
  \begin{equation*}
    \sum_{i_1}^n\lambda_ip_i\mapsto \sum_{i_1}^n\lambda_i\phi(p_i)
  \end{equation*}
  is therefore a well-defined $*$-homomorphism from $\spa{C}$ to
  $\spc{\Bx_2}$ which maps $A$ to $\phi(A)$ for $A\in C$.
\end{proof}

\section{Crossed products of $\cs$-partial dynamical systems}
\label{sec:parcro} 

A \emph{$\cs$-partial dynamical system} has been defined in 
\cite{MR99a:46121} to be a triple $(A,G,\alpha)$ where $A$ is a
$\cs$-algebra, $G$ is a discrete group and $\alpha$ is a 
\emph{partial action} of $G$ on $A$. That means that $\alpha$ consists
of a family $(D_g)_{g\in G}$ of closed ideals of $A$ and a family
$(\alpha_g)_{g\in G}$ of isomorphisms $\alpha_g:D_{g\inv}\to D_g$ such that 
\begin{align}
  &D_e=A,\\
  &\parbox[t]{0.8\textwidth}{$\alpha_{hi}$ extends $\alpha_h\alpha_i$
    for all $h,i\in G$ (where the domain of $\alpha_h\alpha_i$ is
    $\alpha_i\inv(D_{h\inv})$).}
\end{align}
$\cs$-partial dynamical systems have been studied in
\cites{MR2003f:46108,MR99a:46121,MR1331978} (the definition of $\cs$-partial
dynamical systems in \cite{MR1331978} is a bit different from the
above mentioned, 
but it is showed in \cite{MR99a:46121}*{Remark 1.9} that the two
definitions are equivalent). 

A \emph{covariant representation} of a $\cs$-partial dynamical system
$(A,G,\alpha)$ on a Hilbert space $\H$ is a pair $(\pi,u)$ where $\pi$
is a non-degenerate representation of $A$ on $\H$, and $u$ is a partial
representation (cf. Appendix \ref{sec:parrep}) of $G$ on $\H$ such
that for each $g\in G$, $u(g)u(g)^*$ is the projection of $\H$ onto
the subspace $\spc\pi(D_g)\H$, and
$\pi\bigl(\alpha_g(a)\bigr)=u(g)\pi(a)u(g\inv)$ for $a\in D_g$. As it
is the case 
with $\cs$-partial dynamical systems, the definition of a covariant
representation in \cite{MR1331978} is a bit different from the above
mentioned, but 
it is shown in \cite{MR99a:46121}*{Remark 1.12} that the two
definitions are equivalent.

The \emph{crossed product} $A\rtimes_\alpha G$ of a $\cs$-partial
dynamical system 
$(A,G,\alpha)$ is a $\cs$-algebra which is generated by a copy of $A$
and a family $(\delta_g)_{g\in G}$ of elements such that there exists
a bijective map $(\pi,u)\mapsto \pi\times u$ 
between covariant representations of $(A,G,\alpha)$ on $\H$ and
non-degenerated representations of $A\rtimes_\alpha G$ on $\H$ such
that $(\pi\times u)(a\delta_g)=\pi(a)u(g)$ for $g\in G$ and $a\in D_g$.

\end{document}